\documentclass[3p,times]{elsarticle}

\usepackage{array}
\usepackage{color}
\usepackage{tabularx}
\usepackage{graphicx}
\usepackage{amsmath}
\usepackage{amssymb}
\usepackage{amsfonts}	
\usepackage{moreverb}
\usepackage{dsfont}
\usepackage{bm}
\usepackage{multirow}
\usepackage{soul}
\usepackage{microtype} 

\newcommand\BibTeX{{\rmfamily B\kern-.05em \textsc{i\kern-.025em b}\kern-.08em
T\kern-.1667em\lower.7ex\hbox{E}\kern-.125emX}}

\newcommand{\x}{\mbf{x}}

\newcommand{\mbf}[1]{\mathbf{#1}}			%

\newcommand{\Q}{\mathbf{Q}}

\newcommand{\w}{\mathbf{w}}
\newcommand{\q}{\mathbf{q}}
\newcommand{\F}{\mathbf{F}}
\newcommand{\f}{\mathbf{f}}
\newcommand{\g}{\mathbf{g}}

\renewcommand{\v}{\mathbf{v}}
\newcommand{\B}{\mathbf{B}}

\newcommand{\be}{\begin{equation}}
\newcommand{\ee}{\end{equation}}
\newcommand{\bdm}{\begin{displaymath}}
\newcommand{\edm}{\end{displaymath}}

\newcommand{\apriori}{\textit{a priori }}

\renewcommand{\phi}{\varphi}
\usepackage{epstopdf}  %???
\usepackage{subfig}
\newcommand{\de}[2]{\frac {\partial #1}{\partial#2}}
\newcommand{\talg}[1]{\texttt{#1}}

\newcommand{\Nl}{N_{\text{edge}}^{\,n}}
\newcommand{\Nn}{N_{\text{node}}^{\,n}}
\renewcommand{\S}{\mathbf{S}}
\newcommand{\urho}{u_r}
\newcommand{\uphi}{u_\phi}
\newcommand{\sphi}{\sin \phi \,}
\newcommand{\cphi}{\cos \phi \,}

\newcommand{\nVertex}{N\mathsf{v}(i)}

\newcommand{\Vn}{\mathcal{V}}

\journal{Computers and Fluids}

\newfont{\numerikEleven}{ecrm1000}
\newfont{\numerikTen}{cmss10}
\newfont{\numerikNine}{cmss9}
\newfont{\numerikEight}{cmss8}

\begin{document} 
%!=========================================================================
%!
%!      F R O N T    M A T T E R 
%!
\begin{frontmatter}
%-------------------------------------------------------
% TITLE
\title{Direct Arbitrary-Lagrangian-Eulerian finite volume schemes on moving nonconforming unstructured meshes} 
%-------------------------------------------------------
%-------------------------------------------------------
% AUTHORS
\author[UniTNMath]{Elena Gaburro}
\ead{elena.gaburro@unitn.it}

\author[UniTN]{Michael Dumbser$^{*}$}
\ead{michael.dumbser@unitn.it}
\cortext[cor1]{Corresponding author}

\author[Uma]{Manuel J. Castro}
\ead{castro@anamat.cie.uma.es}

%-------------------------------------------------------
% INSTITUTIONS
\address[UniTNMath]{Department of Mathematics, University of Trento, Via Sommarive, 14 - 38123 Trento, Italy}  
\address[UniTN]{Department of Civil, Environmental and Mechanical Engineering, University of Trento, Via Mesiano, 77 - 38123 Trento, Italy.}
\address[Uma]{Department of Mathematical Analysis, Statistics and Applied Mathematics, University of M\'alaga, Campus de Teatinos s/n, 29071, M\'alaga, Spain.}
%-------------------------------------------------------

%-------------------------------------------------------
% ABSTRACT
\begin{abstract}
	In this paper, we present a novel second-order accurate Arbitrary-Lagrangian-Eulerian (ALE) finite volume scheme on moving 
	\textit{nonconforming} polygonal grids, in order to avoid the typical mesh distortion caused by \textit{shear flows} 
	in Lagrangian-type methods.
	In our new approach the nonconforming element interfaces are \textit{not} defined by the user, 
	but they are automatically detected by the algorithm  if the tangential velocity difference across an element interface 
	is sufficiently large. 
	The grid nodes that are sufficiently far away from a shear wave are moved with a standard node solver,  
	while at the interface we insert a new set of nodes that can \textit{slide along} the interface in a nonconforming
	manner. In this way, the elements on both sides of the shear wave can move with a different velocity, without producing
	highly distorted elements. 
	
	The core of the proposed method is the use of a space-time conservation formulation in the construction of the final finite volume 
	scheme, which completely avoids the need of an additional remapping stage, hence the new method is a so-called \textit{direct} ALE 
	scheme. For this purpose, the governing PDE system is rewritten at the aid of the space-time divergence operator and then a fully 
	discrete
	one-step discretization is obtained by integrating over a set of closed space-time control volumes. The nonconforming sliding of 
	nodes along an edge requires the insertion or the deletion  of nodes and edges, and in particular the space-time faces of 
	an element can be shared between more than two cells.    
	
	Due to the space-time conservation formulation, the geometric conservation law (GCL) is \textit{automatically} 
	satisfied by \textit{construction}, even on moving nonconforming meshes. Moreover, the mesh quality remains high and, as a direct 
	consequence, also the time step remains almost constant in time, even 
	for highly sheared vortex flows. In this paper we focus mainly on logically straight slip-line interfaces, but we show also 
	first results for general slide lines that are \textit{not} logically straight.  
	Second order of accuracy in space and time is obtained by using a MUSCL-Hancock strategy, together with a Barth and Jespersen 
	slope limiter. 

The accuracy of the new scheme has been further improved by incorporating a special \textit{well balancing} technique that is 
able to maintain particular stationary solutions of the governing PDE system up to machine precision. In particular, we consider 
steady vortex solutions of the shallow water equations, where the pressure gradient is in equilibrium with the centrifugal forces.
	
	A large set of different numerical tests has been carried out in order to check the accuracy and the robustness of the new method for both smooth and discontinuous problems. In particular we have compared the results for a steady vortex in equilibrium solved with a standard \textit{conforming} ALE method (without any rezoning technique) and with our new \textit{nonconforming} ALE scheme, to show that the new nonconforming scheme is able to avoid mesh distortion in vortex flows even after very long simulation times.

\end{abstract}
%-------------------------------------------------------

%-------------------------------------------------------
% KEY WORDS
%-------------------------------------------------------
\end{frontmatter}
%!=========================================================================

\section{Introduction} 
\label{sec.intro}

Arbitrary-Lagrangian-Eulerian (ALE) finite volume schemes are characterized by a moving computational mesh: at each time step the new position of all the nodes has to be recomputed according to a  prescribed mesh velocity, which generally is chosen as close as possible to the local fluid velocity (as it is in the purely Lagrangian framework), but it can also be set to zero (to reproduce the 
Eulerian case), or it can be chosen arbitrarily. The aim of these methods is to reduce the numerical dissipation errors due to the convective terms, hence to capture contact 
discontinuities sharply and to precisely identify and track material interfaces. For these reasons already in the 1950's John von Neumann and Richtmyer were working on Lagrangian schemes 
\cite{Neumann1950} for one-dimensional flows, and Wilkins proposed a two-dimensional extension in 1964, see \cite{Wilkins1964}. 
Since the fluid velocity is required at each node and at each time step, a natural approach is a \textit{staggered} discretization, where the momentum is defined at the grid vertices and 
all the other flow variables are defined at the cell center.  
Despite some drawbacks of the initial version of staggered Lagrangian schemes, which was not conservative and which produced some spurious modes in the numerical solution, it was widely used in 
the last forty years and many improvements have been made in the meantime; for further details one can refer to the papers of Caramana and Shashkov \cite{Caramana1998, CaramanaShashkov1998} 
and of Loub\`ere at al. \cite{LoubereMaire2010, LoubereMaire2013}.  Moreover, examples on general polygonal grids have been presented 
in \cite{LoubereShashkov2004}.  

An alternative consists in a \textit{conservative cell-centered} discretization, which was first introduced by Godunov in \cite{godunov}. An early application of conservative 
cell-centered Godunov-type schemes to the compressible Euler equations of gas dynamics in a Lagrangian framework on moving grids was provided by Munz in \cite{munz94}, using Roe-type and 
HLL-type approximate Riemann solvers. 
In many recent papers, see for example Despr\'es et al. \cite{DepresMazeran2003,Despres2005,Despres2009} and Maire et al. \cite{Maire2007, Maire2009, Maire2011, Maire2010},    
the conservative cell-centered Godunov-type approach is used both with structured and unstructured moving grids and in two and three space dimensions, respectively.  
Successively also better than second-order accurate schemes were introduced: \textit{high order} of accuracy in \textit{space} was first achieved by Cheng and Shu \cite{chengshu1,chengshu2} 
by means of a non linear ENO reconstruction, and high order in \textit{time} was obtained either by the use of Runge-Kutta type methods or by adapting the ADER-WENO schemes to the Lagrangian 
framework, see for example Dumbser et al. \cite{Lagrange1D} and Cheng \& Toro \cite{chengtoro}. Recent work on high order Lagrangian discontinuous Galerkin finite element methods can be 
found in the papers of Vilar et al. \cite{Vilar1,Vilar2,Vilar3}, Yu et al. \cite{Yuetal} and Boscheri \& Dumbser \cite{LagrangeDG}, while high order Lagrangian continuous finite elements 
have been studied by Scovazzi et al. \cite{scovazzi1,scovazzi2} and Dobrev and Rieben et al. \cite{Dobrev1,Dobrev2,Dobrev3}. 

For all the cell-centered methods an important step is the computation of the fluid velocity at the nodes, since this information is not directly available in the scheme, but it has to be 
extrapolated from the adjacent cells. To obtain these values three different types of \textit{node solvers} can be employed. 
The simplest one is that proposed in the above mentioned papers by Cheng and Shu \cite{chengshu1,chengshu2}, somehow employed also in this work, where the node velocity is obtained as arithmetic average 
among the near states; another possibility, suggested by Despr\'es et al. (GLACE scheme) \cite{Despres2009} and Maire (EUCCLHYD scheme) \cite{Maire2009}, is to solve multiple one-dimensional 
half-Riemann problems around a vertex, in order to get an approximate solution of the multi-dimensional (generalized) Riemann problem at the node; the most recent method introduced by 
Balsara et al. \cite{Balsara2012,BalsaraMultiDRS,MUSIC1,MUSIC2} consists in solving approximately a multidimensional Riemann problem at the nodes, using 
a new family of genuinely multidimensional HLL-type Riemann solvers. They are all compared with each other within a high-order ADER-WENO ALE scheme in the recent paper of Boscheri et al. 
\cite{LagrangeMHD}. 

Although all these different schemes are widely used, especially to describe compressible multi-material flows, a common problem that affects all Lagrangian methods is the severe mesh 
\textit{distortion} or the 
mesh tangling that happens in the presence of shear flows and that may even destroy the computation. Hence, all Lagrangian methods must be in general combined with an algorithm to (locally) rezone 
the mesh at least from time to time and to remap the solution from the old mesh to the new mesh in a conservative manner. 
Lagrangian remesh and remap ALE schemes are very popular and some recent work on that topic can be found in \cite{MaireMM2, ShashkovRemap1, MaireMM3,  ShashkovRemap3, ShashkovMultiMat2}. 
Extensions of the remesh-remap approach to better than second-order of accuracy can be found in \cite{Loubere2014,Blanchard2016}. 
In contrast to indirect ALE schemes (purely Lagrangian phase, remesh and subsequent remap phase) there are the so-called direct ALE schemes, where the local rezoning is performed before the 
computation of the numerical fluxes, hence changing directly the chosen mesh velocity of the ALE approach, see for example \cite{LagrangeMDRS,Lagrange2D,Lagrange3D} for recent work in that 
direction based on high order ADER-WENO schemes.  

Moreover, when dealing with shear flows at material interfaces in realistic applications, see e.g. \cite{Kuchariv14Impact}, a special treatment of \textit{slide lines} may be required.  
The introduction of slide lines goes back to an idea of Wilkins \cite{Wilkins1964}, successively studied and refined by Caramana \cite{Caramana2009}, Barlow et al. \cite{Barlow2000} and Loub\`ere at al. \cite{LoubereSL2013}; the main ideas adopted in their papers regard the subdivision of the nodes at the interface in master and slave nodes and the study of the forces between the two sides of the 
slide lines. Another very interesting approach to slide lines was presented by Clair et al. in \cite{Clair201356,Clair2014315} and by Del Pino et al. in \cite{DelPino2010,Bertoluzza2016}. 
In \cite{Caramana1998} a staggered Lagrangian code has been presented, where the internal interfaces are handled with a special type of boundary condition. A very original solution to the problem of shear flows in Lagrangian simulations has been recently proposed by Springel
in \cite{springel}, where the connectivity of the moving mesh is dynamically regenerated via a moving unstructured but conforming Voronoi  tesselation of the domain.     

This paper presents a novel second-order accurate \textit{direct} cell-centered ALE scheme on unstructured polygonal grids, which deals with shear flows in an original and effective way. 
The sliding element interfaces are automatically \textit{detected} during the computation (not fixed \apriori), and nodes along such sliding edges are allowed to move in a \textit{nonconforming} 
way by the insertion and deletion of new nodes and new edges. The algorithm can handle both the insertion and the deletion of nodes and edges, using completely \textit{local} procedures. 
In the straight slip-line case no distinction  between master and slave nodes is required and the mesh movement is only based on the corner-extrapolated values of the fluid velocity. 
%, eventually projected onto the nearest edge to avoid the creation of holes or the superposition of elements. (troppo per l'abstract)

At this point we also would like to refer to some recent works on high order Eulerian and ALE schemes on moving meshes with time-accurate \textit{local time stepping} (LTS) presented in 
\cite{ALELTS1D,ALELTS2D,LTSTransport}, where each element is allowed to run at its own optimal local time step according to a \textit{local} CFL stability condition. 
The resulting algorithms use a \textit{conforming} grid in space, but naturally produce a \textit{nonconforming} mesh in time. Therefore, the new nonconforming ALE method 
presented in this paper, which produces a nonconforming mesh in both space and time, is naturally related to some of the ideas forwarded in \cite{ALELTS1D,ALELTS2D} in the 
context of local time stepping. 

Finally, we incorporate in our scheme a new well balancing technique that is useful for preserving particular steady state solutions of the governing PDE. In this paper, we consider 
an equilibrium vortex, where the pressure gradient is exactly balanced by the centrifugal forces. 

The rest of the paper is organized as follows. At the beginning of Section~\ref{sec.method} we outline the main stages of the proposed numerical method, then we describe the employed second order 
finite volume scheme emphasizing the novelties due to the nonconforming mesh motion. In particular in Section~\ref{subsec.motion} we explain how to deal with the moving nonconforming (hanging) 
nodes and the corresponding \textit{local} update of the connectivity tables of the unstructured mesh. In Section~\ref{sec.results} some numerical test problems are presented in order to check 
the efficiency and the robustness of the proposed approach in maintaining a high quality mesh, local and global volume conservation, and in satisfying the GCL condition. 
The algorithm presented here is not necessarily limited to logically straight slipe lines. In Section \ref{sec.general} we therefore show first preliminary results for general, logically non-straight
slide lines. Finally, in Section \ref{sec.WB}, we introduce a so-called \textit{well balancing} technique, which allows our scheme to be more accurate near  the steady states solutions. These methods,  presented in details in particular in \cite{Pares2006}, are adapted to the shallow water equations  in polar coordinates and to the context of a nonconforming moving meshes. 
The paper is closed by some conclusions and an outlook to future work.  

%----------------------------------------------
\section{Numerical method}
\label{sec.method}
We consider two-dimensional non linear hyperbolic systems of conservation laws that can be cast in the following general form
\be
\label{eq.generalform}
\de{\Q}{t} + \nabla \cdot \F(\Q) = \mathbf{S}(\Q), \quad \x \in \Omega(t) \subset \mathbb{R}^2, \Q \in \Omega_{\Q} \subset \mathbb{R}^{\nu},
\ee
where $\x = (x,y)$ is the spatial position vector, $t$ represents the time, $\Q = (q_1,q_2, \dots, q_{\nu})$ is the vector of conserved variables defined in the space of the admissible states $\Omega_{\Q} \subset \mathbb{R}^{\nu}$, $ \F(\Q) = (\,\f(\Q), \g(\Q)\,) $ is the non linear flux tensor, and $\mathbf{S}(\Q)$ represents a non linear algebraic source term.

To discretize the moving domain, we consider a total number $N_E$ of polygonal elements $T_i^n$ (the \textit{spatial control volumes}), each one with an arbitrary number of vertices $\nVertex$, $i=1, \dots, N_E$: the union of all these elements results in an unstructured mesh $\mathcal{T}_\Omega^n$ which covers the computational domain $\Omega(\x, t^n) = \Omega^n$ at the current time $t^n$ and which
contains a total number $\Nn$ of nodes and a total number $\Nl$ of edges. 

At each time step the algorithm evolves the cell averages via a discrete form of the equation \eqref{eq.generalform} and computes the new node positions through the following intermediate stages: 
\begin{enumerate}[I.]
	\item
	First, the edges along relevant shear flows are detected and the nodes on these edges are marked as problematic. 
	\item
	Then the new node positions are computed according to the type of the considered node, in particular
	\begin{enumerate}[a)]
		\item
		regular non-hanging nodes that are not in regions of relevant shear flow (i.e. they have not been marked as problematic) are evolved using a mass-weighted Cheng and Shu node solver; 
		\item
		regular non-hanging nodes that are in regions of relevant shear flow (i.e. they have been marked as problematic) are \textit{doubled}; their new position is projected along the nearest  
		interface edge, and they subsequently change their type from regular non-hanging nodes to \textit{hanging} nodes; 
		\item
		hanging nodes on an edge are allowed to slide only along that edge, and if they get too close to other nodes, they are \textit{merged} together (deleted), eventually changing back 
		their type from hanging nodes to regular non-hanging nodes. 
	\end{enumerate}
	\item \label{stage.FVscheme}
	Finally, we apply a MUSCL-Hancock strategy to produce a space-time reconstruction of the data from the known cell averages, and to evolve the solution to the new geometry (without any remapping 
	stage) using a space-time conservation formulation of the governing PDE system. 
\end{enumerate}
While this is the natural execution order in the computer program, for the sake of clarity we first start this section by briefly summarizing our space-time finite volume scheme, focusing in particular 
on the computation of the \textit{numerical fluxes} at the \textit{nonconforming interfaces}. Only later, after having introduced the connectivity tables employed in the algorithm (Section \ref{ssec.ConnectTables}) and having described the interface detector (Section \ref{ssec.DiscDetector}), we will discuss in detail the procedure to move each kind of node and the corresponding \textit{local} update of the connectivity tables.

\subsection{Finite volume scheme: reconstruction and time-evolution}
\label{sec.Reconstruction}

As usual in a classical cell-centered finite volume scheme, at the beginning of each time step $t^n$ we dispose of  the cell averages $\Q_i^n$ of the conserved variables for each spatial control volume $T_i^n$ , defined as
\[
\Q_i^n = \frac{1}{\vert T_i^n \vert} \int_{T_i^n} \Q(\x,t^n)\,d\x, 
\] 
where $\vert T_i^n \vert$  denotes the area of $T_i^n$. These are the data computed and stored at the previous time, and which will be used to evolve the solution during the current time step. To construct a method which is better than first order accurate we cannot compute the numerical fluxes directly with these piecewise constant data, but we have to reconstruct for each $T_i^n$ a piecewise  space-time polynomial $\q_h(\x,t^n)$, exploiting the cell averages of the cell and its neighbors, combined with a time-evolution procedure. 

Here, second order of accuracy in space and time is achieved by using the MUSCL-Hancock strategy that was for the first time proposed by van Leer in \cite{leer2} and which is very well explained
in \cite{toro-book}, slightly adapted to our context of a nonconforming moving mesh. 

For the \textit{spatial} reconstruction, let us define a polynomial $\w_h(\x,t^n)$ of the form   
\[
\w_h(\x,t^n) \, \vert_{T_i^n} = \w_i(\x,t^n) = \Q_i^n + \nabla \Q_i (\x - \x_i), \quad \ \x \in T_i^n,
\]
where $\x_i$ is the barycenter of cell $T_i^n$. We denote by $\mathcal{S}_i^n$ the set of neighbors of $T_i^n$ that share a common edge with $T_i^n$ 
(the set $\mathcal{S}_i^n$ may change at each time step). 
To compute $\nabla \Q_i$, integral conservation is imposed on each element of $\mathcal{S}_i^n$ 
\be
\label{eq.IntegralCons}
\frac{1}{\vert T_j^n \vert} \int_{T_j^n} \w_h(\x,t) \, d \x = \Q_j^n \quad \ \forall T_j^n \in \mathcal{S}_i^n.
\ee
The above system is in general over-determined, so we add the constraint that equation \eqref{eq.IntegralCons} holds exactly at least for $T_i^n$. This is easily satisfied by rewriting the equations as 
\be
\label{eq.IntegralCons2}
\frac{1}{\vert T_j^n \vert} \int_{T_j^n} \nabla \Q_i (\x - \x_i) \, d \x = \Q_j^n - \Q_i^n \quad \  \forall T_j^n \in \mathcal{S}_i^n,
\ee
then we solve \eqref{eq.IntegralCons2} via a classical least-squares approach using the normal equation of \eqref{eq.IntegralCons2}, and we thus obtain the \textit{non-limited} slope $\nabla \Q_i$.
To ensure that new extrema are not created in the reconstruction process, we employ the classical \textit{slope limiter} function $\Phi_i$ presented by Barth and Jespersen in \cite{BarthJespersen}.  

%The idea is to find the largest admissible $\Phi_i$ in such a way that 
%\[
%\tilde{\w}_h(\x,t^n) = \Q_i^n + \Phi_i \nabla \Q_i (\x - \x_i) % + \partial_t \Q_i(t-t^n)
%\]
%satisfies
%\[
%\min \limits_{j \in \mathcal{V}_i^n} \Q_j^n = \Q_i^{min} \le \tilde{\w}_h(\x,t^n) \le \Q_i^{max} = \max \limits_{j \in \mathcal{V}_i^n} \Q_j^n,
%\]
%where $\Q_i^{max}$ and $\Q_i^{min}$ are the componentwise maximum and minimum among the cell-averages of the set $\mathcal{V}_i^n$, respectively.
%The set  $\mathcal{V}_i^n$ contains \textit{all} the vertex neighbors of $T_i^n$ and the element $T_i^n$ itself. 
%Since $\w_h$ is obtained as a piecewise linear reconstruction of the data, its extrema occur at the vertices of $T_i^n$. Hence, to compute the limiter for each conserved variable, it suffices to find for all vertices $\x_j$ of $T_i^n$ 
%\[
%\begin{aligned}
%& \Phi_{i,j} = 
%\begin{cases}
%\min \left (1, \frac{\Q_i^{max} - \Q_i^n }{\w_{h,j} - \Q_i^n }\right ), \qquad & \text{if}   \ \w_{h,j} - \Q_i^n > 0  \\[4pt]
%\min \left (1, \frac{\Q_i^{min} - \Q_i^n }{\w_{h,j} - \Q_i^n }\right ), \qquad & \text{if} \  \w_{h,j} - \Q_i^n < 0  \\[4pt]
%1  &\text{if} \    \w_{h,j} - \Q_i^n = 0.
%\end{cases} 
%\end{aligned}
%\]
%with $\w_{h,j}= \w_h(\x_j, t^n)$ (ratios and inequalities are to be understood component-wise). Then, the slope limiter is defined as 
%\[
%\Phi_i = \min_{j} \, ( \Phi_{i,j}).
%\]

Finally, second order of accuracy in \textit{time} is achieved by an element-local predictor stage that evolves the reconstructed polynomials ${\w}_i(\x,t^n)$ within each element $T_i^n(t)$ during the time interval $[t^n,t^{n+1}]$. The piecewise space-time polynomials are denoted by $\q_h(\x,t)$, and are of the form
\be
\q_h (\x,t) \vert_{T_i^n} = \q_i(\x,t^n) = \Q_i^n + \Phi_i \nabla \Q_i (\x - \x_i) + \partial_t \Q_i(t-t^n), \quad \x \in T_i(t), \  t \in [t^n,t^{n+1}].
\ee
The value of $\partial_t \Q_i$ can be easily computed from a discrete integral form of \eqref{eq.generalform} by summing 
over the set of edges of $T_i^n$, denoted by $\mathcal{E}_i$: %is recovered through the strong form of the PDE
\be
\label{eq.strongForm}
%\Q_t  = - \mathbf{f}_x(\Q) - \mathbf{g}_y(\Q) + \mathbf{S}(\Q),
\partial_t \Q_i = - \sum \limits_{e \in \mathcal{E}_i} \left( \lambda_e \mathbf{F}\left( {\q}_i(\x_e, t^n) \right) \cdot \mathbf{n}_e \right) + \mathbf{S}\left( {\q}_i(\x_i, t^n) \right). 
\ee
Here, $\x_e$ denotes the midpoint of edge $e$, $\lambda_e$ is the edge length and $\mathbf{n}_e$ is the outward-pointing unit normal vector; $\x_i$ is the barycenter of cell $T_i^n$. 
 
%where the r.h.s of \eqref{eq.strongForm} can be easily computed. Indeed the fields $\F$ and $\mathbf{G}$ over $T_i^n$ are approximated as linear fields
%\[
%\begin{aligned}
%&\F(x,y) = f_0 + f_1(x-x_i) + f_2(y-y_i),  \\
%&\mathbf{G}(x,y) = g_0 + g_1(x-x_i) + g_2(y-y_i),  \quad \ \text{with} \ \x_i=(x_i,y_i)  \  \text{the barycenter of} \ T_i^n, \\
%\end{aligned}
%\] whose coefficients $f_i$ and $g_i $ are determined interpolating the values of the fields computed at the vertices $j$ of $T_i^n$, i.e.
%\[
%\F(\tilde{\w}_h(\x_j,t^n)) \quad \text{ and } \quad  \mathbf{G}(\tilde{\w}_h(\x_j,t^n)), \quad \ \forall j=1:\nVertex.
%\]
%Then $f_1 = \F_x(\Q)$ and $g_2 = \mathbf{G}_y(\Q)$. Besides, the source $\mathbf{S}(\Q)$ is computed at the barycenter $\x_i$ of $T_i^n$, i.e. $S(\tilde{\w}_h(\x_i,t^n))$.

\subsection{Finite volume scheme: space-time conservation form}
\label{sec.Finite volume scheme: one step element update}
\begin{figure}[p]
	\centering
	\includegraphics[width=0.36\textwidth]{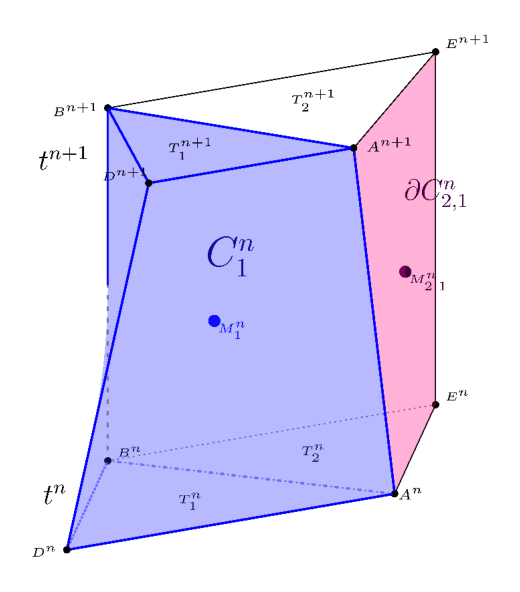}
	\includegraphics[width=0.41\textwidth]{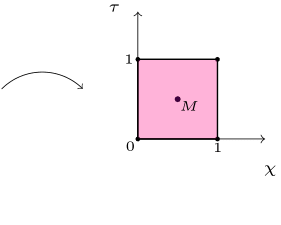}
	\caption{Left. In blue we show the physical space--time control volume $C^n_1$ obtained by connecting via straight line segments each vertex of $T_1^n$ with the corresponding vertex of $T_1^{n+1}$, and its space-time midpoint $M_1^n$. In pink we show one of the lateral surfaces of $C_2^n$, $\partial C_{2,1}^n$, together with its space--time midpoint $M_{2,1}^n$. Right. The reference system $(\chi,\tau)$ adopted for the bilinear parametrization of the lateral surfaces $\partial C^n_{ij}$.}
	\label{fig.Cin}
\end{figure}

Once $\q_h(\x,t)$ has been computed for each $T_i^n$, we are in the position to introduce the one-step space-time finite volume scheme. 
As proposed in \cite{Lagrange2D}, the governing PDE \eqref{eq.generalform} is first reformulated in a space--time divergence form as
\be
\tilde \nabla \cdot \tilde{\F} = \mathbf{S}(\Q), 
\label{PDEdiv3D}
\ee 
with
\be
\tilde \nabla  = \left( \frac{\partial}{\partial x}, \, \frac{\partial}{\partial y}, \, \frac{\partial}{\partial t} \right)^T,  \qquad 
\tilde{\F}  = \left( \mathbf{F}, \, \Q \right) = \left( \mathbf{f}, \, \mathbf{g}, \, \Q \right),
\ee
and it is then integrated over the space--time control volume $C_i^n = T_i(t) \times [t^n,t^{n+1}] $
\begin{equation}
\int_{t^{n}}^{t^{n+1}} \int_{T_i(t)} \tilde \nabla \cdot \tilde{\F} \, d\mathbf{x} dt = \int_{t^{n}}^{t^{n+1}} \int_{T_i(t)} \S \, d\mathbf{x} dt.
\label{STPDE}
\end{equation} 
The space--time control volumes $C_i^n$ are obtained by connecting each vertex of the element $T_i^n$ via \textit{straight} line segments with the corresponding vertex   
of $T_i^{n+1}$. For a graphical interpretation one can refer to Figure \ref{fig.Cin}, where we have reported an example of a control volume and of the parametrization of the 
lateral space--time surfaces. A lateral space--time  surface is denoted by $\partial C_{ij}^n$, where the index $i$ refers to the space-time element $C_i^n$ and the index $j$ 
refers to the number of the neighbor control volume $C_j^n$ that shares $\partial C_{ij}^n$ with $C_i^n$. 

\medskip 

Now, for each control volume we have to compute its barycenter and the areas, the normal vectors, and the space--time midpoints of all its sub--surfaces.
 
The upper space--time sub--surface $T_i^{n+1}$ and the lower space--time sub--surface $T_i^{n}$ are the simplest, since they are orthogonal to the time coordinate. The space--time unit normal vectors are respectively $\mathbf{\tilde n} = (0,0,1)$ and $\mathbf{\tilde n} = (0,0,-1)$. To compute their areas we can use the \textit{shoelace formula} or \textit{Gauss's area formula} which is valid for any type of polygonal element 
\be
\label{eq.shoelaceFormula}
\vert\, T_i^n\, \vert = \frac{1}{2} \left \vert  x_{\nVertex}^n y_1^n - x_1^n y_{\nVertex}^n + \sum_{j=1}^{\nVertex-1} \left( x_j^n y_{j+1}^n - x_{j+1}^n y_{j}^n \right) \right \vert,
\ee
where $\mbf{x}_j^n = (x_j^n, y_j^n), \ j = 1, \dots, \nVertex,$ are the coordinates of the vertices of element $T_i^n$ numbered in a counterclockwise order.
%The only thing that has to be taken into account is the periodicity of the domain: so if the domain is periodic w.r.t the $y$ coordinate we rewrite formula \eqref{eq.shoelaceFormula} in the form
%\[
%\vert \,T_i^n\, \vert= \frac{1}{2} \left \vert \sum_{j=1}^{\nVertex} x_j^n \left(y_{j+1}^n-y_{j-1}^n\right) \right \vert , 
%\]
%which depends only on the signed distance between the $y$ coordinates of the vertices.
Moreover, the space--time barycenter $M_i^n$ of each control volume $C_i^n$ reads
\[
M_i^n = \left ( \frac{\x_i^n + \x_i^{n+1}}{2}, \frac{t^n+t^{n+1}}{2} \right ), 
\]
where the spatial barycenter $\x_i^n=(x_i^n, y_i^n)$ of $T_i^n$ is given by the explicit formula
\begin{equation}
\begin{aligned}
\x_i^n & = \frac{1}{6 \, \vert T_i^n \vert} \sum_{j=1}^{\nVertex} \left (\x_j^n + \x_{j+1}^n \right) \left (x_j^n y_{j+1}^n - x_{j+1}^n y_{j}^n \right ), \\
%y_i^n & = \frac{1}{6 \, \vert T_i^n \vert} \sum_{j=1}^{\nVertex} \left (y_j^n + y_{j+1}^n \right ) \left (x_j^n y_{j+1}^n - x_{j+1}^n y_{j}^n %\right ).  
\end{aligned}
\end{equation}
with the convention that $j =\nVertex + 1$ coincides with $j= 1$.

Next, the lateral space--time surfaces of $C_i^n$ are parametrized using a set of bilinear basis functions  
\be
\partial C_{ij}^n = \mathbf{\tilde{x}} \left( \chi,\tau \right) = 
\sum\limits_{k=1}^{4}{\beta_k(\chi,\tau) \, \mathbf{\tilde{X}}_{ij,k}^n },	
\qquad 0 \leq \chi \leq 1,  \quad	0 \leq \tau \leq 1, 										 
\label{eq.SurfParBeta}
\ee
where the $\mathbf{\tilde{X}}_{ij,k}^n$ represent the physical space--time coordinates of the four vertices of $\partial C_{ij}^n$, and the $\beta_k(\chi,\tau)$ functions are defined as follows 
\begin{equation}
\beta_1(\chi,\tau) = (1-\chi)(1-\tau), \qquad 
\beta_2(\chi,\tau) = \chi(1-\tau), \qquad 
\beta_3(\chi,\tau) = \chi\tau, \qquad  
\beta_4(\chi,\tau) = (1-\chi)\tau.
\label{BetaBaseFunc}
\end{equation}
The mapping in time is given by the transformation 
\be
t = t_n + \tau \, \Delta t, \qquad  \tau = \frac{t - t^n}{\Delta t}, 
\label{eq.timeTransf}
\ee 
hence the Jacobian matrix $J_{\partial C_{ij}^n}$ of the parametrization is 
\be
J_{\partial C_{ij}^n} = \left( \begin{array}{ccc} \vec{e}_x & \vec{e}_y & \vec{e}_t \\ 
	\frac{\partial x}{\partial \chi} & \frac{\partial y}{\partial \chi} & \frac{\partial t}{\partial \chi} \\ 
	\frac{\partial x}{\partial \tau } & \frac{\partial y}{\partial \tau } & \frac{\partial t}{\partial \tau } 
\end{array} \right) = \left( \begin{array}{c} \mathbf{\tilde{e}} \\ \frac{\partial \mathbf{\tilde{x}}}{\partial \chi} \\ \frac{\partial \mathbf{\tilde{x}}}{\partial \tau} \end{array} \right). 
\label{JacSTsurf}
\ee
The space--time unit normal vector $\mathbf{\tilde n}_{ij}$ can be evaluated computing the normalized cross product between the transformation vectors of the mapping \eqref{eq.SurfParBeta}, i.e.
\be
| \partial C_{ij}^n| = \left| \frac{\partial \mathbf{\tilde{x}}}{\partial \chi} \times \frac{\partial \mathbf{\tilde{x}}}{\partial \tau} \right|, 
\quad 
\mathbf{\tilde n}_{ij} = \left( \frac{\partial \mathbf{\tilde{x}}}{\partial \chi} \times \frac{\partial \mathbf{\tilde{x}}}{\partial \tau}\right) / | \partial C_{ij}^n|,
\label{n_lateral}
\ee
where $| \partial C_{ij}^n|$ is the determinant of the Jacobian matrix $J_{\partial C_{ij}^n}$ and represents also the area of the lateral surfaces.
Moreover, exploiting the parametrizations in \eqref{eq.SurfParBeta}-\eqref{eq.timeTransf} and choosing 
$\chi = 0.5$ and $\tau =0.5$ we recover the coordinates $M_{i,j}^n$ of the space--time midpoint of the lateral surfaces.

\begin{figure}[p]
	\centering
	\subfloat[Insertion of a new node]{
		\includegraphics[width=0.32\linewidth, height=0.32\textheight]{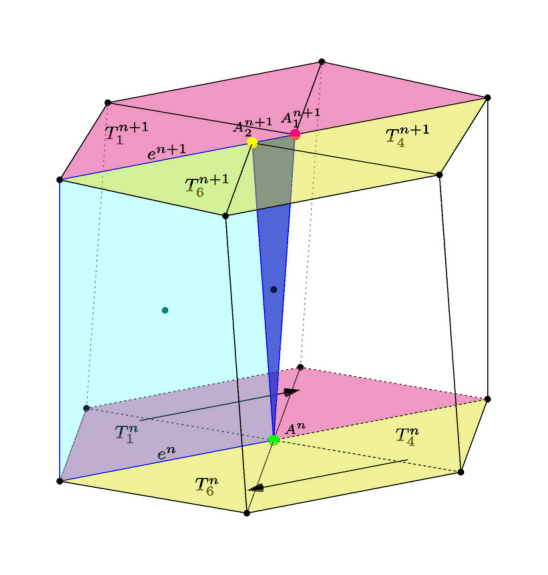}}
	\subfloat[Motion of hanging nodes]{
		\includegraphics[width=0.32\linewidth, height=0.32\textheight]{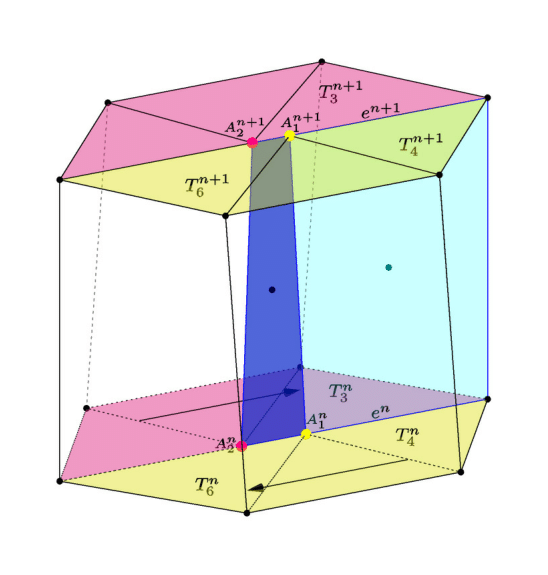}}
	\subfloat[Fusion of two existing nodes]{
		\includegraphics[width=0.32\linewidth, height=0.32\textheight]{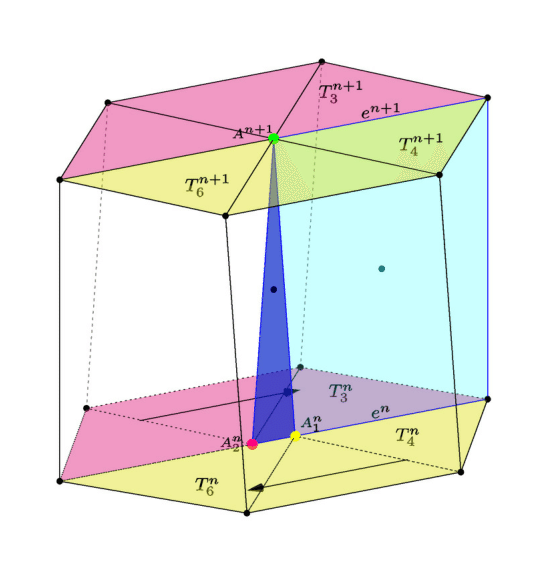}}
	\caption{Suppose that at time $t^n$ across the pink and the yellow elements the tangential fluid velocity changes sharply, as suggested by the arrows; at $t^{n+1}$ the pink elements will move in one  direction and the yellow ones will move in the opposite direction. In (a) at time $t^n$ we have a conforming mesh, but in order to avoid a severe mesh distortion in the following time steps we decide to double the green node $A^n$. So at time $t^{n+1}$ there are both $A_1^{n+1}$ and $A_2^{n+1}$: $A_1^{n+1}$ is a vertex for the pink elements and $A_2^{n+1}$ is a vertex for the yellow elements. Moreover $A_2^{n+1}$ is hung to edge $e^{n+1}$. So the blue lateral face of $T_i^n$, which has $e^n$ and $e^{n+1}$ as bases, is composed by two pieces: the one in light blue which is trapezoidal and touches elements $T_1$ and $T_6$, and the one in dark blue which is triangular and touches elements $T_1$ and $T_4$. Note in particular that we need to compute the flux between $T_1$ and $T_4$ during the interval $[t^n, t^{n+1}]$ even if at time $t^n$ they were not in contact. In (b) we show the intermediate situation where a hanging node slides along an edge. In this case the blue surface is still divided into two parts and it is shared between three elements $T_3, T_4$ and $T_6$, so the computation of two fluxes will be required. In order to compute the fluxes and to maintain the second order of accuracy of the entire method the reconstruction polynomial $\q_h(\x,t)$ will be evaluated at the midpoints of each of the part of the lateral surface. Finally, in (c) we report the last limiting case: $A_1^n$ and $A_2^n$ are close and at $t^{n+1}$ will be even closer since they are moving one towards the other, so we decide to merge them and to restore the conforming and simpler situation, in particular to avoid  that $A_1^n$ will leave edge $e^{n+1}$ at time $t^{n+1}$. Eventually $A^{n+1}$ could be doubled again at $t^{n+2}$ if the tangential velocity difference across the interface is sufficiently large.}
	\label{fig.NonConformingControlVolumes}
\end{figure}

Finally, when we allow a node to slide along an edge in a \textit{nonconforming way}, the lateral space--time surfaces have to be treated slightly differently, refer to Figure~\ref{fig.NonConformingControlVolumes} for a graphical interpretation. 

Consider the case of  $\partial C_{i,j}^n$ with the four standard vertices and two more hanging nodes on the edges orthogonal to the time coordinate (as in the middle of Figure~\ref{fig.NonConformingControlVolumes}). Then the lateral surface is shared between three (and not two, as usual) control volumes. However it can be subdivided into two pieces, each one shared between only two  control volumes, which are still trapezoidal, so each of them can be mapped into the reference element using the standard map in \eqref{eq.SurfParBeta}, just taking care to select in a correct way the vertices of each piece. Hence areas, normal vectors, and space--time midpoints can be computed exactly as in the conforming case but on each part, and so we will recover these data for each piece.
Next, on the left and on the right of Figure~\ref{fig.NonConformingControlVolumes} we have reported the two limiting cases: first, at time $t^{n+1}$ a new node has been inserted, which at the previous  time $t^n$ did not yet exist; or vice-versa, at time $t^{n+1}$ a hanging node is merged together with one of the other vertices and hence it disappears. In these cases the lateral surfaces can still be   divided into two parts, and even if one of them is triangular it can still be treated as a degenerate quadrilateral face, so all the computations can be performed, once again, as above. The coordinates of a hanging node at the moment of its creation or destruction will be set equal 
to those of the vertex from which the hanging node was born, or those of the vertex to whom it was merged, respectively. 

Note that the treatment of the nonconforming lateral space--time surfaces basically requires only to repeat the computation of the necessary geometric information over each piece and the same will hold for the flux computation, which will be simply split in several parts.

\medskip 

Once having computed all the relevant geometric information about each control volume $C_i^n$ and its space--time surface 
\begin{equation}
\partial C^n_i = \left( \bigcup \limits_{j} \partial C^n_{ij} \right) 
\,\, \cup \,\, T_i^{n} \,\, \cup \,\, T_i^{n+1}, 
\label{dCi}
\end{equation}
%and eventually of each sub--piece of $ \partial C^n_{ij}$, 
%we can go on with the space--time divergence form of the PDE 
%\be
%\int_{t^{n}}^{t^{n+1}} \int_{T_i(t)} \tilde \nabla \cdot \tilde{\F} \, d\mathbf{x} dt = \int_{t^{n}}^{t^{n+1}} \int_{T_i(t)} \S \, d\mathbf{x} %dt.
%\label{STPDE2}
%\ee 
we can apply the Gauss theorem to the integral with the space-time flux divergence in \eqref{STPDE} and obtain 
\be
\label{eq.I1}
\int_{\partial C^n_i} \tilde{\F} \cdot \ \mathbf{\tilde n} \, \, dS = 
\int_{t^{n}}^{t^{n+1}} \int_{T_i(t)} \S \, d\mathbf{x} dt,   
\ee    
where $\mathbf{\tilde n} = (\tilde n_x,\tilde n_y,\tilde n_t)$ is the outward pointing space--time unit normal vector on the space--time surface $\partial C^n_i$. Substituting the physical boundary fluxes $\tilde{\F} \cdot \ \mathbf{\tilde n}$ with appropriate numerical fluxes
leads to a consistent and conservative finite volume discretization. In principle, the entire structure of the numerical scheme is already 
given by \eqref{eq.I1}.     
The final one--step direct ALE finite volume scheme is then obtained from equation \eqref{eq.I1} as
\be
	|T_i^{n+1}| \, \Q_i^{n+1} = |T_i^n| \, \Q_i^n - \sum \limits_{j} \,\, {\int_0^1 \int_0^1 
		| \partial C_{ij}^n| \tilde{\F}_{ij} \cdot \mathbf{\tilde n}_{ij} \, d\chi d\tau}
	+ \int_{t^{n}}^{t^{n+1}} \int_{T_i(t)} \S(\mathbf{q}_h) \, d\mathbf{x} dt, 
	\label{eq.PDEfinal}
\ee
where the discontinuity of the solution at the space--time sub--face $\partial C_{ij}^n$ is resolved by an ALE numerical flux function $\tilde{\F}_{ij} \cdot \mathbf{\tilde n}_{ij}$, which computes the flux between two neighbors across the intermediate space--time lateral surface. In particular when the lateral surface is shared between more than two control volumes (as shown in Figure~\ref{fig.NonConformingControlVolumes}) we have to compute the fluxes across each sub-piece and sum each contribution. 
The results presented in this paper are obtained using a Rusanov--type or an Osher--type ALE flux. Note that in time we have used the upwind
flux due to the causality principle, which naturally leads to the terms $|T_i^{n}| \, \Q_i^{n}$ and $|T_i^{n+1}| \, \Q_i^{n+1}$.  

Let $\q_h^-(\x,t)$ be the reconstructed numerical solution inside the element $T_i(t)$ and $\q_h^+(\x,t)$ be the reconstructed numerical solution inside one of the neighbors of $T_i^n$ through $\partial C_{i,j}^n$; let $\q_h^-$ and $\q_h^+$ the values of these polynomials evaluated at the space-time midpoint $M_{i,j}^n$ of the considered piece of the lateral surface. 
Define the ALE Jacobian matrix w.r.t. the normal direction in space
\begin{equation} 
\label{eq.ALEjacobianMatrix}
\mathbf{A}^{\!\! \mathbf{V}}_{\mathbf{n}}(\Q)=\left(\sqrt{\tilde n_x^2 + \tilde n_y^2}\right)\left[\frac{\partial \mathbf{F}}{\partial \Q} \cdot \mathbf{n}  - 
(\mathbf{V} \cdot \mathbf{n}) \,  \mathbf{I}\right], \qquad    
\mathbf{n} = \frac{(\tilde n_x, \tilde n_y)^T}{\sqrt{\tilde n_x^2 + \tilde n_y^2}},  
\end{equation} 
with $\mathbf{I}$ representing the identity matrix and $\mathbf{V} \cdot \mathbf{n}$ denoting the local normal mesh velocity. 

Then the expression for the Rusanov flux is given by
\begin{equation}
\tilde{\F}_{ij} \cdot \mathbf{\tilde n}_{ij} =  
\frac{1}{2} \left( \tilde{\F}(\q_h^+) + \tilde{\F}(\q_h^-)  \right) \cdot \mathbf{\tilde n}_{ij}  - 
\frac{1}{2} s_{\max} \left( \q_h^+ - \q_h^- \right),  
\label{eq.rusanov} 
\end{equation} 
where $s_{\max}$ is the maximum eigenvalue of $\mathbf{A}^{\!\! \mathbf{V}}_{\mathbf{n}}(\q_h^+)$ and $\mathbf{A}^{\!\! \mathbf{V}}_{\mathbf{n}}(\q_h^-)$.

The Osher--type flux formulation has been proposed in the Eulerian framework in \cite{OsherUniversal} and has been subsequently extended to moving meshes in one and two space dimensions in \cite{Lagrange1D} and \cite{Lagrange2D}, respectively. It is defined as
\begin{equation}
\tilde{\F}_{ij} \cdot \mathbf{\tilde n}_{ij} =  
\frac{1}{2} \left( \tilde{\F}(\q_h^+) + \tilde{\F}(\q_h^-)  \right) \cdot \mathbf{\tilde n}_{ij}  - 
\frac{1}{2} \left( \int_0^1 \left| \mathbf{A}^{\!\! \mathbf{V}}_{\mathbf{n}}(\boldsymbol{\Psi}(s)) \right| ds \right) \left( \q_h^+ - \q_h^- \right),  
\label{eqn.osher} 
\end{equation} 
where we choose to connect the left and the right state across the discontinuity using a simple straight--line segment path  
\begin{equation}
\boldsymbol{\Psi}(s) = \q_h^- + s \left( \q_h^+ - \q_h^- \right), \qquad 0 \leq s \leq 1.  
\label{eqn.path} 
\end{equation} 
The absolute value of the dissipation matrix in \eqref{eqn.osher} is evaluated as usual as 
\begin{equation}
|\mathbf{A}| = \mathbf{R} |\boldsymbol{\Lambda}| \mathbf{R}^{-1},  \qquad |\boldsymbol{\Lambda}| = \textnormal{diag}\left( |\lambda_1|, |\lambda_2|, ..., |\lambda_\nu| \right),  
\end{equation}
where $\mathbf{R}$ and $\mathbf{R}^{-1}$ denote, respectively, the right eigenvector matrix and its inverse of the ALE Jacobian 
$\mathbf{A}^{\!\! \mathbf{V}}_{\mathbf{n}} = \frac{\partial \mathbf{F}}{\partial \Q} \cdot \mathbf{n} - (\mathbf{V}\cdot \mathbf{n}) \mathbf{I}$. 
In \eqref{eq.PDEfinal} the time step $\Delta t$ is given by 
\begin{equation}
	\Delta t = \textnormal{CFL} \, \min \limits_{T_i^n} \frac{d_i}{|\lambda_{\max,i}|}, \qquad \forall T_i^n \in \Omega^n, 
	\label{eq:timestep}
\end{equation}
where $\textnormal{CFL}$ is the Courant-Friedrichs-Levy number, $d_i$ represents the encircle diameter of element $T_i^n$ and $|\lambda_{\max,i}|$ is the maximum absolute value of the eigenvalues computed from the solution $\Q_i^n$ in $T_i^n$. As stated in \cite{toro-book} in Chapter 16, for linear stability in two space dimensions the Courant number must satisfy $\textnormal{CFL} \leq 0.5$.

We underline that the integration over a closed space--time control volume, as done above, automatically satisfies the so-called geometric conservation law (GCL),  
since from the Gauss theorem follows 
\begin{equation}
	\int_{\partial \mathcal{C}_i^n} \mathbf{\tilde n} \, dS = 0. 
	\label{eqn.gcl} 
\end{equation} 
The relation between \eqref{eqn.gcl} and the usual form of the GCL that is typically employed in the community working on Lagrangian schemes has been established in the appendix of \cite{Lagrange3D}. 
For all the numerical test problems shown later in this paper it has been explicitly verified that property \eqref{eqn.gcl} holds for all elements and for all time steps up to machine precision, even on moving nonconforming meshes. 

We would like to emphasize that the direct ALE scheme presented here does in general \textit{not} lead to a vanishing mass flux across element boundaries, similar to previous work 
on direct ALE schemes presented in \cite{Lagrange2D,Lagrange3D}. The mass flux is exactly zero only for isolated contact discontinuities moving in uniform flow when using appropriate Riemann solvers 
that resolve contact waves, like the Godunov method, or the Roe, HLLC, HLLEM and Osher flux.

\subsection{Connectivity matrices} 
\label{ssec.ConnectTables}

Since the core of the proposed method is the motion and the changing of the nonconforming mesh topology in time, we have to know all the connectivities of the mesh and to maintain them updated. 
In this way we will have enough information both to rearrange the mesh after the insertion of a new node, or the fusion of two existing nodes, and to know all the neighbors of 
each space--time lateral surface during the numerical flux computation.

%Since we allow our mesh to be nonconforming the connectivity tables are a little bit more complicate than in the conforming case, for example because more than two triangles can share the same edge and more than two points can belong to it. 

As in the standard conforming case for each element $T_i^n$ we save the global numbering of its vertices $V_1, \dots, V_{\nVertex} $ in row $i$ of a matrix called \talg{tri} 
%n $N_E \times \max(\nVertex$ matrix 
in counterclockwise order, and in matrix \talg{Elem2Edge} 
%with the same dimensions as the previous one 
we store the global numbering of its edges $E_1, \dots, E_{\nVertex}$. 
However, in the nonconforming case, additional connectivity tables are needed, since more than two elements can share the same edge and more than two points can belong to it. 
For each edge $E_j^n$, we store the elements that share it in row $j$ of matrix \talg{Edge2Elem},
% with $\Nl$ rows, 
and all the nodes that belong to $E_j^n$ in row $j$ of matrix \talg{Edge2Vertex} in such a way that the first two entries of each row contain the endpoints of the corresponding edge. 
Then, for each node we memorize the edge to which it belongs in \talg{Vertex2Edge} (both if this node is an endpoint of the edge or an intermediate point) and the elements for which it 
is a vertex in $\talg{Vertex2Elem}$ (note that if a node $N_i$ belongs to an edge of an element but it is not one of its vertices, that element will not appear in the row $N_i$ of this 
last matrix). Moreover, each node has a label that indicates whether the node is free to move everywhere, if it has been doubled, or if it is constrained to slide along a particular edge, 
i.e. if it is a hanging node. 

Besides, we allow our data structures to be completely dynamic in such a way that nodes and edges can appear and disappear in time: so rows can be added to our matrices or be nullified, 
and the information regarding which global numbering of nodes and edges is currently used is always available.

\subsection{Shear interface detector}
\label{ssec.DiscDetector}
Since the sliding interfaces are not defined \apriori  by the user, at the beginning of each time step the algorithm has first to identify along which \textit{edges} the shear interfaces 
lie, and mark the corresponding edges and nodes. Basically an edge $e$ will be considered at the interface if the tangential velocity difference $\Delta V_e$ across $e$ exceeds a 
certain threshold value  $\kappa_e$. So for each edge we need to compute $\Delta V_e$ and $\kappa_e$. 

Given the set of nodes $S_{\!\! e}^{\! n}$ over the edge $e$, and the set of neighbors $S_{\!\! j}^{\! n}$ of each node $j$, the threshold value $\kappa_e$ is computed as
\be
\label{eq.threshold}
\kappa_e = \min_{j \in S_{\!\! e}^{\! n}} \kappa_j,  \quad \text{ with} \quad \kappa_j = \max_{i\in S_{\!\! j}^{\! n}} \left ( \frac{\alpha  d_i}{ \vert \vert J_i \vert \vert} \right ),
\ee
where $d_i$ is the encircle diameter of element $T_i^n$, $J_i$ is the Jacobian of the transformation that maps element $T_i^n$ in the corresponding reference element, the norm is 
the two-norm of Frobenius divided by $\sqrt{2}$ (other matrix norms could also be used), and $\alpha$ is chosen in $[0,1]$ according to the desired sensitivity of the detector. 
If the velocity jump at the interface is very large, the value of $\alpha$ does not matter. Instead, where the velocity field changes smoothly, the number of interfaces, 
and as a consequence the number of new nodes, will be dependent on $\alpha$. Moreover, in the limit $\alpha \rightarrow +\infty$ we recover the standard conforming algorithm. 

Once  the threshold value has been fixed we loop over all the edges of the mesh: for each edge $e$ we consider all its neighbors and we  
compute their tangential velocity with respect to $e$. Say, for example, that two elements $A=T_a^n$ and $B=T_b^n$ with area $ \lvert T_a^n  \rvert$ and $\lvert T_b^n \rvert$ share the same edge $e$ and their tangential velocities are $v_{t,a}^{\,n}$ and $v_{t,b}^{\,n}$. If the quantity $\Delta V_e$  exceeds $\kappa_e$ 
\be
\label{eq.threshold2}
\Delta V_e = \frac{ \left \lvert \, v_{t,a}^{\,n} \, \lvert T_a^n \rvert - v_{t,b}^{\,n} \,\lvert T_b^n \, \, \rvert \right \rvert }{ \left ( \, \lvert  v_{t,a}^{\,n} \rvert \lvert T_a^n  \rvert + \lvert v_{t,b}^{\,n} \rvert \lvert T_b^n  \rvert \rvert + \epsilon \right ) }  \ {{\ge \kappa_e}} ,
\ee
with $\epsilon = 10^{-14}$ to avoid division by zero, then edge $e$ is marked as an edge at a shear interface, and the elements $A$ and $B$  are divided into two different groups: the elements on the left and the ones on the right with respect to this  particular edge $e$. 
Afterwards, we also need to find the \textit{nodes} that have \textit{to be doubled} and to separate their Voronoi neighbors  (the elements stored in \talg{Vertex2Elem}) into two groups. 
So we loop over the nodes considering the ones which belong to an interface edge. If in their list of Voronoi neighbors there are elements from both the sides of the interface, according to the 
previous subdivision, we mark them and we separate their Voronoi neighbors into two groups which are stored in two matrices.
% called \talg{Vertex2RightElem} and \talg{Vertex2LeftElem}. 

Note that the two cycles, the one over the edges and the other over the nodes, are not nested one into the other, but are run one after the other.

\subsection{Node motion}
\label{subsec.motion}

At this point we are able to distinguish between nodes far away from the interfaces, hanging nodes and nodes which lie at the interface. So we loop over the nodes and according to their labels we choose what to do.
%
%%%%%%%%%%%%%%%%% NORMAL NODES %%%%%%%%%%%%%%%%%
% 
First, consider a regular non-hanging node $k$ located in a smooth region. We compute its coordinates at the new time level $t^{n+1}$ simply by 
\be
\x^{n+1}_k = \x^n_k + \Delta t \overline{\mathbf{V}}^n_k,
\ee 
where $\overline{\mathbf{V}}^n_k$ is obtained using  the node solver of Cheng and Shu.  
Cheng and Shu introduced a very simple and general formulation for obtaining the final node velocity, which is chosen to be the arithmetic average velocity among all the contributions coming from the Voronoi neighbor elements $\mathcal{V}_k^n$. Moreover, following the ideas presented in \cite{LagrangeMHD} 
we take a mass weighted average velocity among the neighborhood $\mathcal{V}_k^n$, that is,
\be
\overline{\mathbf{V}}^n_k = \frac{1}{\mu_k} \, \sum_{T_j^n \in \mathcal{V}_k} \mu_{k,j}  \overline{\mathbf{V}}_{k,j}
\ee 
with
\be
\mu_k = \sum_{T_j^n \in \mathcal{V}_k} \mu_{k,j}, \quad \mu_{k,j} = \rho_j^n \lvert T_j^n\rvert.
\ee
The local weights $\mu_{k,j}$, which are the masses of the elements $T_j^n$, are defined by multiplying the cell averaged value of density $\rho_j^n$ (or water depth $h_j^n$ for shallow water flows) 
with the cell area $\lvert T_j^n \rvert$. The local  contributions $\overline{\mathbf{V}}_{k,j}$ in a pure Lagrangian context represent the fluid velocity in the  $j^{th}$ neighbor of vertex $k$, 
while in the ALE framework they can be obtained either according to an arbitrary, prescribed mesh velocity function or by the local fluid velocity. 
% 
%%%%%%%%%%%%%%%%%%  INTERFACE NODES %%%%%%%%%%%%%%%%%
%
Now let us consider the nodes at the interfaces. The following considerations are carried out by supposing for the moment that each interface is separated from the others and lies on a straight line. 
% even if this is a rigid constraint, already with this configuration some interesting test cases can be studied.
% and the generalization to the general case of piece-wise linear interfaces will be the object of another work. 
%Note, however, that the only constraint we have is a straight line interface, but we do not need to choose its position \apriori.

%%% --- FIRST: INSERTION --- %%%%
\subsubsection{Insertion of a new node}

The first situation we encounter is a node $k$ that has some of its Voronoi neighbors on the left of the interface, call them \textit{left  neighbors}, $\mathcal{V}_{k,\text{left}}$, and the others on the right of the same interface, call them \textit{right neighbors}, $\mathcal{V}_{k,\text{right}}$; these two sets of neighbors have been provided by the \textit{interface detector}  described above.

We apply the node solver of Cheng and Shu at the two sets of neighbors obtaining two different new coordinates 
\be
\label{eq.newTildeCoord}
\tilde{\x} ^{n+1}_{k,\,\text{left}} = \x^n_k + \Delta t  \!\! \sum_{T_j^n \in \mathcal{V}_{k,\text{left}}} \frac{\mu_{k,j} }{\mu_k}  \overline{\mathbf{V}}_{k,j}, 
\quad \text{and} \quad  \tilde{\x} ^{n+1}_{k,\,\text{right}} = \x^n_k + \Delta t \!\! \sum_{T_j^n \in \mathcal{V}_{k,\text{right}}} \frac{\mu_{k,j} }{\mu_k}  \overline{\mathbf{V}}_{k,j} \, .
\ee
We allow this kind of nodes to move only along the interface, so basically according to their averaged tangential velocity with respect to the interface: for this reason we need to find the nearest interface edges and to project over them the coordinates in (\ref{eq.newTildeCoord}) obtaining thus 
${\x} ^{n+1}_{k,\,\text{left}}$ and ${\x} ^{n+1}_{k,\,\text{right}}$.

Call the nearest interface edges belonging to the left elements $e_1^{\ell}$ and $e_2^\ell$, and the nearest interface edges belonging to the right elements $e_1^r$ and $e_2^r$ (suppose also that $e_1^{\ell,r}$ are closer to $\tilde{\x} ^{n+1}_{k,\,\text{left}}$ than to $\tilde{\x} ^{n+1}_{k,\,\text{right}}$).

We decide to assign $\x ^{n+1}_{k,\text{left}}$ as new coordinate to the old node $k$ 
\be
\x_k^{n+1} = \x ^{n+1}_{k,\text{left}}
\ee
and to create a new node with global number $k_{\text{new}}$ and coordinates (at time $n$ and $n+1$)
\be
\x^n_{k_{\text{new}}} = \x^n_k  \quad \text{and} \quad  \x^{n+1}_{k_{\text{new}}} = \x ^{n+1}_{k, \,\text{right}}.
\ee
The global number $k_{\text{new}}$ can be larger than $\Nn$ if all the numbers between $1$ and $\Nn$ are currently used, otherwise we choose the first of the unused numbers (indeed if two nodes have been merged together one of their global numbers is no more used, see section \ref{sssec.fusion}).

\begin{figure}
	\centering
	\includegraphics[width=0.48\textwidth]{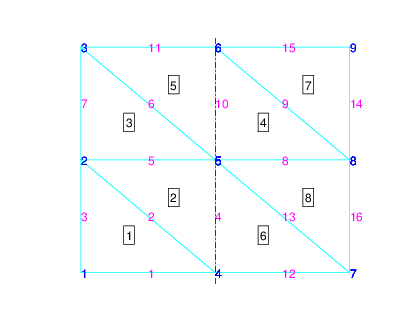} \\
	\includegraphics[width=0.48\textwidth]{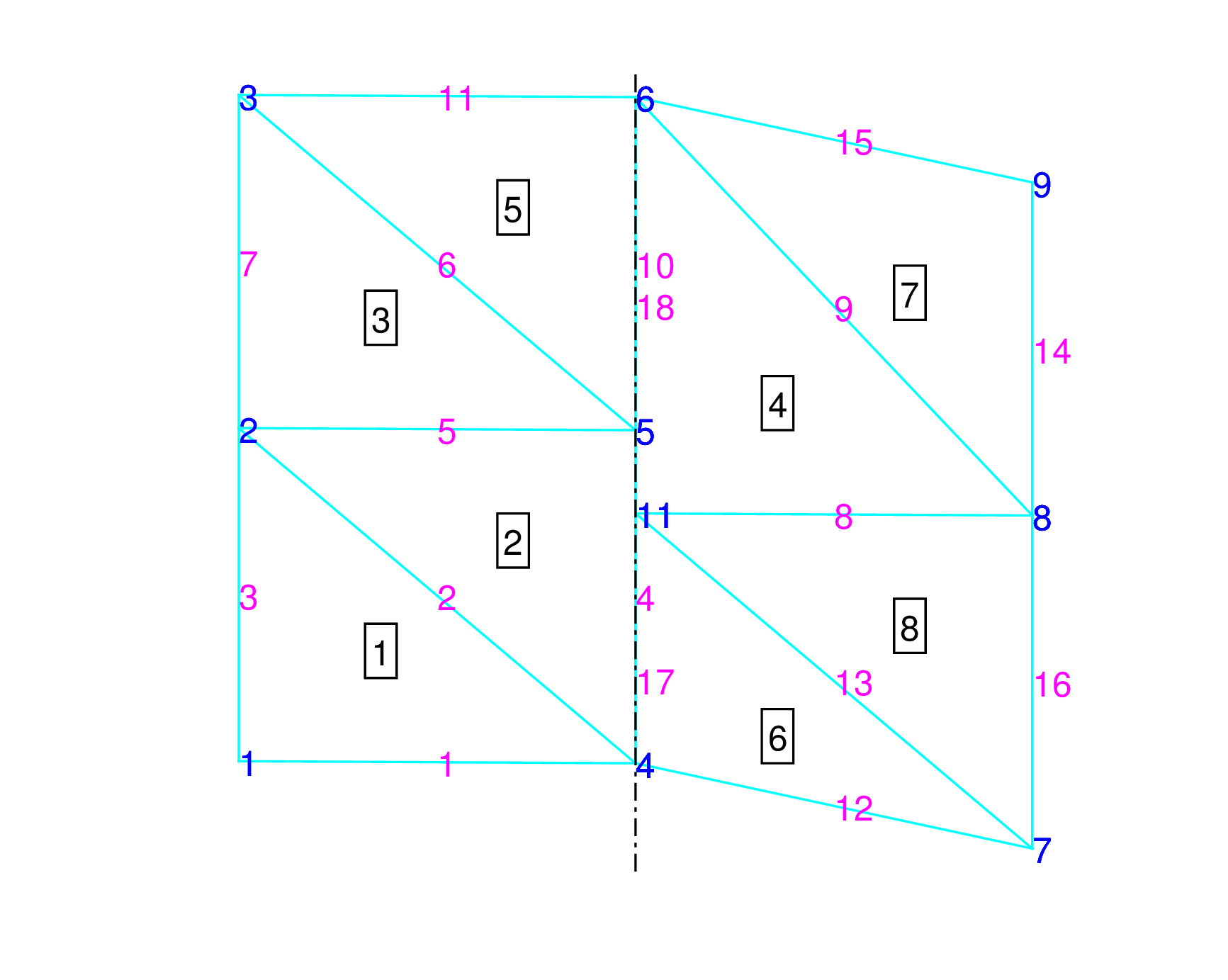}
	\includegraphics[width=0.48\textwidth]{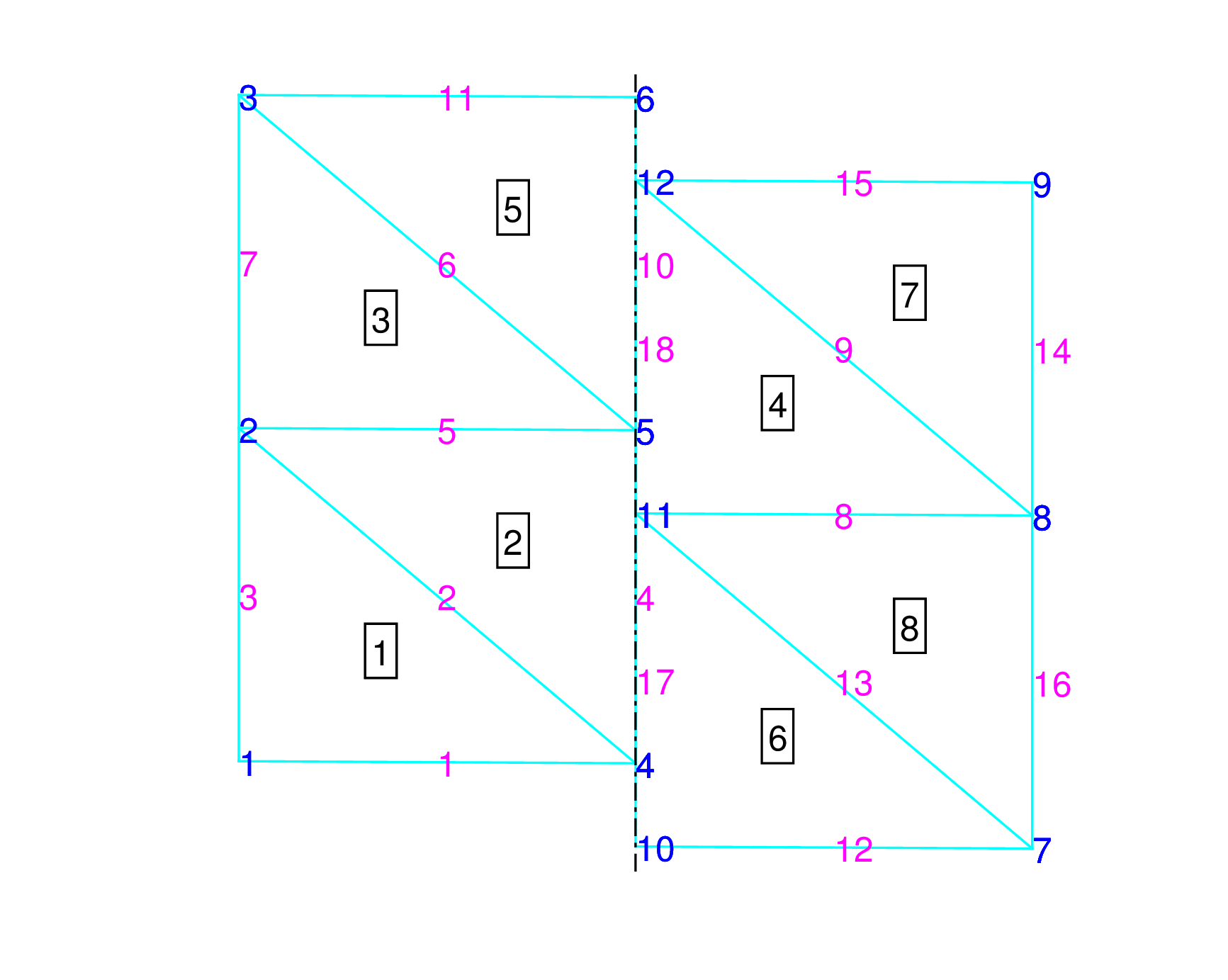}
	\caption{Example of how to double a node.  At the top we show the situation before a nonconforming motion, and at the bottom  after the motion and the corresponding update of the connectivity matrices. Precisely at the bottom on the left we have supposed to move in a nonconforming way only one of the nodes at the interface (for this reason the mesh is deformed, but we did it only to explain clearly one step of the algorithm), while the realistic motion of all the nodes at the interface is shown on the right.		
		The black vertical dotted line represents the interface: suppose that the elements on the left $\{1,2,3,5\}$ move with velocity $\v = (0,2)$ and the elements on the right $ \{4,6,7,8\}$ move with velocity $\v=(0,-2)$. We want to double vertex number $k=5$, so we insert a new node $ k_{\text{new}} = 11$.
		The nearest interface edges on which we project the new coordinates of node $5$ are $e_1^\ell=e_1^r = 10$  and $e_2^\ell = e_2^r =4$. Note that edges $e_1^{\ell,r}$ are closer to $k$ than to $k_\text{new}$. Since the edges from the left and from the right are equal we create two new edges $e_{1_{\text{new}}}^r = 18$ and $e_{2_{\text{new}}}^r =17$. The endpoints of edges $10$ and $4$ remain untouched. 
		Edge 4 gains an intermediate point, the node $11$, and edge $18$ gains the node $5$. 
		To better understand we list now the vertices of each edge at the end of the updating process (first we write the endpoints and then, if existing, the intermediate points): $e_1^\ell =10 \rightarrow \{5,8\}$, $e_2^\ell = 4\rightarrow \{4,5,11\}$, $e_1^r = 18\rightarrow \{11,6,5\}$ and $e_2^r =17\rightarrow\{4,11\}$.
		Finally, elements $\{1,3,5,6,7,8\}$ maintain the same edge neighbors, while the neighbors of elements $2$ and $4$ are augmented: indeed edge $4$ has neighbors $\{2,6,4\}$ and edge $18$ has neighbors $\{4,5,2\}$. 	
		Note that the situation on the right appears to be more complicated only because also nodes $4$ and $6$ have been doubled and so the corresponding update of the connectivity matrices has been done.
	}
	\label{fig.insertion}
\end{figure}

Now we have to update the connectivity tables taking into account the insertion of this new node. See also Figure~\ref{fig.insertion} to follow our construction.

First, in matrix \talg{tri} we substitute $k$ with $k_{\text{new}}$ in all the right elements; moreover, we consider matrix \talg{Vertex2Elem} and in row $k$ we leave only the left elements and we put the others in row $k_{\text{new}}$ (because now $k$ is no more a vertex for the right neighbors).

Then we have to deal with the edges: if $e_1^\ell=e_1^r$ we need to substitute $e_1^r$ with a new edge $e_{1_{\text{new}}}^r$. 
In matrix \talg{Elem2Edge} all the right neighbors change $e_1^r$ with $e_{1_{\text{new}}}^r$, and in matrix \talg{Edge2Elem} we insert a new row $e_{1_{\text{new}}}^r$ equal to row $e_1^r$ (the new edge inherits all the characteristics from the old one). 
The same has to be done if $e_2^\ell = e_2^r$. 
The endpoints of these new edges are the endpoints of the substituted edges seen from the right (so basically there is $k_\text{new}$ instead of $k$). The endpoints of the left edges do not change.  
Besides we add $k$ as intermediate point in $e_1^r$ and $k_\text{new}$ as intermediate point of $e_2^\ell$, (note that an edge is allowed to have more than one intermediate point).
In this way also matrix \talg{Edge2Vertex} has been updated. Matrix \talg{Vertex2Edge} is easily modified at the same time.

Finally, we have to revise the list of neighbors: in particular the edges that gained an intermediate point ($e_1^r$ and $e_2^\ell$)  gain also neighbors. In particular the new neighbors of $e_1^r$ are the left neighbors of $e_2^\ell$ and the new neighbors of $e_2^\ell$ are the right neighbors of $e_1^r$. This allow us to update \talg{Edge2Elem} and \talg{Elem2Edge}.

At the end we mark with a label the nodes which are intermediate for an edge: we call them \textit{hanging nodes} and they are constrained to move along that edge.
Note that in the case of straight slip-lines no distinction between master and slave nodes is required, since both will move along the same straight interface. To extend the algorithm to the case of piece-wise linear interfaces, this distinction is introduced in such a way that only slave nodes will be constrained to slide along edges, while the master nodes can move freely. For some first preliminary 
results concerning the extension to completely general slide lines, see Section \ref{sec.general} of this paper.

%%% --- SECOND: sliding --- %%%

\subsubsection{Hanging nodes}

Consider a hanging node $k$ which lies on the edge $e$: it is at the interface and it is a vertex only of elements lying on one side of the interface, so its Voronoi neighbors are in the same smooth region. However it is not free to move everywhere but it must slide along that edge, to avoid creation of holes or superposition of elements in the mesh.

We apply the averaged node solver of Cheng and Shu among its Voronoi neighbors, we find its new coordinates $\tilde{\x}^{n+1}_k$ and we project them over edge $e$, obtaining $\x_k^{n+1}$.
Now, we compute also the new coordinates of the other points over edge $e$. If two of them are sufficiently close, we decide to merge them (see section \ref{sssec.fusion}), otherwise the computed coordinates $\x_k^{n+1}$ are the new coordinates of such a node and no update of the connectivity matrices is required.

%%% --- THIRD: fusion --- %%%
\subsubsection{Fusion of two existing nodes}
\label{sssec.fusion}

\begin{figure}[h]
	\centering
	\includegraphics[width=0.49\textwidth]{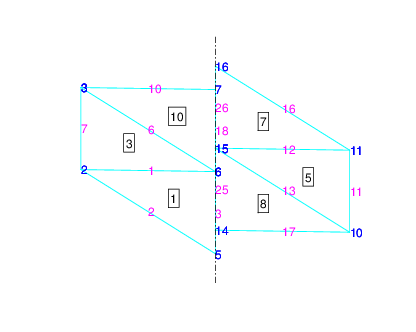} 
	\includegraphics[width=0.49\textwidth]{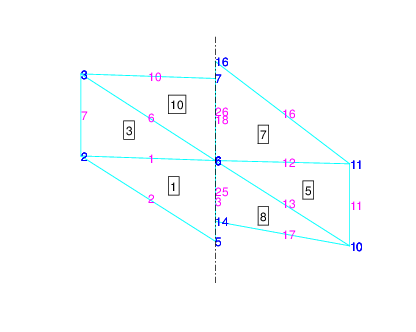} 
	\caption{Example of how to merge two existing nodes. The black dotted line represents the interface: suppose that the elements on the left $\{1,3,10\}$ move with positive velocity and the elements on the right $\{5,7,8\}$ move with negative velocity. On the left we show the situation at time $t^n$ and on the right at time $t^{n+1}$. Nodes $6$ and $15$ at $t^{n+1}$ will be so close that we decide to merge them (as in the previous example, for the sake of clarity, we present on the right the situation after the fusion of only two nodes).	We maintain the smallest global number so $k_\text{fn} =6$ and we remove $k_\text{dn} =15$. In \talg{triNew} elements $\{5,7,8\}$ substitute $k_\text{dn} =15$ with  $k_\text{fn} =6$. Note that in \talg{tri} nothing changes, so some elements refer to node $6$ and some other to node $15$, but everything works because at time $t^{n+1}$ they have the same new coordinates $X_{k_\text{dn}}^{n+1} = X_{k_\text{fn}}^{n+1} $ and at the successive time step $t^{n+2}$ \talg{tri} will no longer exist because it will be overwritten by \talg{triNew}. In row $k_\text{fn}$ of matrix \talg{Vertex2ElemNew} there are listed elements $\{1,3,5,7,8,10\}$, while row $k_\text{dn}$ is empty. In row $k_\text{fn}$ of matrix \talg{Vertex2EdgeNew} there are edges $\{1,3,6,12,13,18,25,26\}$, while row $k_\text{dn}$ is empty.
		List \talg{Edge[dn-fn]} contains edges $\{18,25\}$ 	and list  \talg{Elem[dn-fn]} contains elements $\{8,10\}$. Knowing these lists we can update  matrices \talg{Edge2ElemNew} because we remove element $8$ from the neighbor of edge $18$ and element 10 from the neighbors of edge $25$. In this case even if we removed the segment $\overline{6,15}$ no edge becomes equal so we do not need to merge edges neither to update \talg{Elem2EdgeNew}.  		
	}
	\label{fig.fusion}
\end{figure}

Suppose we computed the new coordinates at time $t^{n+1}$ of all the nodes $k_i$ over the same edge $e$ denoted by $X_{k_i}^{n+1}$, which are assumed to be already projected 
onto the straight line spanned by edge $e$. If the new coordinates of two of them, say $k_1$ and $k_2$, are too close, we decide to merge them. Moreover, if one intermediate 
node of edge $e$ falls outside the edge, we decide to merge it with the closest endpoint of the edge. 

Since the loop over the nodes is carried out according to the increasing global numbering of the nodes, we decide to remove the node with the largest global number (we call it \textit{dead node}, $k_{\text{dn}}$) because we have not worked with it yet, and to maintain the one with the smallest global numbering (call it \textit{fusion node}, $k_{\text{fn}}$) assigning to it as new coordinates the average between $\x_{k_1}^{n+1}$ and $\x_{k_2}^{n+1}$
\be
\x_{k_\text{fn}}^{n+1} = \frac{\x_{k_1}^{n+1} + \x_{k_2}^{n+1} }{2}.
\ee
We assign the same coordinates also to the dead node 
\be 
\label{eq.CoordDeadNode}
\x_{k_\text{dn}}^{n+1} = \frac{\x_{k_1}^{n+1} + \x_{k_2}^{n+1} }{2}.
\ee
Now, we need to update the connectivity tables. See also Figure~\ref{fig.fusion} to follow our construction.

This process is somehow more complicated than the nodes splitting. Indeed when we insert a new node at time $t^{n+1}$ we only add information without losing anything about the previous time step, and even if it is true that the right neighbors of a doubled node $k$ change their node $k$ with a new one $k_\text{new}$, we can dispose of all its reference simply by giving to $k_\text{new}$ at time $t^n$ the same coordinates of $k$, see also (\ref{eq.newTildeCoord}). On the contrary, when we remove a node we lose all the reference to it, reference that, only for time $t^{n+1}$, we still need during the computation of the interface fluxes in the finite volume scheme (it is for this reason that in \eqref{eq.CoordDeadNode} we have assigned the coordinates  $\x_{k_\text{dn}}^{n+1}$ even to the dead node).
So we decide to duplicate some of the connectivity tables, creating \talg{triNew}, \talg{Elem2EdgeNew}, \talg{Edge2ElemNew}, \talg{Edge2VertexNew}, and \talg{Vertex2ElemNew}.
During the insertion procedure we modify in the same way both the old and the new matrices, while during the fusion we modify only the \talg{new} matrices. Hence we can use the old ones in the finite volume scheme, because they store all the needed information (for example they refer both to the dead and the fusion node which have the same coordinates at the new time $t^{n+1}$), while when we advance in time, to $t^{n+2}$, we maintain updated only the new ones because the information about two previous time steps are no longer necessary and so we can overwrite the old connectivity 
matrices with the new ones. 

First, in matrix \talg{triNew} we substitute $k_\text{dn}$ with $k_\text{fn}$ in all the neighbors of the dead node; moreover, we consider matrix \talg{Vertex2ElemNew}, in row $k_\text{fn}$ we put both the neighbors of the dead and the fusion node and we nullify row $k_\text{dn}$. We do the same with matrix \talg{Vertex2EdgeNew}: we nullify row $k_\text{dn}$ and we put in row $k_\text{fn}$ all the edges that contain $k_\text{fn}$ or $k_\text{dn}$.

Then all the edges that contain $k_\text{dn}$ substitute it with  $k_\text{fn}$ (in matrix \talg{Edge2VertexNew}), whereas the edges with both $k_\text{dn}$ and $k_\text{fn}$ (that we memorize in a list \talg{Edge[dn-fn]}) remove $k_\text{dn}$.
We note that merging $k_\text{dn}$ and $k_\text{fn}$ we are removing the segment in between, so we look for the edges that contain it (listed in \talg{Edge[dn-fn]}) and its neighbor elements that we list in \talg{Elem[dn-fn]}.
We update now matrix \talg{Edge2ElemNew}  because the edges in \talg{Edge[dn-fn]} have no more one of the neighbors in \talg{Elem[dn-fn]} . 

Afterward we check if the absence of this segment makes some edges in \talg{Edge[dn-fn]} equal: in this case we remove one of them (the one with the largest global number) and we update correspondingly the new connectivity matrices.

Besides we modify the labels telling us if a node is hung to some edges and which nodes and edges are currently existing. This last passage prevents us to work again with disappeared nodes and allows us to reuse their global numbering when we want to insert a new node or a new edge.

\section{Numerical results}
\label{sec.results}

In this section, we solve a large set of numerical tests in order to validate the presented nonconforming direct ALE scheme. The robustness of the method is checked both on smooth and discontinuous problems related to the shallow water equations written in Cartesian and in polar coordinates. The test cases are carried out using either the Rusanov or the Osher type flux, the value of $\alpha$ in \eqref{eq.threshold} is always taken equal to $\alpha =1$ unless otherwise specified, and the CFL number is chosen as CFL$=0.3$.
Furthermore, the order of convergence is verified first fixing for the mesh motion an arbitrary velocity, then in the case of a steady vortex in equilibrium using the local fluid velocity.

\subsection{Sanity checks: pure sliding}

\begin{figure}
	\centering
	\subfloat[initial mesh t = 0]{ \includegraphics[width=0.3\textwidth]{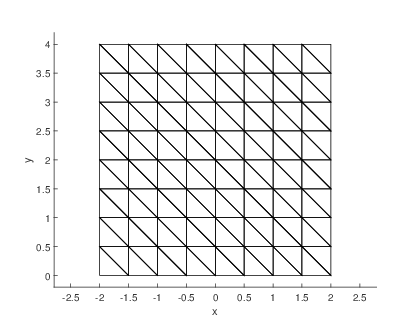} 	}
	\subfloat[mesh at t = 0.5]{    \includegraphics[width=0.3\textwidth]{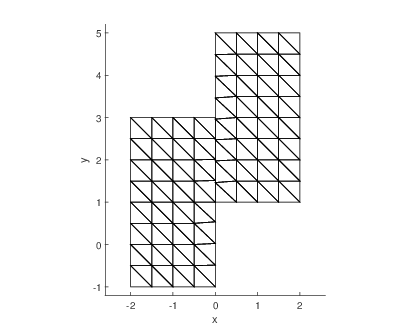} }	
	\subfloat[mesh at t = 1.1]{    \includegraphics[width=0.3\textwidth]{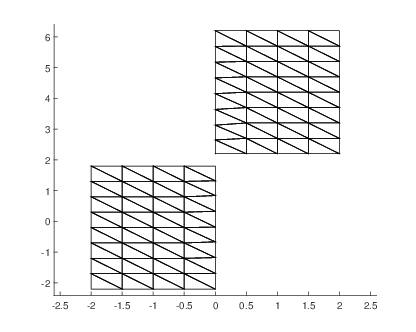}}  \\
	\subfloat[initial mesh t = 0]{ \includegraphics[width=0.3\textwidth]{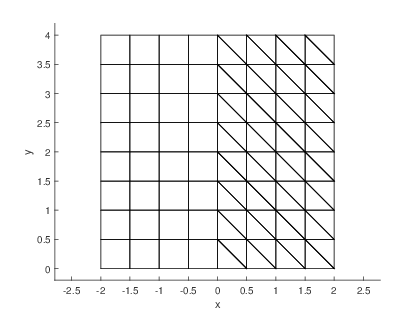} }	
	\subfloat[mesh at t = 0.4]{ 	 \includegraphics[width=0.3\textwidth]{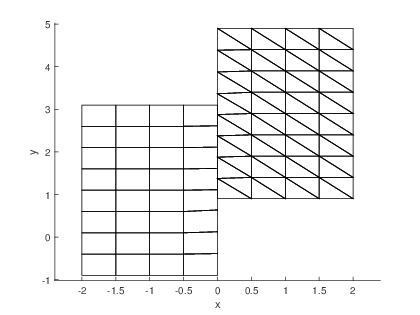} 	}
	\subfloat[mesh at t = 1.1]{ 
		\includegraphics[width=0.32\textwidth]{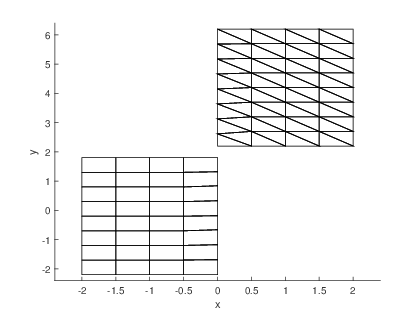} }
	\caption{Slide lines test case with initial condition as in Equation \eqref{eq.SanityCheck1}. The mesh is moved with the local fluid velocity, which at $x=0$ is discontinuous: so nodes over there are handled in a nonconforming way. At the top we show the results obtained employing a triangular mesh and at the bottom using a mesh made of both triangular and quadrilateral elements. We report the mesh at three different computational times: note that the computational domain can also be split in two non connected parts. The level of the water, the total area and the total volume are conserved at any time step, and the solution coincides with the exact one up to machine precision.}
	\label{fig.tc1}
\end{figure}

The numerical examples reported in this section are sanity checks testing the ability of the method to detect and maintain straight slip-line 
interfaces. We consider the shallow water equations, which can be cast into form (\ref{eq.generalform}) with 
\be
\label{eq.shallowWater}
\Q = \left( \begin{array}{c} h \\ hu \\ hv \end{array} \right), \quad \f = \left( \begin{array}{c} h u \\ hu^2 + \frac{1}{2} g h^2 \\ h uv \end{array} \right), \quad \g = \left( \begin{array}{c} h v \\ h uv \\ h v^2 + \frac{1}{2}gh^2  \end{array} \right).  
\ee
The initial computational domain is given by $\Omega(t_0) = [-2,2]\times[0,4]$. 
First, we take the initial condition
\begin{equation}
\label{eq.SanityCheck1}
	\Q(\x, 0) = \left\{ \begin{array}{ccc} \left( 1, 0, -2     \right) & \textnormal{ if } & x \leq 0, \\
		\left( 1, 0, 2 \right) & \textnormal{ if } & x > 0,        
	\end{array}  \right. 
\end{equation} 
which also coincides with the exact solution at any time. 
We impose \textit{wall} boundary conditions on the left and on the right side of the domain, respectively, whereas at the top and at the bottom of the domain we impose \textit{transmissive} boundary conditions.  
In Figure~\ref{fig.tc1} we show the numerical results over a triangular mesh and then over a mixed mesh composed of both, triangular and quadrilateral elements. The chosen mesh velocity coincides 
exactly with the fluid velocity, as in a pure Lagrangian context. At each time step we have verified that the total water volume is conserved up to machine precision both locally and globally and 
that relation (\ref{eqn.gcl}), the GCL, is verified also up to machine precision. 

\begin{figure}[!htbp] 
	\centering
	\subfloat[initial mesh t = 0]{            \includegraphics[width=0.3\textwidth]{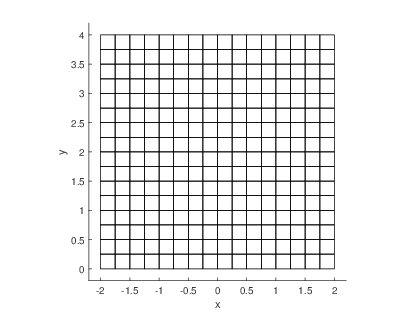} }	
	\subfloat[final mesh with $\alpha=1$]{ 		\includegraphics[width=0.3\textwidth]{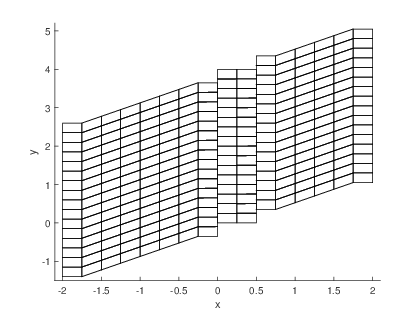} 	}
	\subfloat[final mesh with $\alpha = 0.4$]{ 
		\includegraphics[width=0.32\textwidth]{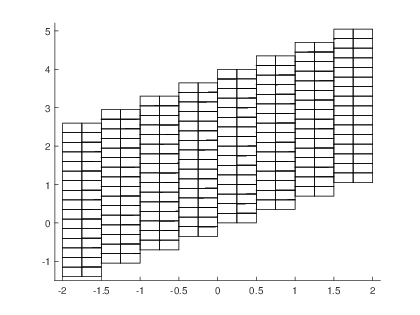} }
	\caption{Slide lines test case with initial condition as in Equation \eqref{eq.initialCondFloor}. We start with a conforming quadrilateral mesh; using a value of $\alpha=1$ in \eqref{eq.threshold} we obtain only two slip-lines (at $x=0$ and $x=0.5$), whereas using $\alpha = 0.4$, which makes the detector more strict, the mesh slides along each straight line where the fluid velocity changes.}
	\label{fig.tc2_moreInterfaces}
\end{figure}

\begin{figure}[!htbp] 
	\centering
	\subfloat[initial mesh t = 0]{        \includegraphics[width=0.3\textwidth]{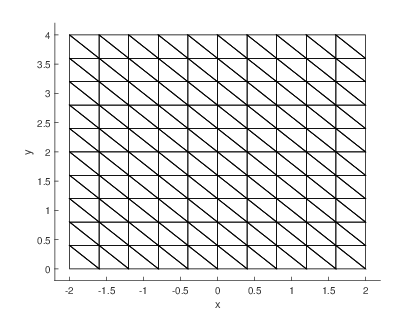} }	
	\subfloat[mesh at t = 0.35]{ 					\includegraphics[width=0.3\textwidth]{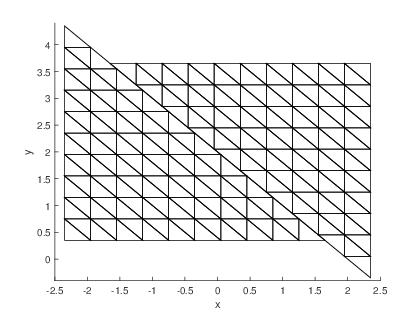} 	}
	\subfloat[mesh at t = 0.7]{ 
		\includegraphics[width=0.32\textwidth]{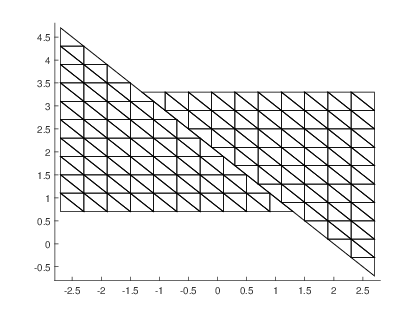} }
	\caption{Oblique slide line. We show the discretization of the computational domain at three different times. The corresponding numerical solution matches the exact one.}
	\label{fig.tc3_obliqua}
\end{figure}

Next, we consider as initial condition 
\be
\label{eq.initialCondFloor}
\Q(\x,0) = \left ( 1, 0, 0.5\, \text{floor}\,(2x) \right ), \quad -2 \le x \le 2,
\ee
with $\text{floor}(x) = \lfloor x \rfloor$ denoting the lower Gauss bracket, and we run our algorithm until a final time $t=0.7$ with different \textit{threshold} values, see (\ref{eq.threshold}), in such a way that there will be a different number of interfaces along which nodes have to be doubled and merged in time. The discretization of the computational domain is reported in Figure \ref{fig.tc2_moreInterfaces}. Also in this case we reach the exact solution (that is the initial condition translated in the motion direction), the total volume of water is conserved and relation (\ref{eqn.gcl}) is verified up to machine precision at each time step and on each element.

Finally, we want to show that the interface can be along any straight line (provided that edges lie over this line): we take as initial condition
\begin{equation}
\label{eq.DiscObliqua}
\Q(\x, 0) = \left\{ \begin{array}{ccc} \left( 1, -1, 1     \right) & \textnormal{ if } & x+y-2 \leq 0, \\
\left( 1, 1, -1 \right) & \textnormal{ if } & x+y-2 > 0,        
\end{array}  \right. 
\end{equation} 
and in Figure \ref{fig.tc3_obliqua} we report the computational domain at different times. Again, the numerical solution matches the exact one and as expected, the total volume 
is conserved and equation \eqref{eqn.gcl} is satisfied up to machine precision.

\subsection{Periodic boundary conditions}

\begin{figure}[!htbp] 
	\centering
  \includegraphics[width=0.35\textwidth]{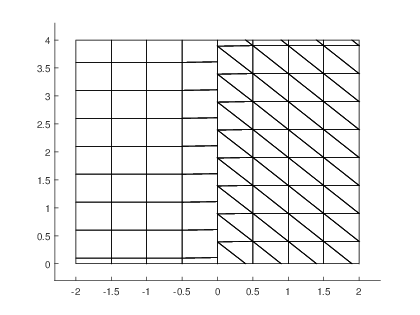} 	
	\includegraphics[width=0.35\textwidth]{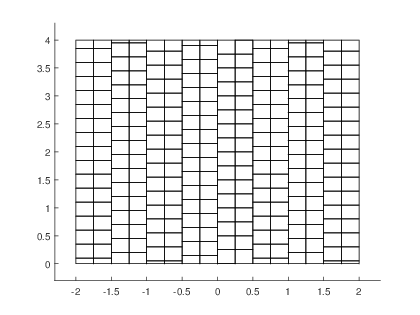} 	
	\caption{Slide lines with periodic boundary conditions. We report the final computational domain at time $t=100.2$ corresponding to the initial condition in \eqref{eq.SanityCheck1} on the left, and the one corresponding to the initial condition in \eqref{eq.initialCondFloor} on the right. No distortion of the computational domain appears neither at the interfaces, and the numerical solution coincides with the exact one.}
	\label{fig.tc_periodica}
\end{figure}

The tests reported in the previous section can be run also by imposing periodic boundary conditions on the top and at the bottom of the computational domain. In Figure \ref{fig.tc_periodica} we show the discretization of the computational domain at time $t = 100.2 $ for the initial conditions in \eqref{eq.SanityCheck1} and in \eqref{eq.initialCondFloor}.  
We would like to underline that no distortion of the mesh elements appears even after a very long computational time, and as a direct consequence the time steps remain almost constant during the computation. As always in this type of test cases the volume conservation holds and the numerical solution is equal to the exact one up to machine precision.

\subsection{Riemann problem}
\begin{figure}
	\begin{tabular}{ccc} 
	\includegraphics[width=0.3\textwidth]{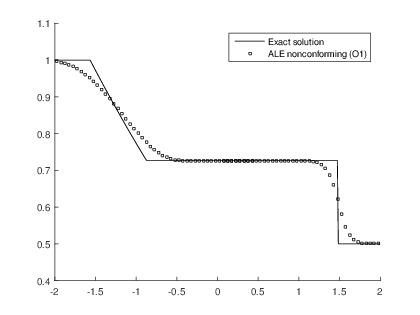} 	&
	\includegraphics[width=0.3\textwidth]{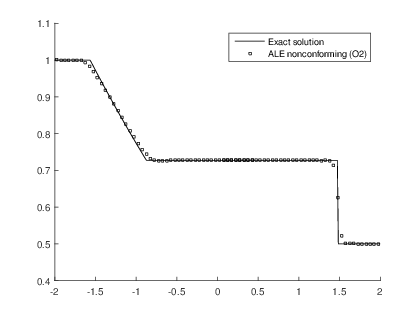} 	& 
	\includegraphics[width=0.3\textwidth]{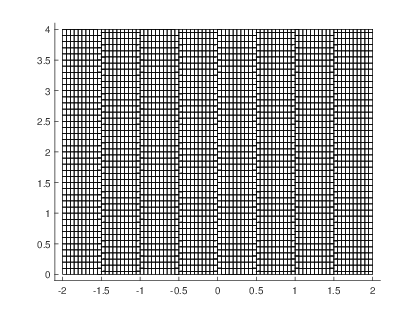} 
	\end{tabular} 
	\caption{Riemann problem with an arbitrary mesh velocity. Taking $\alpha =0.4$ in \eqref{eq.threshold} the algorithm identifies $7$ interfaces which are then handled in a nonconforming way. In the  figure we report the final discretization of the computational domain, and the comparison between the exact solution and the numerical solutions obtained with our nonconforming method showing first
	order results (left), second order results (center) and the mesh at the final time (right). }
	\label{fig.RiemannPb}
\end{figure}

Let us now consider as initial condition a Riemann problem with a discontinuity in the water level
 \begin{equation}
 \label{eq.DiscObliqub}
 \Q(\x, 0) = \left\{ \begin{array}{ccc} \left( 1, 0, 0     \right) & \textnormal{ if } & x \leq 0, \\
 \left( 0.5, 0, 0 \right) & \textnormal{ if } & x > 0,        
 \end{array}  \right.
 \end{equation} 
that originates a left-traveling rarefaction fan and a right-moving shock wave. We decided to move the mesh with an arbitrary mesh velocity function 
\[
\mathbf{V} = \left(\, 0, \, 0.5\,\text{floor}\,(2x)\, \right)   \quad -2 \le x \le 2,
\]
in order to check the robustness of the algorithm also in the presence of discontinuities. We imposed periodic boundary conditions on the top and on the bottom of the square,  and wall boundary conditions on the left and on the right.
The final discretization of the computational domain together with the comparison between the numerical and the exact solution are depicted in Figure \ref{fig.RiemannPb} both for the first order 
accurate scheme (i.e. without the MUSCL-Hancock strategy for the reconstruction) and the second order accurate scheme. 

\subsection{Convergence test}
\begin{figure}
	\centering
	\begin{tabular}{ccc} 
    \includegraphics[width=0.3\textwidth]{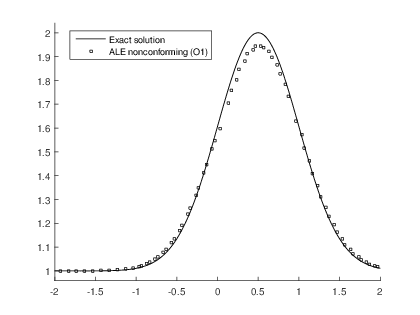} & 
   	\includegraphics[width=0.3\textwidth]{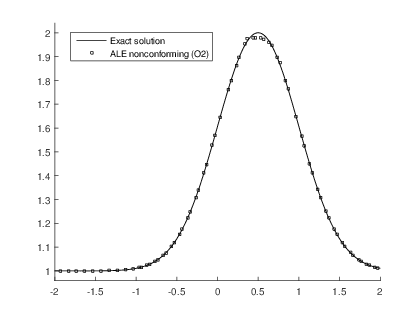} & 
		\includegraphics[width=0.3\textwidth]{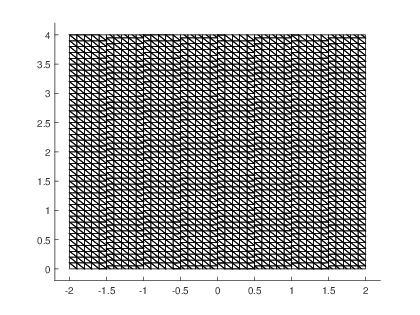} 		
	\end{tabular} 
	\caption{Comparison of the exact solution for the quantity $c$ with the numerical solution obtained on moving nonconforming meshes. The results obtained with the first order algorithm are shown on   
	         the left, while those obtained with the second order MUSCL-Hancock method are presented in the center. The comparison is done at time $t=0.5$ taking a cut of the profile of $c$ corresponding  
	         to $y=2$. On the right we show the discretization of the computational domain at time $t=0.5$. }  
	\label{fig.GaussianProfile}
\end{figure}

\begin{table} 
	\begin{center} 	
		\begin{tabular}{cccc|cccc} 
			\hline 
			&   $\mathcal{O}1$  &		& &  & $\mathcal{O}2$  &   &    \\ 
			\hline
			mesh points &  $h(\Omega(t_f))$ & $\epsilon_{L_2}$ & $\mathcal{O}(L_2)$ &  mesh points &  $h(\Omega(t_f))$ & $\epsilon_{L_2}$ & $\mathcal{O}(L_2)$  \\ 
			\hline
			12 $\times$ 12 & 1.95E-01 & 1.44E-01 &	-	 &  12 $\times$ 12   & 1.95E-01 & 4.96E-02 &  -       \\
			24 $\times$ 24 & 9.78E-02 &	7.58E-02 &  0.93 &	24 $\times$ 24   & 9.78E-02 &	1.23E-02 &  2.02	  \\	
			40 $\times$ 40 & 5.88E-02 &	4.69E-02 & 	0.94 & 	40 $\times$ 40   & 5.88E-02 &	4.24E-03 &  2.10    \\	
			80 $\times$ 80 & 2.95E-02 &	2.41E-02 &	0.97 &	80 $\times$ 80   & 2.95E-02 &	1.01E-03 &	  2.09    \\	
			120 $\times$ 120 & 1.98E-02 &	1.62E-02 &	0.99 &	120 $\times$ 120   & 1.98E-02 &	4.51E-04 &  2.01    \\	
			
			\hline 
		\end{tabular}		
	\end{center}
	\caption{Numerical convergence results for the passive transport of a Gaussian profile on moving nonconforming meshes. The error norms refer to the variable $c$ at time $t=0.5$. On the left we report the result for the first order method (i.e. without using the MUSCL-Hancock reconstruction procedure) and on the right using the second order accurate scheme.} 
	\label{tab.orderOfconvergence}
\end{table}

To verify the order of convergence of the proposed method we study the passive transport of a quantity $c$, that at time $t=0$ is taken equal to a Gaussian profile and then will be passively 
transported in the  direction of the fluid flow without changing its shape. 
The PDE system describing this situation is obtained from the standard shallow water equations (\ref{eq.shallowWater}) with the addition of the concentration $c$ of a passive tracer, 
\be
\label{eq.shallowWaterPassiveTransport}
\Q = \left( \begin{array}{c} h \\ hu \\ hv \\hc \end{array} \right), \quad \f = \left( \begin{array}{c} h u \\ hu^2 + \frac{1}{2} g h^2 \\ h uv \\ huc \end{array} \right), \quad \g = \left( \begin{array}{c} h v \\ h uv \\ h v^2 + \frac{1}{2}gh^2  \\ hvc \end{array} \right).  
\ee
We fix the following initial condition
\be
\label{eq.initialCondGaussian}
\Q(\x,0) = \left (1, u, 0, 1 + \textup{e}^{\frac{- 0.5 \, \left ( x^2 + \, \left (y \,-\, 0.5\, p \right )^2 \right )  }{ 0.5  ^2}} \right ), \quad -2 \le x \le 2, \quad 0 \le y \le p,
\ee
where we use a fluid velocity of $u=1$ and where we have taken the period $p=4$. The mesh is moved with the velocity 
\be
\mathbf{V} = \left (0, \,0.5\, \text{floor}\,(x) \right ) \quad -2 \le x\le 2,
\ee
according to the ALE framework, where the mesh velocity can be chosen arbitrarily. We prescribed periodic boundary conditions on the upper and lower side of the rectangular domain, and wall boundary conditions on the left and right sides.

Since the exact solution is known ($\Q(\x,t)=\Q(\x-u t,0)$) and it is smooth, we can verify the order of convergence of our method. In Table \ref{tab.orderOfconvergence} we report the order of convergence of the basic first order finite volume method, and of its second order extension that uses the MUSCL-Hancock strategy for the reconstruction procedure in space and time. 
Moreover, in Figure \ref{fig.GaussianProfile} we compare the numerical solution for the variable $c$ with the profile of the exact solution and we show the mesh at the final time.

\subsection{Steady vortex in equilibrium}
\label{ssec.Numresults_vortex}
To show that our method is also robust enough for vortex flows, we simulate the case of a steady vortex in equilibrium and we will compare the results obtained with our nonconforming method 
with a standard conforming algorithm (without any rezoning technique) looking at the differences after long simulation times. 
First, we rewrite the shallow water equations \eqref{eq.shallowWater} in polar coordinates. 
Consider the usual relation between polar $(r, \phi)$ and Cartesian $(x,y)$ coordinates  
\be 
	x = r \cphi, \qquad \textnormal{ and } \qquad 
	y = r \sphi, 
\ee 
and the corresponding relations for the derivatives 
\be
\label{eq.changeOfDer}
	\de{}{x} =  \, \cphi  \de{}{r} - \frac{\sphi}{r} \de{}{\phi}, \qquad \textnormal{ and } \qquad 
	\de{}{y} = \, \sphi  \de{}{r} + \frac{\cphi}{r} \de{}{\phi}\,  
\ee
and let $\urho$ and $\uphi$ be respectively the radial and the tangential component of the velocity, linked to $u$ and $v$ by
\be
\label{eq.PolarVelocities}
	u =  \cphi \urho  - \sphi  \uphi,  \qquad   
	v =  \sphi  \urho + \cphi \uphi \ . 
\ee
Now by substituting into \eqref{eq.shallowWater} the expressions given in \eqref{eq.PolarVelocities} and \eqref{eq.changeOfDer}, after some calculations, we derive a new set of hyperbolic equations
%\be
%\begin{cases}
%&\de{r h}{t}  +
%\de{r h \urho}{r} + \de{h\uphi}{\phi}  = 0, \\[3pt]
%%%%%
%&\de{r h \urho}{t}  +  \de{}{r}\left ( r h\urho^2 + \frac{1}{2}g r h^2 \right) 
%+ \de{h\urho \uphi }{\phi} =  h\uphi^2  + \frac{1}{2}gh^2, \\[3pt]
%%%%%
%&\de{r h\uphi}{t} + \de{ r h\urho\uphi}{r} 
%+ \de{}{\phi} \left (  h\uphi^2 + \frac{1}{2} gh^2 \right )  = - h \urho \uphi,
%\end{cases}
%\ee
which, however, does not yet fit into the general form \eqref{eq.generalform}, since the fluxes in the above system depend explicitly on the spatial coordinate $r$ 
(i.e. the system is not autonomous).  Thus, we add the trivial equation $\partial r / \partial t = 0$ to the system,   
%\be
%\label{eq.trivialEqRadius}
%\de{r}{t} = 0,
%\ee
obtaining finally  
\begin{eqnarray}
%\begin{cases}
&&\de{r h}{t}  +
\de{r h \urho}{r} + \de{h\uphi}{\phi}  = 0, \\[3pt]
%%%%
&&\de{r h \urho}{t}  +  \de{}{r}\left ( r h\urho^2 + \frac{1}{2}g r h^2 \right) 
+ \de{h\urho \uphi }{\phi} =  h\uphi^2  + \frac{1}{2}gh^2, \\[3pt]
%%%%
&&\de{r h\uphi}{t} + \de{ r h\urho\uphi}{r} 
+ \de{}{\phi} \left (  h\uphi^2 + \frac{1}{2} gh^2 \right )  = - h \urho \uphi, \\[3pt]
&& \de{r}{t} = 0.
%\end{cases}
\label{eqn.swepolar}
\end{eqnarray}
The vector of the conserved variables, the non linear flux, and the source can now be written as 
\be
\label{eq.ShallowWaterCylindrical}
\Q  = \left( \begin{array}{c} r h\\ r h \urho \\ r h \uphi  \\ r  \end{array} \right), \quad \f = \left( \begin{array}{c} r h \urho \\ r h \urho^2 + \frac{1}{2} g r h^2 \\  r h \urho \uphi  \\ 0 \end{array} \right), \quad \g = \left( \begin{array}{c} h \uphi \\ h \urho \uphi \\ h \uphi^2 + \frac{1}{2}gh^2  \\ 0 \end{array} \right), \quad \mathbf{S} = \left( \begin{array}{c} 0 \\ h\uphi^2+ \frac{1}{2}gh^2 \\ -h\urho \uphi \\ 0  \end{array} \right).
\ee
and the Jacobian matrices, necessary for the computation of the ALE Jacobian matrix in \eqref{eq.ALEjacobianMatrix}, are
\be
\mathbf{A}_1 = \de{\textbf{f}}{\Q} =
\begin{pmatrix}
	0 & 1 & 0 & 0 \\ 
	- \urho^2 + gh & 2 \urho & 0 & - \frac{1}{2} g h^2 \\ 
	- \urho \uphi & \uphi  & \urho  & 0 \\ 
	0 & 0 & 0 & 0
\end{pmatrix}, 
\quad 
\mathbf{A}_2 = \de{\textbf{g}}{\Q} = \begin{pmatrix}
0 & 0  & \frac{1}{r} & - \frac{h \uphi}{r} \\ 
- \frac{\urho \uphi}{r} &  \frac{\uphi}{r} & \frac{\urho}{r} & - \frac{h \urho \uphi}{r} \\ 
- \frac{\uphi^2}{r} + g \frac{h}{r} & 0 & \frac{2 \uphi}{r} & - \frac{h \uphi^2}{r} - g \frac{h^2}{r} \\ 
0 & 0 & 0 & 0
\end{pmatrix}.
\ee

We choose the following initial condition
\begin{equation} 
\label{eq.InitalData_Equilibrium_cylindrical}
	 h(r, \phi,0) = 1 - \frac{1}{2g} \textup{e}^{-(r^2-1)},  \qquad 
	\urho(r, \phi,0) = 0, \qquad 
	\uphi(r, \phi,0) = r \textup{e}^{- \frac{1}{2} (r^2-1)}, 
\end{equation}
which is a stationary solution of \eqref{eqn.swepolar}, and so coincides with the exact solution at any time. We performed our test both with the Osher and the Rusanov fluxes and with a mesh made of triangles, quadrilaterals or both. 

\begin{figure}
	\centering 
	\includegraphics[width=0.8\textwidth]{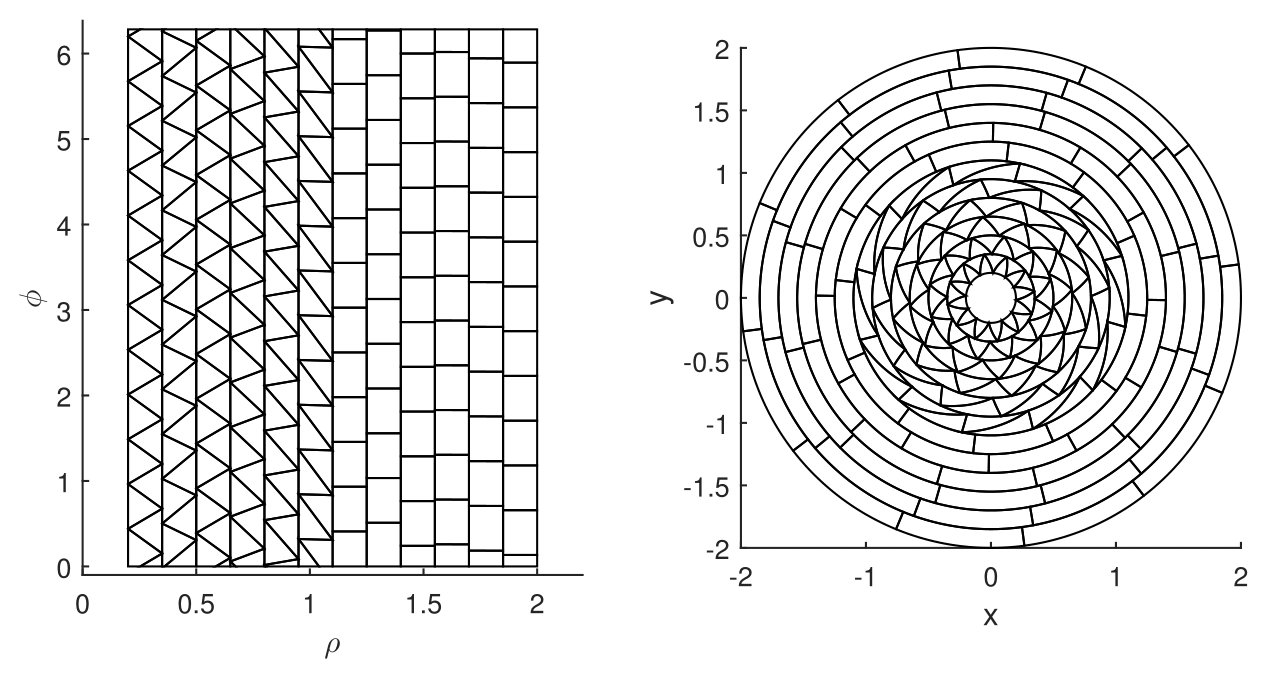} 
	\caption{Domain discretization at time $t=15$. On the left we report the grid in polar coordinates where the shear discontinuities lie over straight lines. On the right the corresponding grid in  Cartesian coordinates.} 
	\label{fig.polarToCart}
\end{figure}

The considered computational domain is $\Omega(r, \phi)= [0.2,2] \times [0, 2\pi]$ which is easily mapped to the annulus with radius $[0.2,2]$. Indeed the choice of considering the shallow water equations in polar coordinates allows us to study the vortex over a rectangular domain with periodic boundary conditions (at $\phi = 0$ and $\phi = 2\pi$) instead of dealing with circles. At $r = 0.2$ and $r = 2$ we have imposed reflective boundary conditions.
In particular using the polar coordinates the detected shear interfaces lie over straight lines and so they are perfectly handled by our algorithm. The images presented in this section  are then obtained by mapping back our results to Cartesian coordinates, as shown in Figure \ref{fig.polarToCart}.

\begin{figure}
	\centering 
	\includegraphics[width=0.8\textwidth]{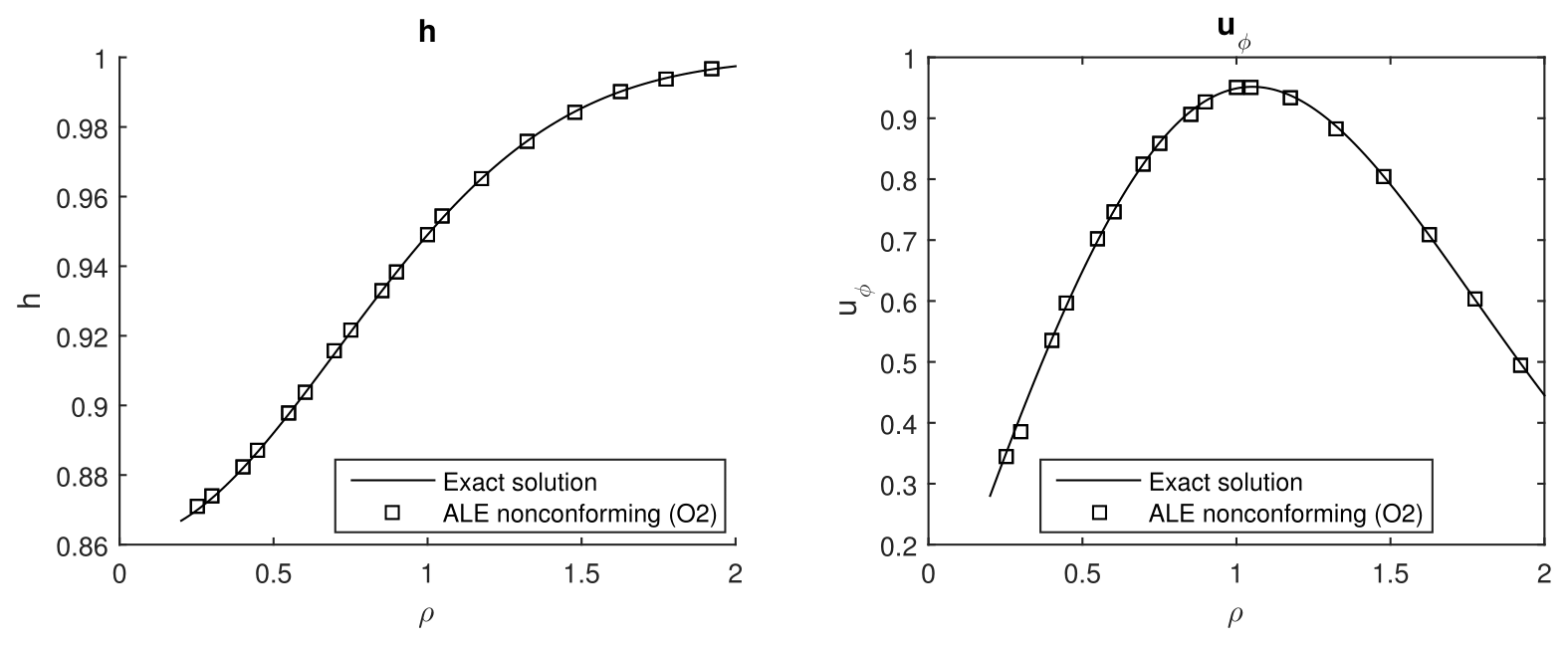} 
	\caption{Comparison between analytical solution and second-order accurate numerical results for the water level $h$ (left) and the tangential component of the velocity $\uphi$ (right), with $\phi=2\pi$ and $r  \in [0.2,2]$.}
	\label{fig.profile}
\end{figure}

First, Table \ref{tab.orderOfconvergence.vort} confirms the designed order of convergence of our algorithm in multiple situations: so primarily we can say that the mesh motion does not affect 
the standard order of convergence of the MUSCL-Hancock strategy and moreover this shows once again that the numerical flux computation, even at the nonconforming interfaces, is carried out
correctly. The numerical solution at $t=15$ is compared with the analytical one in Figure \ref{fig.profile}. 

Then we compare the results with a standard conforming method. First, let us underline that when the velocity changes even within the same element the only way to overcome the mesh distortion would be to split the element itself. For this reason, where the velocity field changes smoothly and as a consequence the shear flow affects all the vertices of the same element, at 
a certain time the mesh will become invalid even in the nonconforming case. This would not happen if the velocity field were uniform within each element, i.e. if each element moved all its 
vertices with the same velocity, e.g. the velocity of the barycenter. 
%In particular in this test the algorithm breaks when one edge of one element is so stretched that its length exceeds the half of the period. The break of the algorithm could be prevented simply by a %rezoning from time to time.

The main difference between the new nonconforming algorithm and a conventional conforming method is the final time at which the computation stops due to an invalid mesh, and the time step restriction that depends on the smallest encircle diameter of the elements. 

In Table \ref{tab.Ntimestep} we report the employed number of time steps and their dimension for different kinds of meshes and at different times. We remark that a larger 
value of $\Delta t$ decreases the required number of time steps and in this way also the total amount of computational time. The last results of each group refer to the 
moment at which the algorithm breaks due to an invalid mesh: one can easily see that the nonconforming method is able to run almost eight times longer than a conventional ALE 
method on conforming grids. 

\begin{figure}
\begin{center} 
\begin{tabular}{cc} 
	\includegraphics[width=0.37\textwidth]{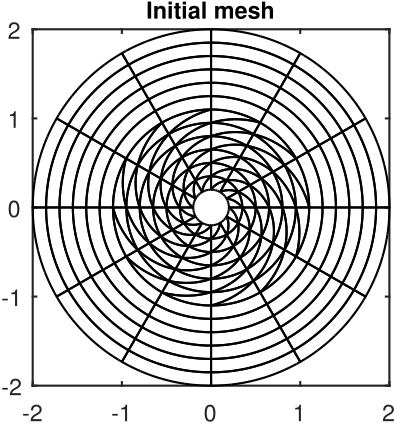} & 
	\includegraphics[width=0.37\textwidth]{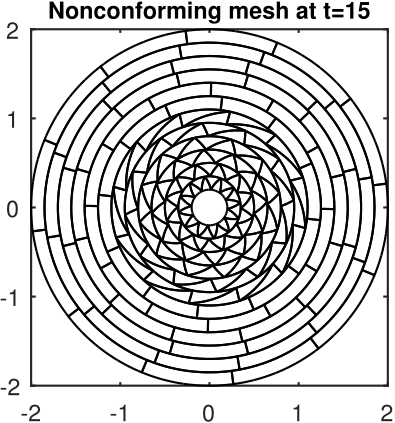} \\ 
	\includegraphics[width=0.37\textwidth]{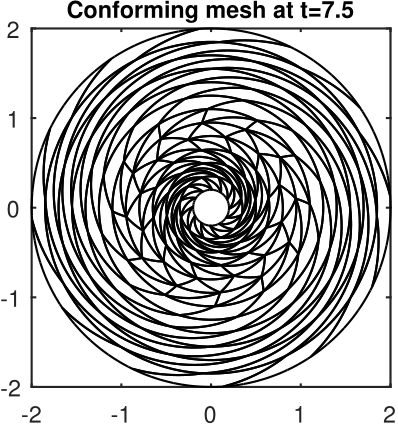} & 
	\includegraphics[width=0.37\textwidth]{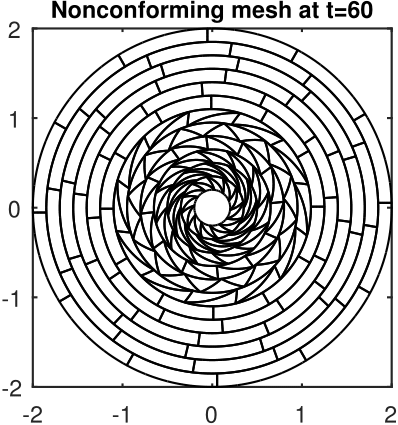}\\
	\includegraphics[width=0.37\textwidth]{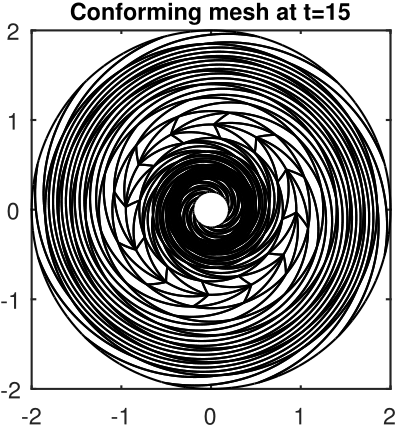} &
	\includegraphics[width=0.37\textwidth]{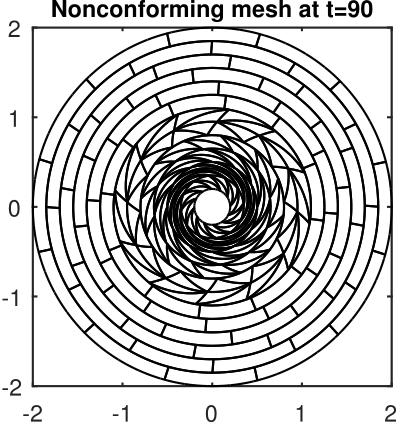} 
\end{tabular} 	
\end{center} 						
	\caption{Steady vortex in equilibrium. We compared the behavior of a standard conforming algorithm (without any rezoning technique) and of our new nonconforming method. Using the conforming algorithm the elements are deformed in a very short time, the time step is heavily reduced and hence the computation is slower. On the contrary, the nonconforming slide lines introduced by our scheme are able to  maintain a good shape of each element and an almost constant time step for a long computational time. Indeed only at time $t=90$ some elements with $r \rightarrow 0$ are deformed because of the presence of shear \textit{inside} the elements, which could be remedied only by subdividing the elements themselves
	or by removing them.} 
	\label{fig.Vortex}
\end{figure}

Finally, looking at Figure \ref{fig.Vortex} one can appreciate that the conforming method destroys the mesh immediately and then breaks, whereas the new nonconforming algorithm 
maintains a high quality mesh for a very long time, even with a very coarse mesh.

%\begin{figure}
%	\centering 
%	\includegraphics[width=0.81\textwidth]{vv7a} 
%	\caption{Steady vortex in equilibrium. We compared the behavior of a standard conforming algorithm (without any rezoning technique) and of our new nonconforming method. Using the conforming algorithm the elements are deformed in a very short time, the time step is heavily reduced and so the computation is slower; moreover the encircled diameter of some elements at $t=15$ surpasses the half of the period, fact that leads to the immediate crash of the periodic algorithm. On the contrary the nonconforming slide lines introduced by our scheme are able to maintain a good shape of each element and an almost constant time step for a long computational time. Indeed only at time $t=90$ some elements with $r \rightarrow 0$ are deformed because of the internal difference of velocity which could be prevented only subdividing the elements.}
%	\label{fig.Vortex}
%\end{figure}

% Immaginine dei vortici 2x3 invece che 3x2 come sopra
%\begin{figure}
%	\centering 
%	\includegraphics[width=\textwidth]{vv6a} 
%	\caption{Steady vortex in equilibrium. We compared the behavior of a standard conforming algorithm (without any rezoning technique) and of our new nonconforming method. Using the conforming algorithm the elements are deformed in a very short time, the time step is heavily reduced and so the computation is slower; moreover the encircled diameter of some elements at $t=15$ surpasses the half of the period, fact that leads to the immediate crash of the periodic algorithm. On the contrary the nonconforming slide lines introduced by our scheme are able to maintain a good shape of each element and an almost constant time step for a long computational time. Indeed only at time $t=90$ some elements with $r \rightarrow 0$ are deformed because of the internal difference of velocity which could be prevented only subdividing the elements.}
%	\label{fig.Vortex}
%\end{figure}

\begin{table} 
	\begin{center} 	
		\begin{tabular}{cccc|cccc} 
			\hline 
			 	\multicolumn{4}{c|}{$\mathcal{O}2$, Osher flux, quadrilateral elements}& \multicolumn{4}{c}{$\mathcal{O}2$, Rusanov flux, triangular elements}       \\ 
			\hline \multicolumn{8}{c}{} \\[-10pt]
			mesh points &  $h(\Omega(t_f))$ & $\epsilon_{L_2}$ & $\mathcal{O}(L_2)$ &  mesh points &  $h(\Omega(t_f))$ & $\epsilon_{L_2}$ & $\mathcal{O}(L_2)$  \\ \multicolumn{8}{c}{} \\[-10pt]
			\hline
			12 $\times$ 12 & 2.33E-01 & 1.36E-03 &	-	 &  20 $\times$ 20   & 7.18E-02 &   5.97E-04 &  -       \\
			24 $\times$ 24 & 1.17E-01 &	3.42E-04 &  1.99 &	30 $\times$ 30   & 5.21E-02 &	2.54E-04 &  2.11	  \\	
			32 $\times$ 32 & 8.74E-02 &	1.94E-04 & 	1.97 & 	40 $\times$ 40   & 3.91E-02 &	1.43E-04 &  2.01    \\	
			44 $\times$ 44 & 6.36E-02 &	1.03E-04 &	1.98 &	55 $\times$ 55   & 2.84E-02 &	7.76E-05 &	1.91    \\	
			60 $\times$ 60 & 4.66E-02 &	5.57E-05 &	1.99 &	60 $\times$ 60   & 2.60E-02 &	6.58E-05 &  1.91    \\	
			
			\hline 
		\end{tabular}		
	\end{center}
	\caption{Numerical convergence results for the steady vortex in equilibrium using nonconforming meshes. In the left table we report the results obtained on a quadrilateral mesh using the Osher type flux. For the results on the right we have employed a triangular mesh and the Rusanov type flux. The error refers to the difference between the computed water level $h$ and the exact one at time $t_f = 0.5$.}
	\label{tab.orderOfconvergence.vort}
\end{table}

\begin{table}  
	\begin{center} 
		\begin{tabular}{c|ccc|ccc|ccc }
	\hline  	
	 $N_E$  $\rightarrow$ & \multicolumn{3}{c|}{216}  & \multicolumn{3}{c|}{264} & \multicolumn{3}{c}{300} \\ 
	\hline \multicolumn{10}{c}{} \\[-8pt]
	\multicolumn{10}{c}{conforming algorithm} \\
	\hline 
	& $t$ & $n$ & $\Delta t$ &  $t$ & $n$ & $\Delta t$ &  $t$ & $n$ & $\Delta t$   \\
	\hline
	& 1 & 110 & 9.58E-03 &  1 & 180 & 5.40E-03 &  1 &  180& 5.71E-03   \\	
	&  8&  1163&  {4.13E-03}  & 8 & 2180 &  {2.52E-03} &  10 &  2071&  {3.11E-03}    \\	
 & 12 &  2370&  2.70E-03&  12 & 4035 & 1.89E-03 &  15 &  4098& 2.04E-03  \\
	stop at  $\rightarrow$ & 15.3 &  \textbf{3773} & 2.06E-03 &  15.5 & \textbf{6072} & 1.54E-03 &  17 & 5190 &  1.78E-03   \\
	\hline \multicolumn{10}{c}{} \\[-8pt]
\multicolumn{10}{c}{nonconforming algorithm} \\
\hline 
& 1    & 110    & 9.58E-03 &  1 & 180 & 5.82E-03 &  1 & 175 & 5.68E-03   \\
 & 8    & 851    &  {9.50E-03} & 8 & 1410 &  {5.52E-03} &  10 & 1720 &  {5.92E-03}   \\
& 30   &  \textbf{3175}  & 9.30E-03 &  30 & \textbf{6033} & 4.06E-03 &  15 & 2565 & 5.94E-03   \\
&  60  & 7757   & 4.90E-03 &  60 & 15010 & 2.84E-03 &  80 & 15979 & 3.34E-03   \\
stop at $\rightarrow$ &  119 & 26430  & 2.24E-03 &  129 & 35791 & 1.94E-03 &  132 & 36275 & 2.13E-03   \\
	\hline 

\end{tabular} 
	\end{center}
	\caption{In this table we report the number of time steps $n$ necessary to reach the time $t$ and the  dimension of the time step $\Delta t$ at that time. We used three different meshes with $N_E$ total number of elements (triangles or quadrilaterals). The results are obtained by applying a standard conforming method and our new nonconforming algorithm. Looking at the bold data one can see that with almost the same number of time steps one reaches a simulation time that is twice as large with the nonconforming algorithm compared to a classical conforming one.		
    Besides the final simulation time that can be reached before obtaining an invalid mesh is almost 8 times larger. } 
	\label{tab.Ntimestep}
\end{table}

%\clearpage 

\section{Extension to general slide lines} 
\label{sec.general} 
\begin{figure}[!ht]
	\centering
	\includegraphics[width=0.85\textwidth]{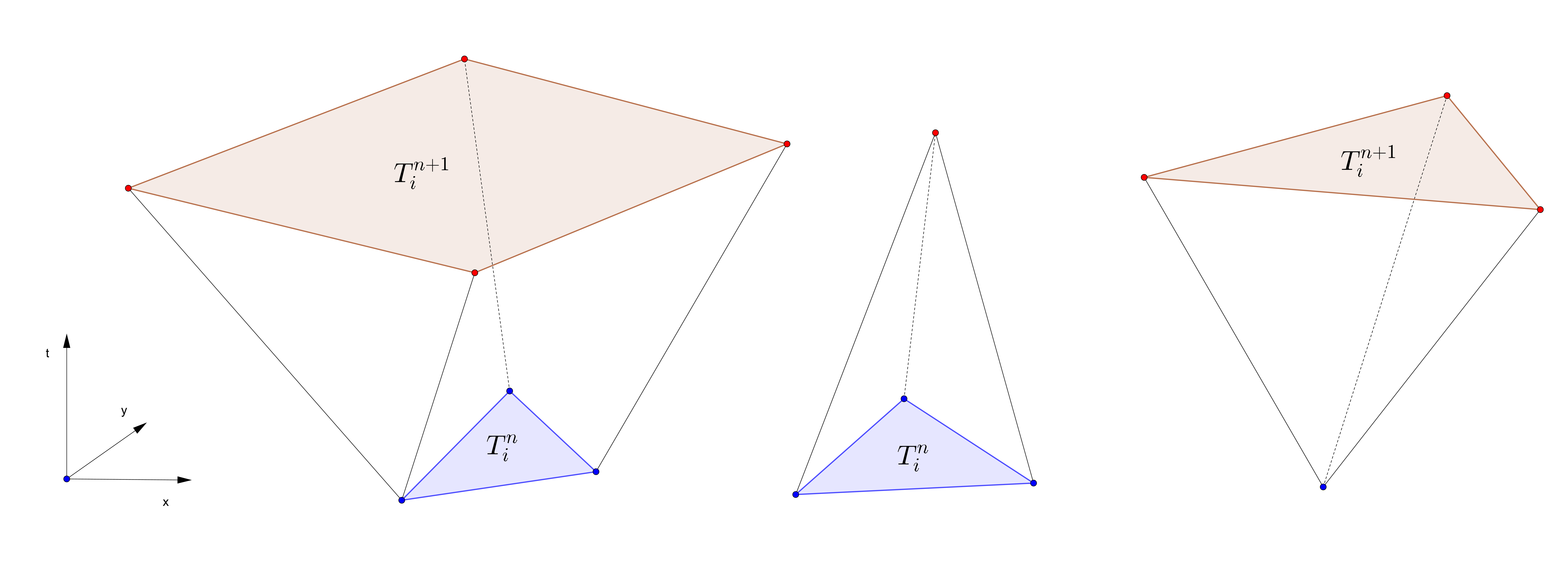}  
	\caption{Dynamic change of element type (left), element removal (center) and element insertion (right) between time $t^n$ and time $t^{n+1}$. 
	         Nodes and element $T_i^n$ at time $t^n$ are highlighted in blue, nodes and element $T_i^{n+1}$ at time $t^{n+1}$ are colored in red.   
	}
	\label{fig.addremove}
\end{figure}

\begin{figure}[!ht] 
	\centering 
	\begin{tabular}{cccc} 
	  \includegraphics*[width=0.22\textwidth]{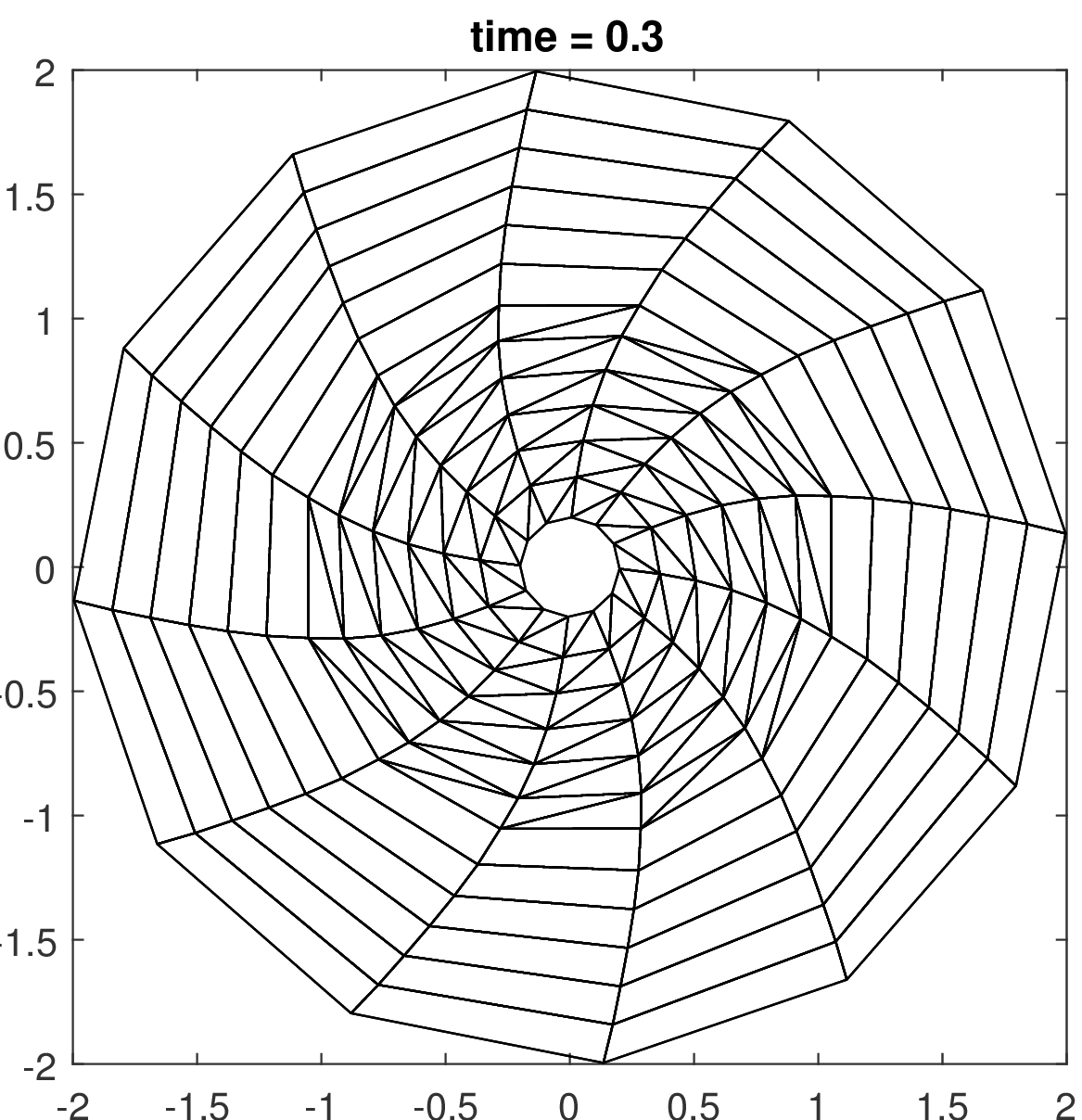} & 
    \includegraphics*[width=0.22\textwidth]{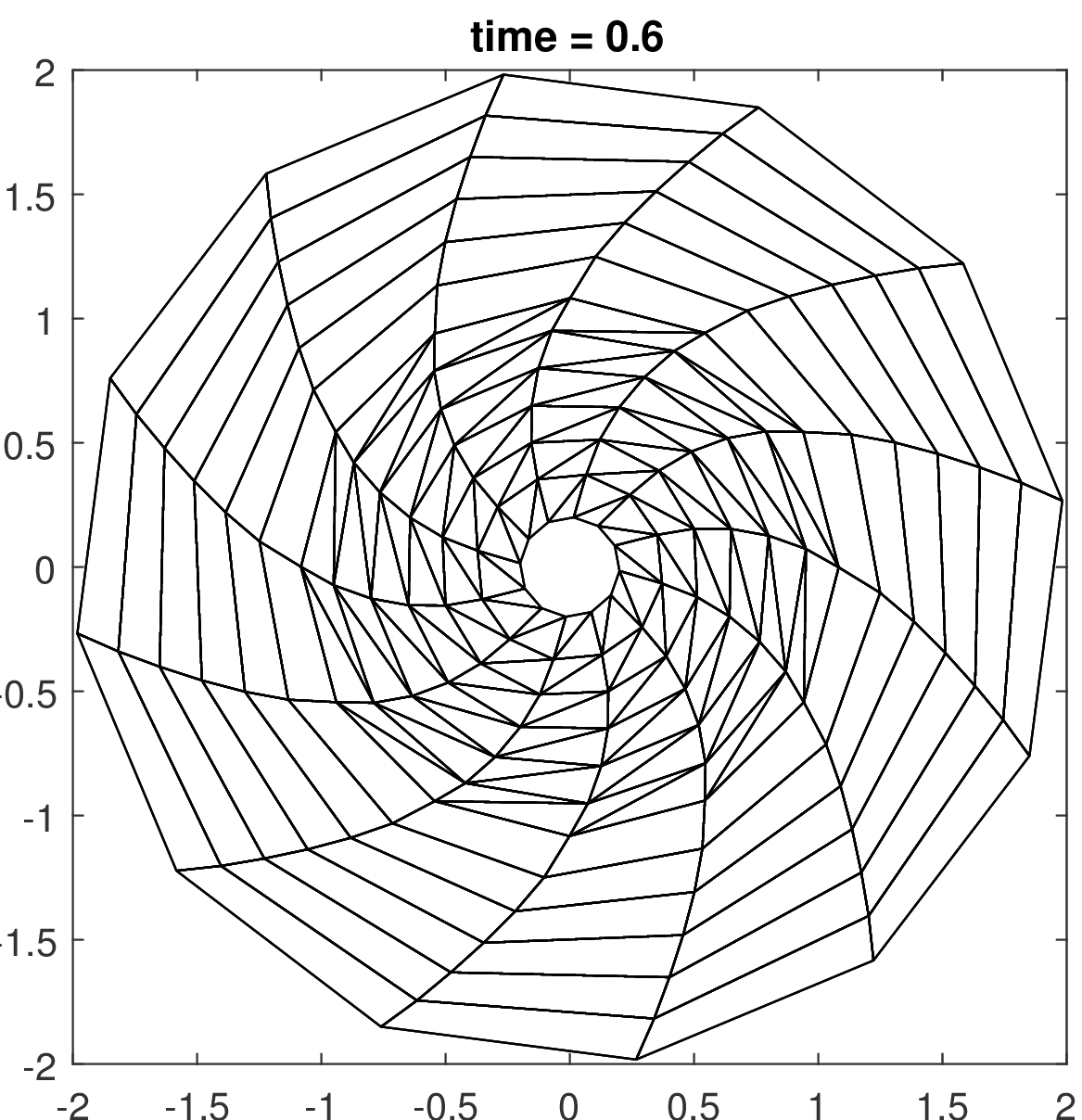} & 
    \includegraphics*[width=0.22\textwidth]{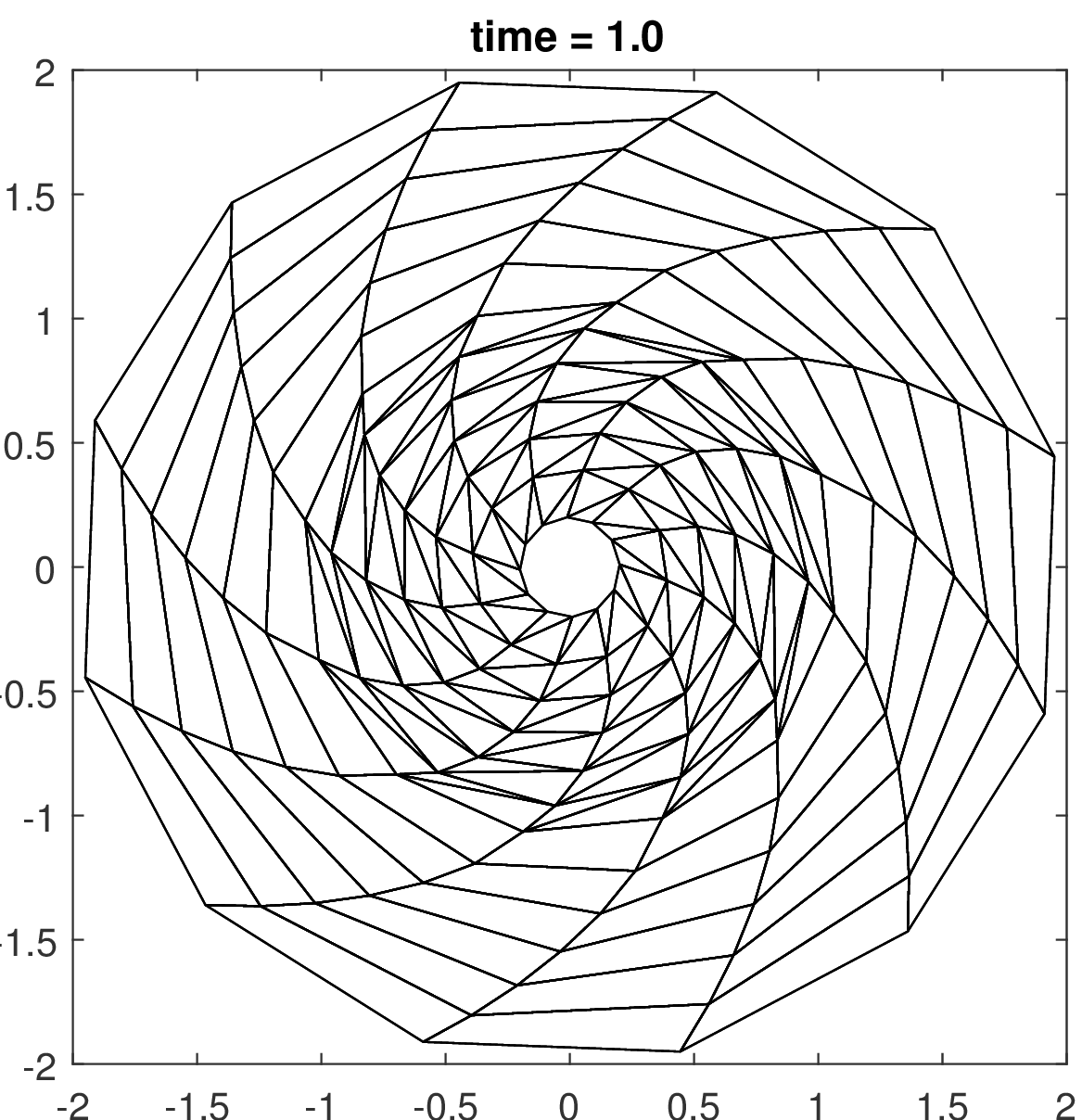} & 
    \includegraphics*[width=0.22\textwidth]{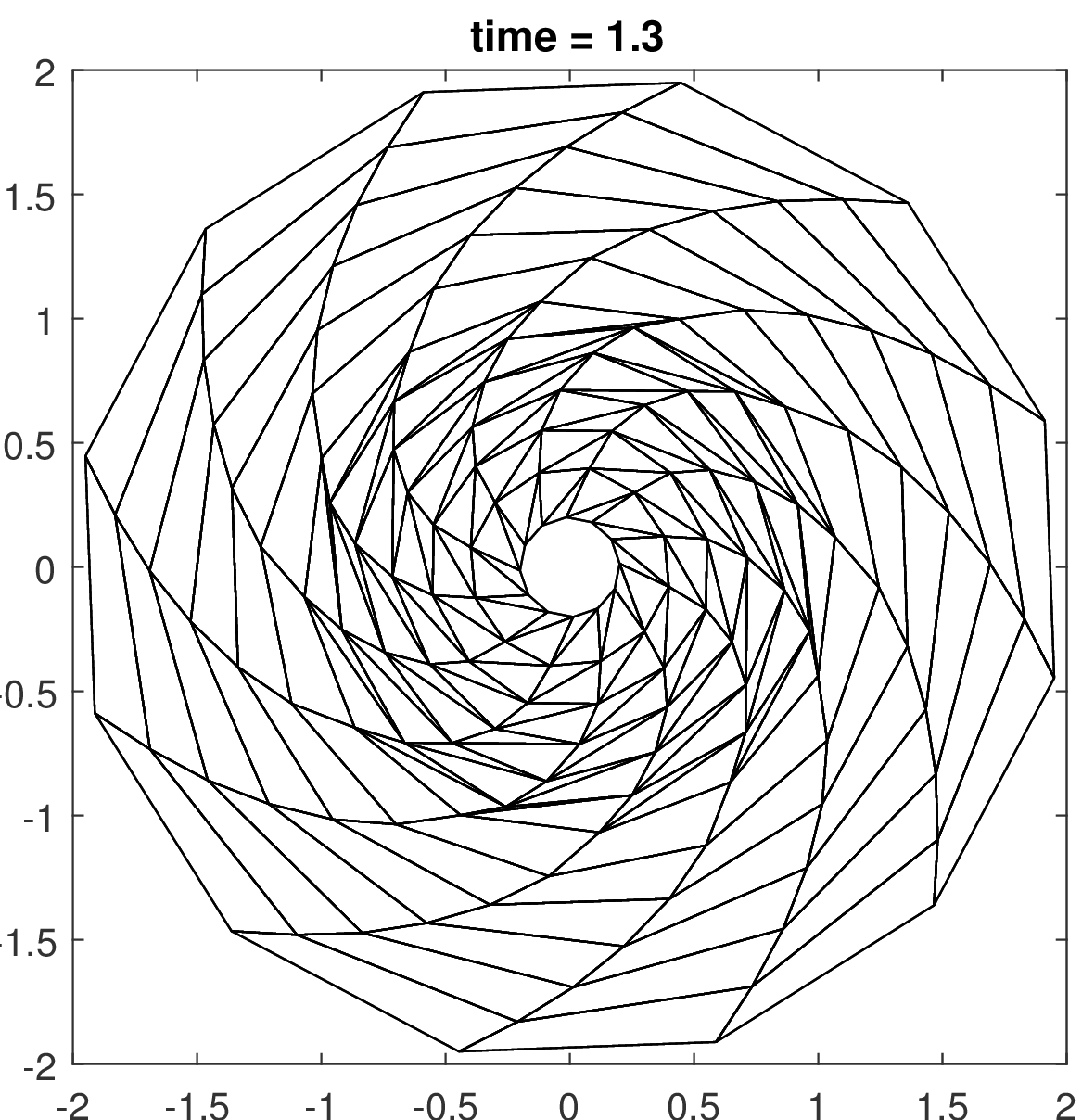} \\ 
 	  \includegraphics*[width=0.22\textwidth]{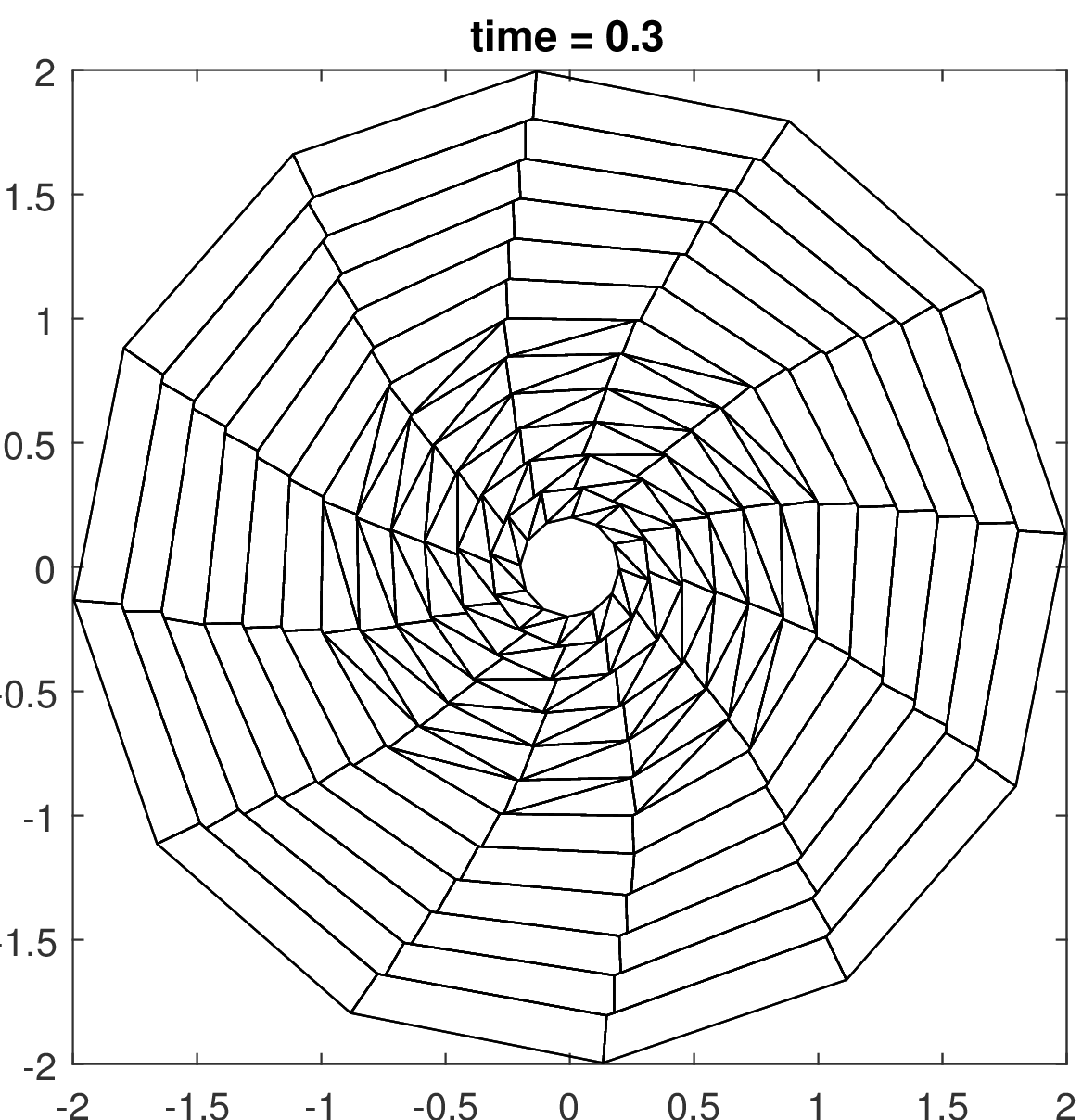} & 
 	  \includegraphics*[width=0.22\textwidth]{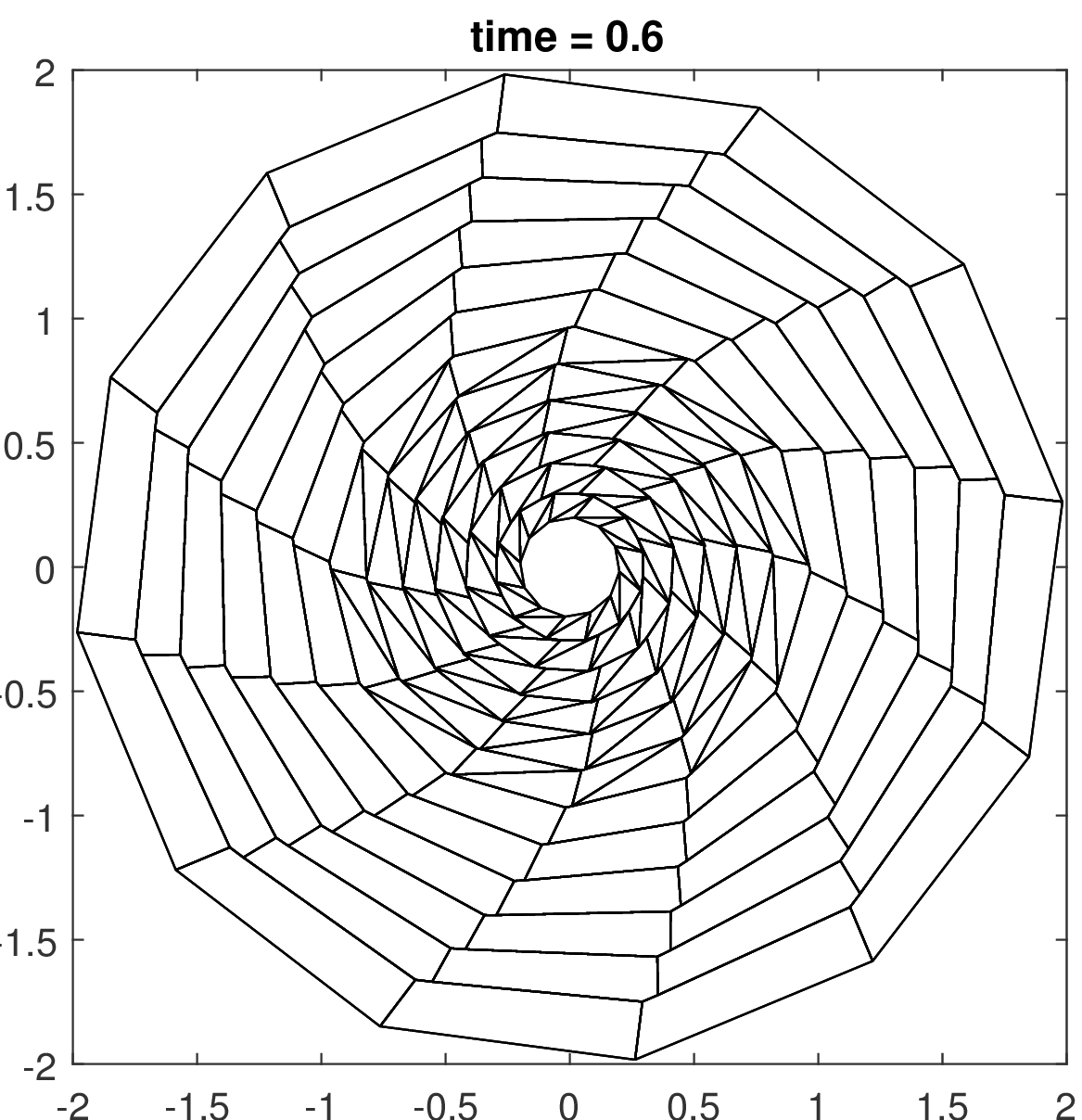} &  
	  \includegraphics*[width=0.22\textwidth]{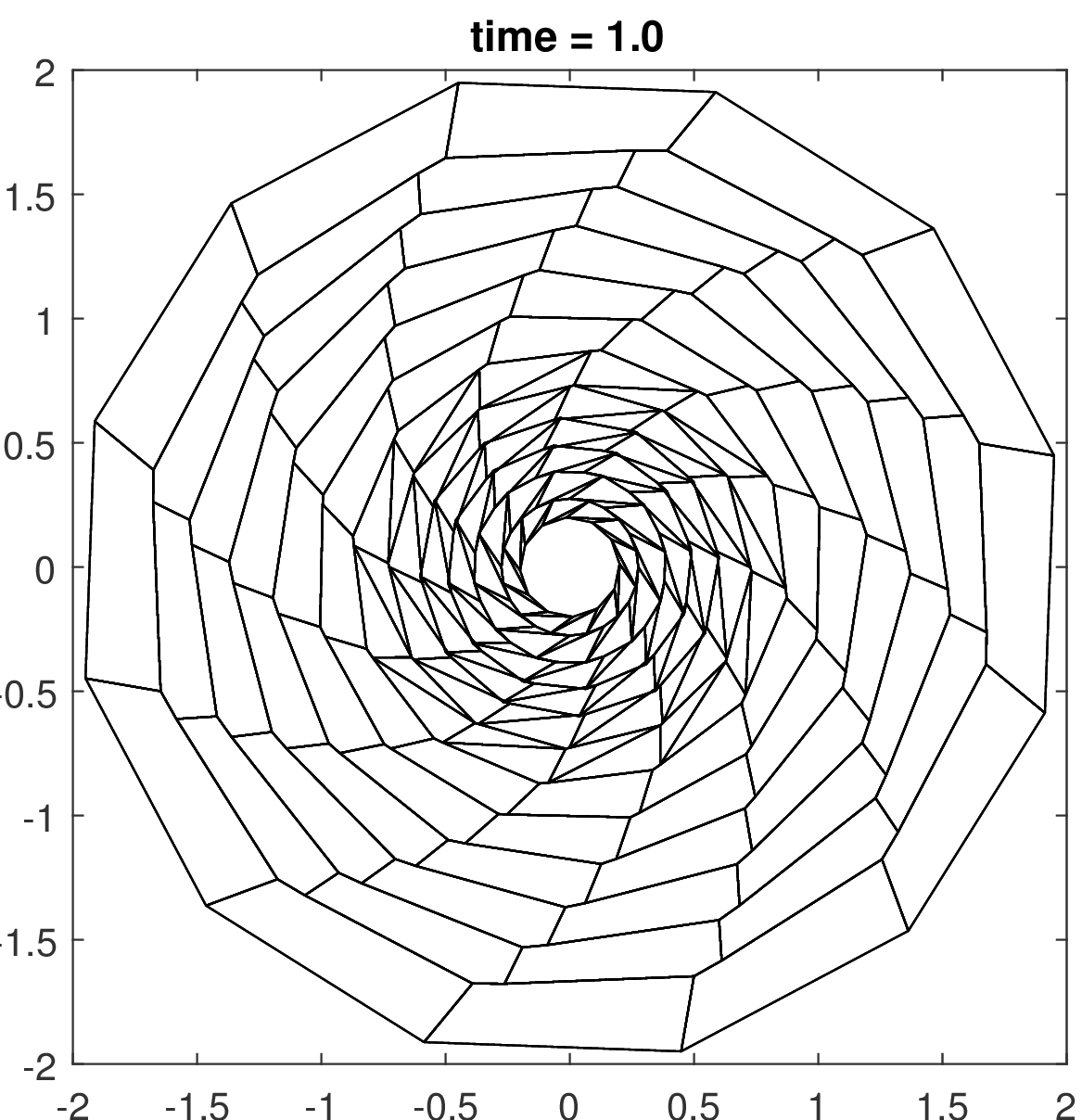} & 
	  \includegraphics*[width=0.22\textwidth]{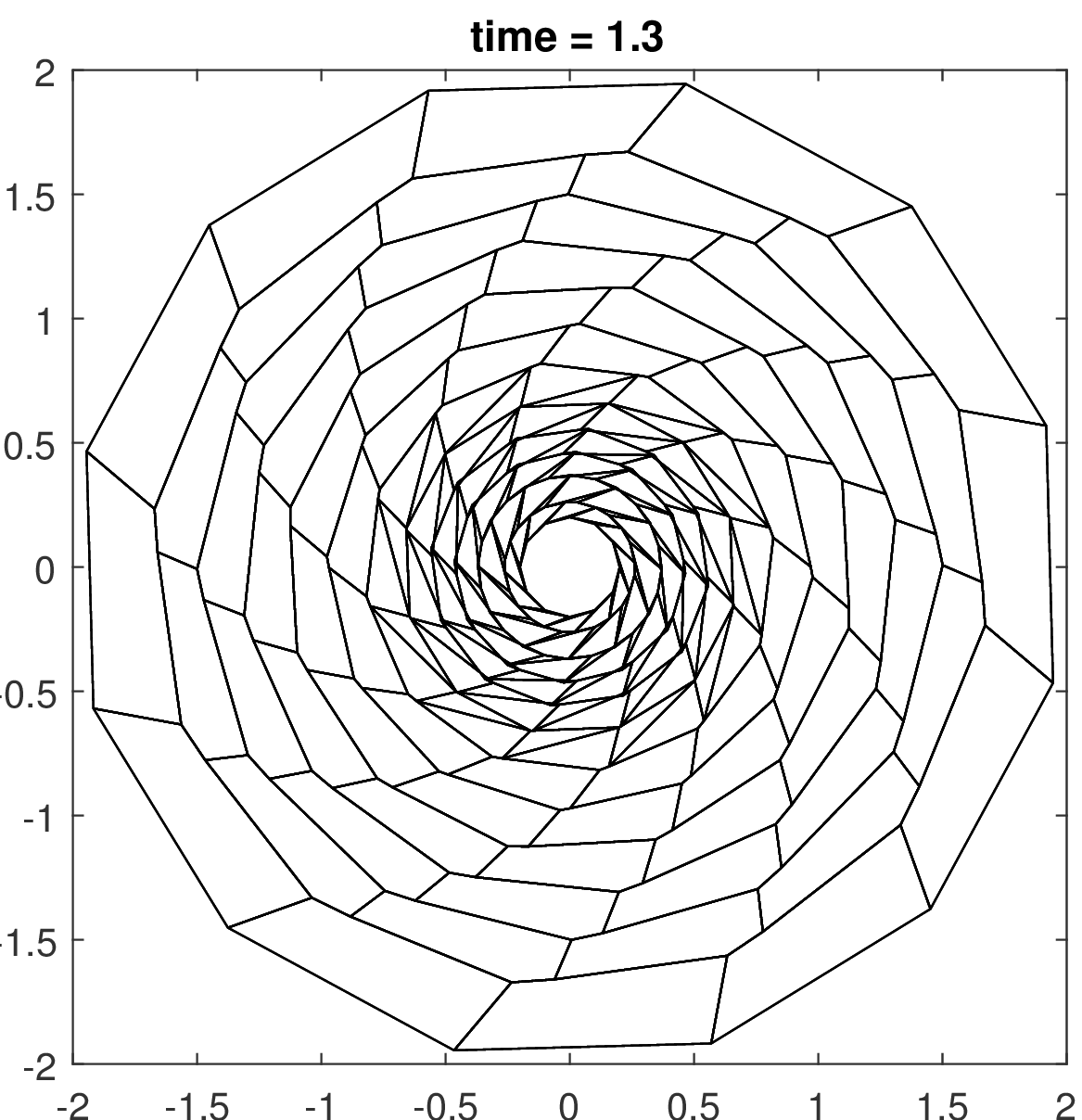} \\ 
	  \includegraphics*[width=0.22\textwidth]{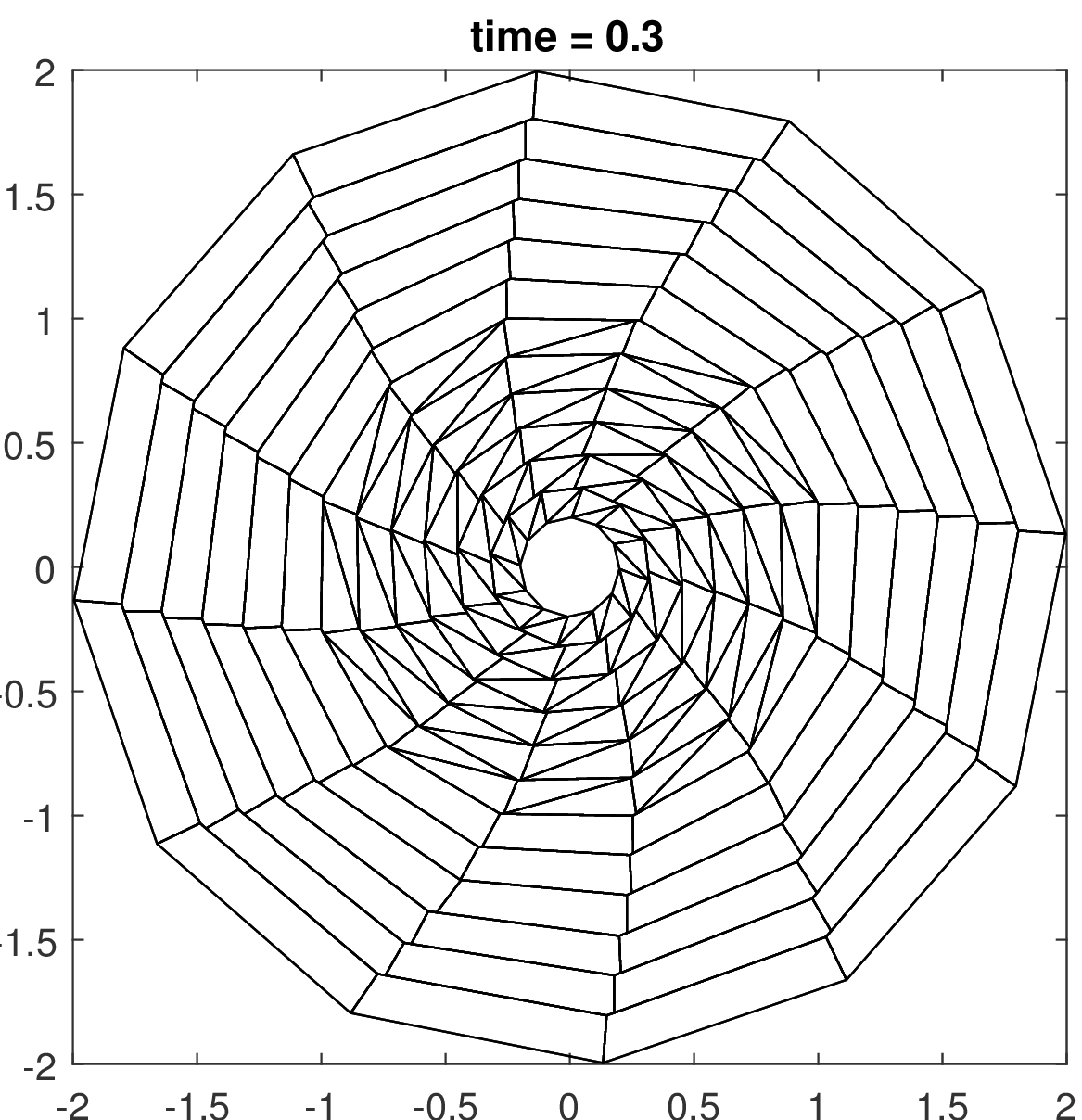} &
	  \includegraphics*[width=0.22\textwidth]{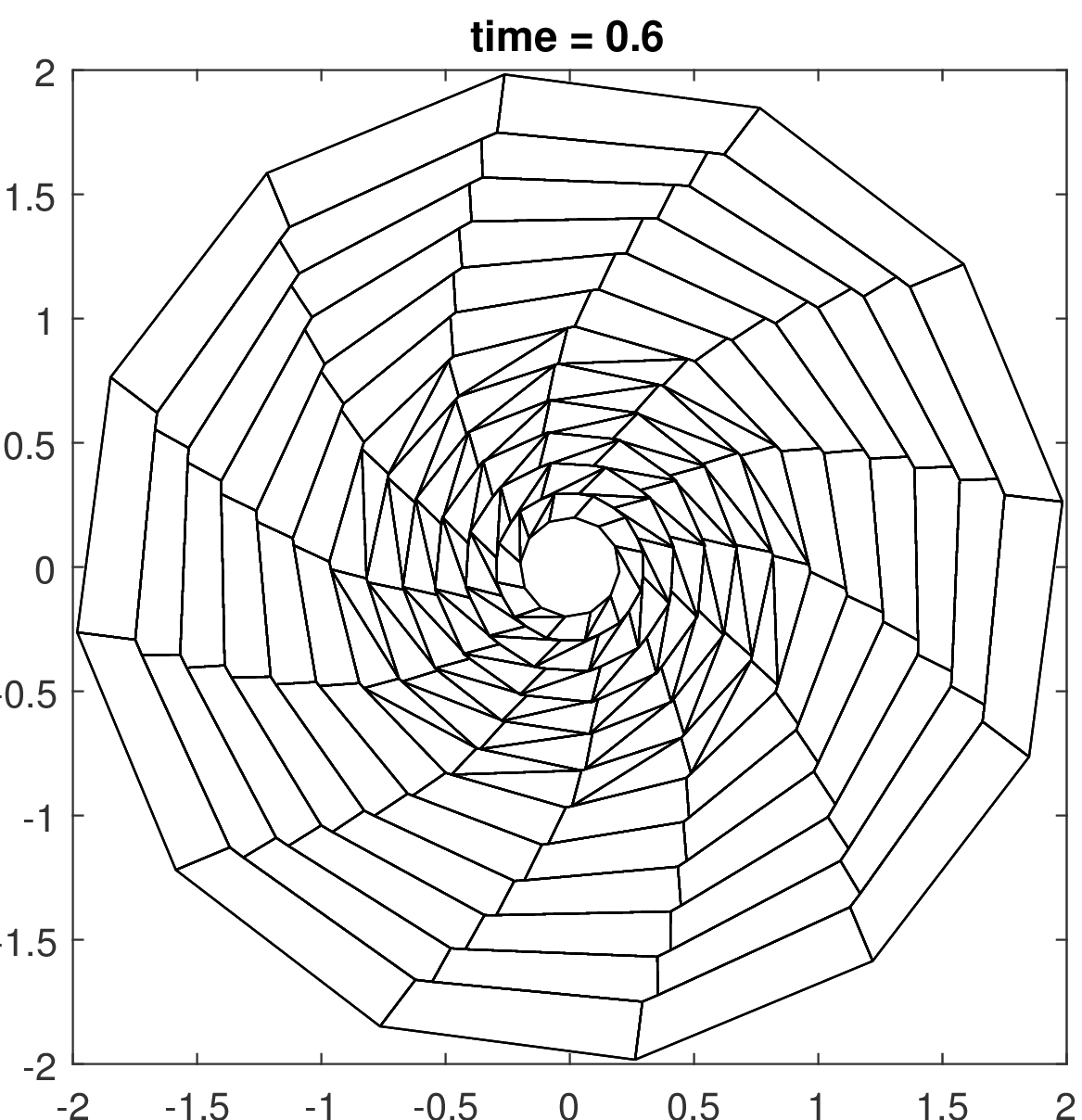} & 
	  \includegraphics*[width=0.22\textwidth]{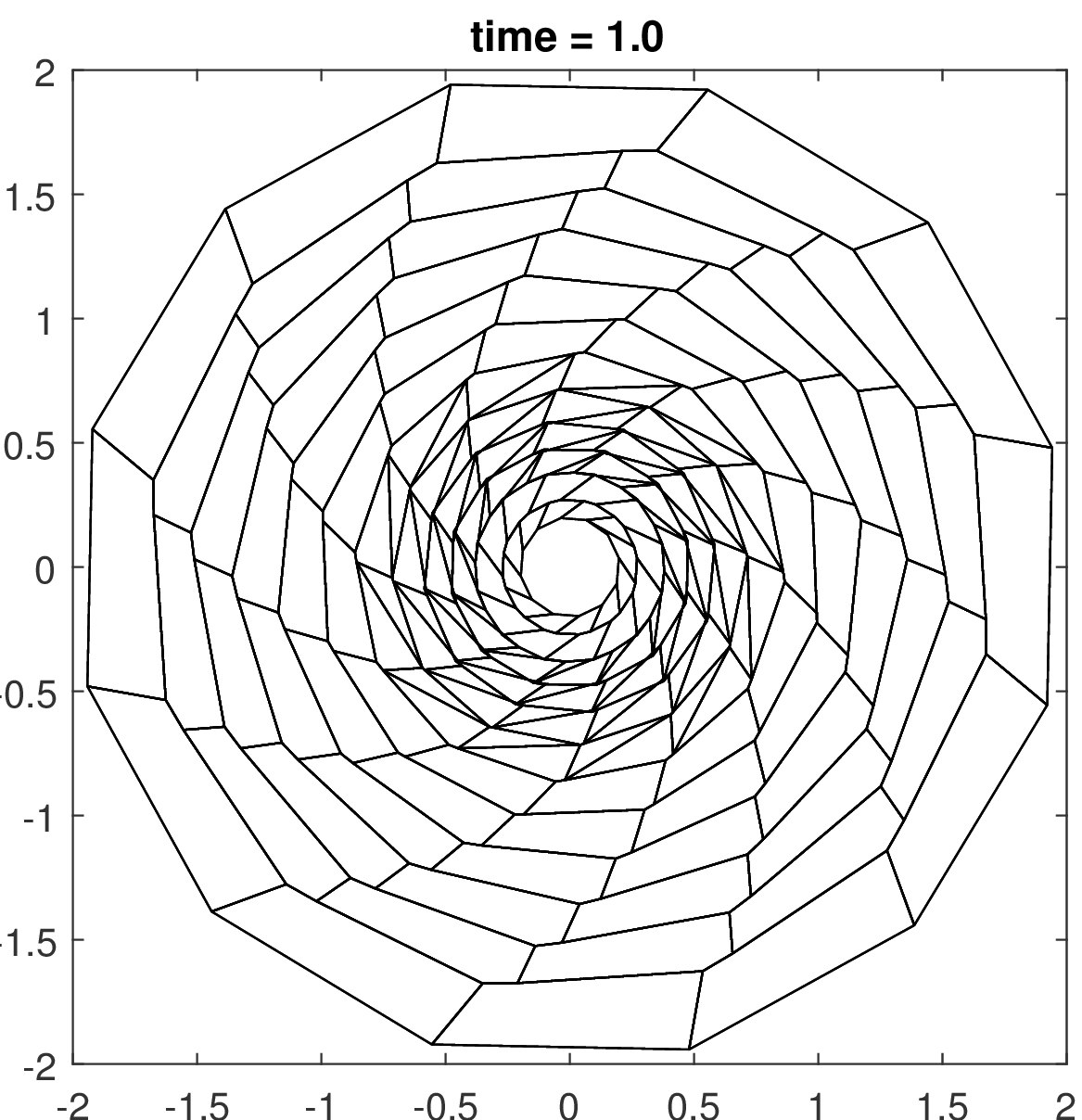} &
	  \includegraphics*[width=0.22\textwidth]{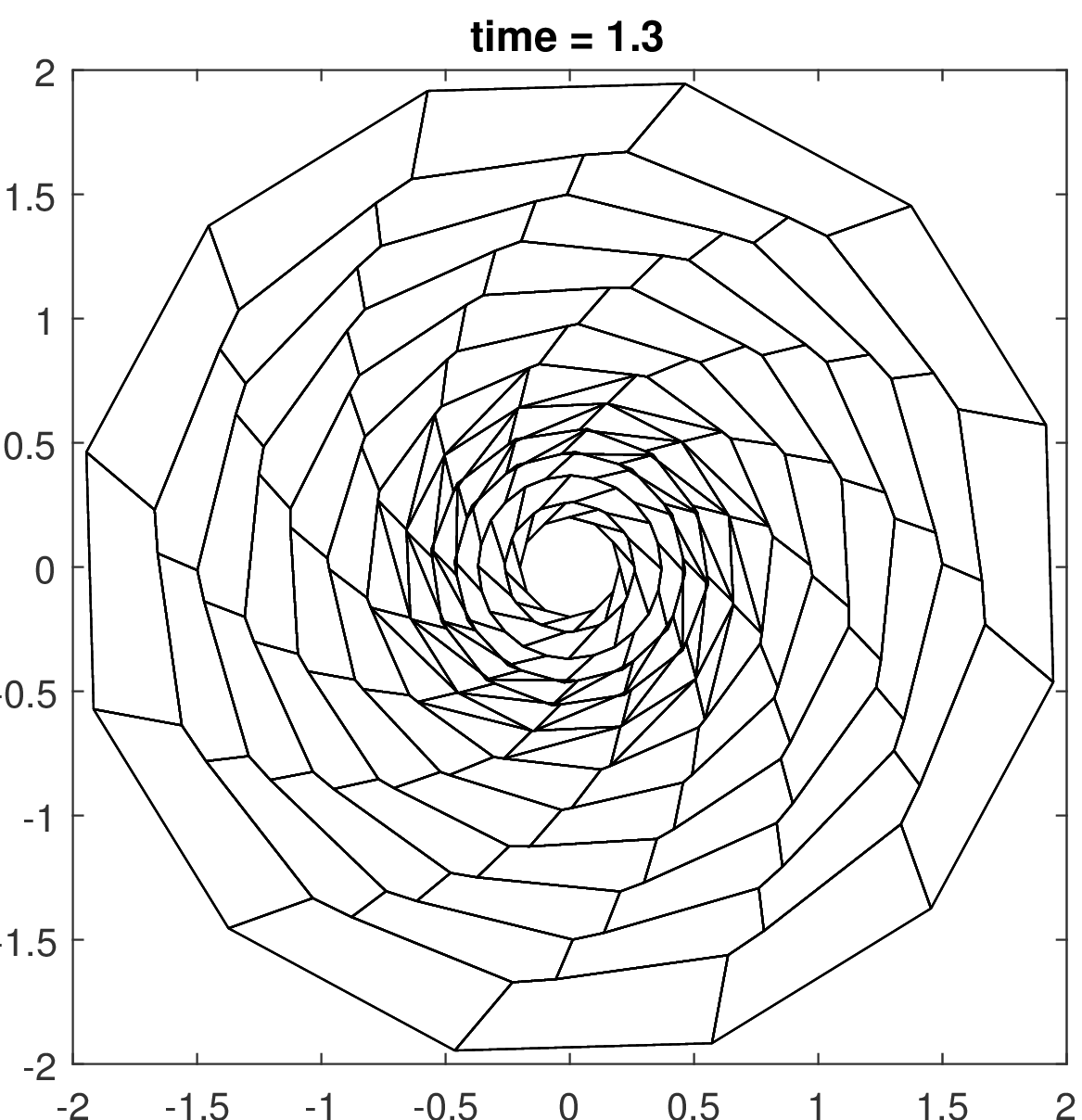} 
		\end{tabular} 
	\caption{ Isolated vortex in Cartesian coordinates. Classical conforming algorithm without any rezoning technique (top). 	          
           	Moving mesh obtained with the new \textit{nonconforming} algorithm at different times (center) without small element removal. 
						Moving \textit{nonconforming} mesh with small element removal (bottom), which allows to control the time step size and to maintain a better mesh quality. 
						The nonconforming algorithms used here use logically \textit{non-straight} slide lines. The sliding edges are automatically detected based on the tangential velocity difference.} 
	\label{fig.NonConf_pw}
\end{figure}

All test problems shown before were limited to logically straight slide lines. However, there is no intrinsic limitation to logically straight slide lines in our algorithm, since the integral 
space-time conservation form \eqref{eq.I1} of the conservation law is valid for \textit{arbitrary} closed space-time control volumes. This simple, elegant but at the same time very powerful 
formulation allows also to 
dynamically \textit{add} and \textit{remove} elements or to change their type during the simulation in a consistent manner that respects the GCL as well as local and global conservation. All these  features are trivially built in \textit{by construction}, due to the integral formulation on closed space-time control volumes. In Figure \ref{fig.addremove} we show examples of space-time  
control volumes that result when elements change type or when elements are dynamically added and removed during a simulation. For logically non-straight slide lines, it is necessary to divide 
elements and nodes into masters and slaves, where the master elements maintain their number of nodes, while the slave elements must in general change their element type during the sliding process. 
Also note that master nodes are free to move anywhere, while slave nodes must slide along the master edges. Furthermore, small elements need to be removed if they lead to  excessively small time 
steps due to the CFL condition. We now repeat the same shallow water vortex test problem as described in the previous section, but using the PDE in Cartesian coordinates.  
This leads to logically non-straight slide lines. The comparison between the classical conforming moving mesh algorithm and the new nonconforming approach presented in this paper is shown in 
Fig. \ref{fig.NonConf_pw} and Table \ref{tab.DeltaTime_comparison}. 
We observe the improved mesh quality and time step size compared to the classical conforming approach, in particular when the moving nonconforming mesh is combined with the removal of small 
elements. The obtained results look promising and justify further research in this direction in the future. 

%\begin{figure}
	%\centering 
	%\begin{tabular}{cc} 
	%\includegraphics[width=0.4\textwidth]{InitialMesh.eps}            &  
	%\includegraphics[width=0.4\textwidth]{WithKillingElement_t1-7}  
	%\end{tabular} 
	%\caption{Initial conforming mesh at time $t=0$ (left) and final mesh at time $t=1.4$ (right) obtained with the new nonconforming algorithm with small element removal. 
	         %The classical conforming mesh algorithm blows up at time $t<1.4$. } 
	%\label{fig.ithKillingElement}
%\end{figure}

\begin{table}[!hb]  
	\begin{center} 	
		\begin{tabular}{c|c|c|c} 
			\hline 
			& \multicolumn{3}{c}{Time step size}     \\ 
			\hline 
			time & conforming & nonconforming & nonconforming + \\
           &            &               & element removal \\  			
			\hline 
			0.3 & 3.8E-3  & 3.2E-3 & 3.2E-2	   \\
			0.6 & 3.6E-3 & 2.1E-3 &	 2.1E-3  \\
			1.0 & 1.9E-3 & 9.0E-4 &	 1.2E-3 \\
			1.3 & 5.8E-4 & 1.2E-4 &	 1.4E-3  \\
			1.7 & - & - & 1.4E-3 \\
			\hline 
		\end{tabular}		
	\end{center}
	\caption{Time step size for three different moving mesh algorithms. The main improvement is achieved when using a nonconforming algorithm combined with small element removal. 
	         This allows to maintain reasonable timesteps also for longer simulation times. }
	\label{tab.DeltaTime_comparison}
\end{table}

%\clearpage

\section{Well balancing for the shallow water equations in polar coordinates}
\label{sec.WB}

At this point we are able to maintain a high quality mesh even in the case of strong shear flows, and to preserve the physical properties of the system (mass, momentum, energy) for very long computational times. Moreover, thanks to the use of a Lagrangian framework, our novel method is little dissipative for contact discontinuities.   

The aim of this last section is to extend the algorithm in such a way that in addition it can preserve also exactly (i.e. up to machine precision) certain relevant and non-trivial equilibrium solutions. 
In particular, our interest is focused on the shallow water equations in cylindrical coordinates given by (\ref{eq.ShallowWaterCylindrical}), where the presence of a source term makes this task quite  challenging. 

The procedures that allow to preserve some equilibrium of interest are called \textit {well balanced} methods. We refer to \cite{Pares2006,Castro2006} and the references therein for a theoretical 
framework. For well-balanced schemes in the presence of gravity, see e.g. \cite{BottaKlein,Klingenberg2015,Kapelli2014}.  

The importance of these techniques is related to the fact that conventional numerical schemes in which the source term may be discretized in a consistent manner are not able to preserve 
certain stationary solutions of the PDE, especially on coarse meshes. This leads to erroneous numerical solutions when trying to compute small perturbations around the stationary solution 
necessitating the need for very fine meshes. 

For example, the initial data considered in the previous test case (\ref{eq.InitalData_Equilibrium_cylindrical}) represents a stationary solution for the system, characterized by the equilibrium 
between the centrifugal and the gravitational forces. But this precise balancing has to be achieved also at the numerical level in order to preserve the stationary solutions. 
And this is not trivial, especially because the two forces appear one in the flux and the other one in the source term, and consequently in classical finite volume schemes, they are discretized in a  different way, making it almost impossible to maintain a precise balancing between them. 
 
Following the ideas presented in \cite{Pares2006} and \cite{Castro2007}, we decide to rewrite the source terms in (\ref{eq.ShallowWaterCylindrical}) by means of nonconservative products in order to take them into account directly in the flux computation.
The general form of the obtained system is 
\be
\label{eq.generalformB}
\de{\Q}{t} + \nabla \cdot \F(\Q)  + \B(\Q) \cdot \nabla \Q = \mathbf{S}(\Q),
\ee
where, with respect to (\ref{eq.generalform}), we have the additional matrix $\B(\Q) = (\, \B_1(\Q), \B_2(\Q) \,) $ which collects the nonconservative terms written using the standard conserved variables and some other trivially conserved variables added to the system.
To clarify our notation, we want to remark that with $\Q = (q_1,q_2, \dots, q_{\nu})$ we mean the vector of conserved variables defined in the space of the admissible states $\Omega_{\Q} \subset \mathbb{R}^{\nu}$, and with $\q_h = (q_{h,1},q_{h,2}, \dots, q_{h,\nu})$ we mean the value of the conserved variables obtained through a reconstruction procedure inside each element. (For a first order method $\q_h$ in the control volume $C_i^n$ simply coincides with the cell average of element $T_i^n$).

System (\ref{eq.ShallowWaterCylindrical}) can be cast in form (\ref{eq.generalformB}) by adding as auxiliary variables the radius $r$ and the bottom topography $b$ such that the free 
surface is $\eta(r, \phi) = b+h(r, \phi)$. The involved terms are the following 
\be
\label{eq.ShallowWaterCylindricalWB}
\Q  \!=\! \left( \begin{array}{c} r h\\ r h \urho \\ r h \uphi  \\ r b \\ r                 \end{array} \right)\!, 
\ \f \!=\! \left( \begin{array}{c} r h \urho \\ r h \urho^2 \\  r h \urho \uphi  \\ 0 \\ 0    \end{array} \right)\!, 
\ \g \!=\! \left( \begin{array}{c} h \uphi \\ h \urho \uphi \\ h \uphi^2 + \frac{1}{2}gh^2  \\ 0 \\ 0  \end{array} \right)\!, 
\ \B_1 \cdot \nabla \Q \!=\! \left( \begin{array}{c}  0 \\ g r h \de{\eta}{r} - h \uphi^2 \de{r}{r} \\ h \urho \uphi \de{r}{r} \\ 0 \\ 0 \end{array} \right)\!, \ 
\B_2 \!=\! 0, \ \mathbf{S}\!=\!0.
\ee

The main difficulty of systems written in this form, both from the theoretical and the numerical point  of view, comes from the presence of nonconservative
products that do  not make sense in the distributional framework when the solution $\Q$ develops discontinuities.
From the theoretical point of view, in this paper we assume the definition of nonconservative products as Borel measures given in \cite*{DalMaso1995}. 
This definition, which depends on the choice of a family of paths in the phase space, allows one to give a rigorous definition of weak solutions of (\ref{eq.generalformB}).

We consider here the discretization of system (\ref{eq.generalformB}) by means of a numerical scheme which is \textit{path-conservative} in the sense introduced in \cite{Pares2006}, and that can be cast in our space--time formulation as follows 
\be
\label{eq.Scheme}
|T_i^{n+1}| \, \Q_i^{n+1} =  |T_i^n| \, \Q_i^n - \sum \limits_{j} \,\,\int_0^1 \!\! \int_0^1 
| \partial C_{ij}^n| \left ( \tilde{\F}_{ij} + \tilde{\mathbf{D}}_{ij} \right ) \cdot \mathbf{\tilde n}_{ij} d\chi d\tau 	- \int_{C_i^n \backslash \partial C^n_i}    \tilde{\B} \cdot \tilde{\nabla} \q_h
\ d\mathbf{x} dt + \int_{C_i^n}  \S d\mathbf{x} dt,
\ee 
where $\tilde{\B}$ represents an extension of the nonconservative matrix $\B$ in the time direction 
\[ \tilde{\B} =  (\B_1, \B_2, \mathbf{0}),\]
and $\left ( {\tilde{\F}_{ij}} + {\tilde{\mathbf{D}}_{ij}} \right ) \cdot \mathbf{\tilde n}_{ij} $ represents a well balanced space--time 
flux function augmented by the jump terms for the non-conservative product, which is explicitly designed to preserve certain equilibrium solutions 
of the PDE \textit{exactly} at the discrete level. For all the other terms the same notation of Section \ref{sec.Finite volume scheme: one step element update} holds. 

The rest of this section is organized as follows. We start giving the expression of the well balanced space--time flux function which is already enough to define a first order scheme. Then will give some hints about the extension of the well balanced techniques to second order and, we conclude presenting some numerical results.

%% SE DOVESSI METTERE ANCHE IL SECONDO ORDINE DETTAGLIATO ... we start giving the expression of the  well balanced space--time flux function, which is already enough to define a first order scheme (indeed in a first order scheme  $\q_h$ is simply the cell average $\Q_i$, which is constant in each cell and so the term involving $\tilde{\B}$ vanishes); then we explain how to modify the reconstruction procedure of Section (\ref{sec.Reconstruction}) to obtain a second order reconstruction $\q_h$ of the cell-average which is able to follow the equilibria and we discretized the remaining terms.

\subsection{Well balanced space--time flux function}

The core of the well balanced method is the design of the well balanced space--time flux function. Its final expression will be 
\be 
\begin{aligned}
	\label{eq.WBspacetimeFlux}
\left ( {\tilde{\F}_{ij}} + {\tilde{\mathbf{D}}_{ij}} \right ) \cdot \mathbf{\tilde n}_{ij}  =   \frac{1}{2} \left  ( {\tilde{\F}(\q_h^+) - \tilde{\F}(\q_h^-)} + { \mathcal{P} \left(  \q_h^+ - \q_h^- \right)  } \right ) \cdot \mathbf{\tilde n}_{ij} - \frac{1}{2} \Vn \left (\q_h^+-\q_h^- \right )
\end{aligned} 
\ee
where the term ${ \mathcal{P} \left(  \q_h^+ - \q_h^- \right)}$ represents a well balanced way to write the nonconservative products, and  $\Vn_{i+\frac{1}{2}} \left (\q_h^+-\q_h^- \right )$ is the viscosity term that we will derive slightly modifying the Osher flux.

$\mathcal{P} \left (\q_h^+-\q_h^- \right )$ and $\Vn \left (\q_h^+-\q_h^- \right )$ are defined in terms of a family of paths $\Phi(s; \q_h^-, \q_h^+)$, $s \in [0,1]$.  
In this paper, we want to choose a family of paths in such a way that stationary solutions of the shallow water equations 
\eqref{eq.ShallowWaterCylindrical} given by 
\be
\label{eq.StatSolForm}
\urho = 0, \quad \de{ \urho }{\phi} = \de{ \uphi }{\phi} = \frac{\partial \eta}{\partial \phi} = 0,  \quad \text{and} \quad \de{\eta}{r} = \frac{\uphi^2}{g r},
\ee
are preserved.

In general,  according to  the theory of \cite{DalMaso1995}, the family of paths should be a Lipschitz continuous family of functions  $\Phi(s;\q_h^-, \q_h^+), \ s \in [0,1]$,  satisfying some regularity and compatibility conditions, in particular
\be  
\label{eq.Pathproperties}
\begin{aligned}
	\Phi(0;\q_h^-, \q_h^+) = \q_h^-, \quad \Phi(1; \q_h^-, \q_h^+) = \q_h^+, \quad \Phi(s; \Q, \Q) = \Q.   
\end{aligned}
\ee
Moreover, according to \cite{Pares2006}, the numerical flux should satisfy the following properties:
\be  
\begin{aligned}   
\left ( {\tilde{\F}_{ij}}  \left(\Q, \Q \right)  + {\tilde{\mathbf{D}}_{ij}}  \left(\Q, \Q \right)  \right ) \cdot \mathbf{\tilde n}_{ij}= \mathbf{0}, \quad \forall \Q \in \Omega_{\Q},
\end{aligned}
\ee
and also for all $\q_h^-$, $\q_h^+ \ \in \Omega_{\Q}$, 
\be   
\label{PC}
\begin{aligned}
	\left  ( {\tilde{\F}(\q_h^+) - \tilde{\F}(\q_h^-)} + { \mathcal{P} \left(  \q_h^+ - \q_h^- \right)  } \right ) \cdot \mathbf{\tilde n}_{ij}  =    \int_0^1        \mathbf{A}^{\!\! \mathbf{V}}_{\mathbf{n}}(\Phi(s;\q_h^-, \q_h^+))\de{\Phi}{s}(s;\q_h^-,  \q_h^+)ds,
\end{aligned}
\ee
where
\be  
\label{eq.ExtendedJacobian}
\begin{aligned}
& \mathbf{A}^{\!\! \mathbf{V}}_{\mathbf{n}}(\Q) = \left( \! \sqrt{\tilde n_x^2 + \tilde n_y^2}\right) \left ( \, \left ( \,  \de{\mathbf{\tilde{F}}}{\Q} + \tilde{\B} \, \right ) \cdot \mathbf{{n}}   - 
(\mathbf{V} \cdot \mathbf{n}) \right ), 
\end{aligned}
\ee 
with the same notations of \eqref{eq.ALEjacobianMatrix}.
In particular $\mathcal{P}(\q_h^+-\q_h^-) $ should satisfy
\be
\label{prop-B}
\mathcal{P}(\q_h^+-\q_h^-)=  \int_0^1  \left ( \tilde{ \B}(\Phi(s;\q_h^-,\q_h^+)) \cdot \mathbf{\tilde{n}} \right ) \de{\Phi}{s}(s;\q_h^-, \q_h^+)ds. 
\ee
Note that if the standard segment path, that is
\[
\Phi(s; \q_h^-, \q_h^+) = \q_h^- + s(\q_h^+ - \q_h^-),
\]
is prescribed for all the variables, then the resulting scheme is \textit{not} well balanced. Here we propose a family of paths that is connected 
to the known equilibrium profiles for the free surface $\eta$ and the angular velocity $\uphi$, whereas for the other variables the 
segment  path is sufficient. 
Let $\Phi^E(s,\q_E^-, \q_E^+)$ be a reparametrization of a stationary solution given by \eqref{eq.StatSolForm} that connects the two equilibrium states $\q_E^-$ with $\q_E^+$, then we define $\Phi(s;\q_h^-, \q_h^+)$ as follows
\be
\label{eq:path}
\Phi(s;\q_h^-,\q_h^+)=\Phi^E(s;\q_E^-, \q_E^+)+\Phi^f(s;\q_f^-, \q_f^+), 
\ee
where $\q_f^-=\q_h^- - \q_E^-$ and $\q_f^+=\q_h^+- \q_E^+$ and
\[
\Phi^f(s;\q_f^-, \q_f^+)=  \q_f^- + s(\q_f^+  - \q_f^-).
\]
That is $\Phi^f$ is a segment path on the {\it fluctuations} with respect to a stationary solution. With this choice, it is clear that if $\q_h^-$ and $\q_h^+$ lie on the same stationary solution satisfying \eqref{eq.StatSolForm}, then $\q_f^-=\q_f^+=\mathbf{0}$ and $\Phi$, reduces to $\Phi^E$. In such situation we have that $\tilde{\F}(\q_h^+)=\tilde{\F}(\q_h^-)=\mathbf{0}$ and
\[
\mathcal{P}(\q_h^+-\q_h^-) = \int_0^1  \left (  \tilde{\B}(\Phi(s;\q_h^-,\q_h^+)) \cdot \mathbf{\tilde{n}} \right ) \de{\Phi^E}{s}(s;\q_h^-, \q_h^+)ds = \mathbf{0}.
\]
Therefore
\[
{\tilde{\F}(\q_h^+) - \tilde{\F}(\q_h^-)} + { \mathcal{P} \left(  \q_h^+ - \q_h^- \right)  } =\mathbf{0}. 
\]
For the sake of simplicity, in the following we will use the notation $\Phi(s)$ instead of $\Phi(s; \q_h^-, \q_h^+)$ when there is no confusion.

Let us now define $\mathcal{P} \left(  \q_h^+ - \q_h^- \right)$ in the general case, where $\q_h^+$ and $\q_h^-$ do not lie on a stationary solution. In such case we have that
\begin{equation}
\label{eq-P}
\mathcal{P} \left(  \q_h^+ - \q_h^- \right)= \left( b^{ij}_1, \ b^{ij}_2, \ b^{ij}_3, \ b^{ij}_4, \ b^{ij}_5 \right)^T.
\end{equation}
It is clear from the definition of $\B$ that $b^{ij}_1=b^{ij}_4=b^{ij}_5=0$.  
What is interesting is the discretization of the second term that can be rewritten as 
\be
\label{eq.bij_2_a}
\left ( g r h \de{\eta}{r} - h \uphi^2 \de{r}{r} \right ) \tilde{n}_x
= 
\left ( \  g r h \de{\eta}{r} - g r h \left [ \int \frac{\uphi^2}{r g} \, dr  \pm \int \frac{ u_{\phi, E}^2}{r g} \, dr \right]_r \  \right ) \tilde{n}_x,
\ee
where $ u_{\phi, E} $ is any known profile for the angular velocity at the equilibrium; moreover call $\zeta(r)$ a primitive of $ \frac{u_{\phi, E}^2}{r g}$, i.e.
$
\zeta(r) = \int \frac{u_{\phi, E}^2}{r g} \, d r \, .
$
In this way we obtain that  
\[
b^{ij}_2 = \int_0^1 \ \left ( \  g\Phi_{r h}(s) \de{\Phi_\eta(s)}{s}  - g \Phi_{r h}(s)  \frac{\Phi_{A}(s)}{r g} \de{\Phi_r(s)}{s} \ \right ) \tilde{n}_x \, ds,
\]
where for variables $r$ and $r h$ we can employ a standard segment path to connect the left and the right states
\[
\begin{aligned}
& \Phi_{r}(s)=\Phi_{r}(s; r^-,r^+) = r^- + s \,( r^+ - r^-), \\ 
& \Phi_{r h}(s)=\Phi_{r h}(s; (r h)^-,(r h)^+) = (r h)^- + s \,( (r h)^+ - (r h)^-).
\end{aligned} 
\]
Instead, following the idea in \eqref{eq:path} and considering the terms in \eqref{eq.bij_2_a} we define
\be
\begin{aligned}
\Phi_A(s) = \Phi_A(s; \uphi^-,\uphi^+) = \Phi_{\zeta_r}^E(s) + \frac{\Phi_{u_{\phi}}^f(s)}{r g }
\end{aligned} 
\ee 
which exploits the reparametrization of $\zeta(r)$ at the equilibrium and approximates with a segment path the fluctuations of the angular velocity
\[
\Phi_{\uphi}^f(s)=\Phi_{\uphi}^f(s; u_{\phi, f}^-,u_{\phi, f}^+) = \frac{1}{r g} \left ( u_{\phi, f}^- + s \,( u_{\phi, f}^+ - u_{\phi, f}^-) \right ).
\]
A similar approach is used for $\Phi_\eta(s)$ defined as 
\[
\Phi_\eta(s) = \Phi_\eta(s; \eta^-,\eta^+) = \Phi_{\eta}^E(s) + \Phi_{\eta}^f(s). 
\]
Taking into account that 
\[
\int_0^1 \left ( \ g \Phi_{r h}(s) \de{\Phi_{\eta}^E(s)}{s}  - g \Phi_{r h}(s)  \frac{\Phi_{\zeta_r}^E(s)}{r g} \de{\Phi_r(s)}{s} \ \right ) \tilde{n}_x \,  ds = 0,
\]
$b_2^{ij}$ could be rewritten as follows
\be
\label{eq.bij_2_d}
b^{ij}_2 = \int_0^1 \left ( \ g \Phi_{r h}(s) \de{\Phi_{\eta}^f(s)}{s}  - g \Phi_{r h}(s)  \frac{\Phi_{u_\phi}^f(s)}{r g} \de{\Phi_r(s)}{s} \ \right ) \tilde{n}_x \, ds.
\ee 
Note that  
\[ \de{\Phi_{\eta}^f(s)}{s}=\eta_f^+ - \eta_f^- = \Delta \eta_f \quad \text{and} \quad  
\de{\Phi_r(s)}{s}=r^+- r^-=\Delta r,
\]
therefore $b_2^{ij}$ reduces to
\[
b^{ij}_2 = \left (  \ g \, (r h)_{ij} \, \Delta \eta_f  - g \, (r h )_{ij}  \, \left( \frac{\uphi^2 - (u_{\phi}^E)^2}{r g} \right )_{ij} \Delta r \ \right ) \tilde{n}_x
\]
where we have employed the mid point rule to approximate the integrals and the following notation holds  $\left( \cdot \right)_{ij} = ( \cdot_i + \cdot_{j} ) / 2$.
Finally, term $b_3^{ij}$ could be approximated in the same way. Nevertheless, as this terms explicitly depends on $\urho$ and we are interested to preserve equilibria with $\urho=0$, a more simple approach could be used. Thus, $b_3^{ij}$ could be defined as
\begin{equation}
\label{eq.bij_3}
b^{ij}_3= \left ( \, \left(r h \urho \right)_{ij} \left(\uphi\right)_{ij} \Delta r \ \right ) \tilde{n}_x \, ,
\end{equation}
which vanishes when $\urho = 0$.
As pointed out in \cite{Pares2006}, a sufficient condition for a first order path-conservative scheme to be well balanced is that 
\be 
\begin{aligned}
	\left ( {\tilde{\F}_{ij}}  \left(\q_h^-, \q_h^+ \right)  + {\tilde{\mathbf{D}}_{ij}}  \left(\q_h^-, \q_h^+ \right)  \right ) \cdot \mathbf{\tilde n}_{ij}= \mathbf{0}
\end{aligned}
\ee
if $\q_h^-$ and $\q_h^+$ lie on the same stationary solution. Therefore, with the previous choice of paths, this quantity is zero  if $\Vn_{ij}(\q_h^+-\q_h^-)=\mathbf{0}$. In the next paragraph we propose a slightly modified version of the Osher flux which results to be well balanced.

\subsection{Osher Romberg viscosity matrix}
The numerical viscosity term associated to the standard path-conservative Osher scheme \cite{OsherNC} reads 
\be
\label{eq.OsherMatrixFormal}
\Vn_{i+1/2}(\q_h^+-\q_h^-) = \int_0^1  \left | \mathbf{A}^{\!\! \mathbf{V}}_{\mathbf{n}}(\Q) \left ( \Phi(s)  \right)  \right | \partial_s \Phi(s) ds, 
\ee 
with $\partial_s \Phi(s) = \partial \Phi / \partial s$. For the numerical approximation of the viscosity matrix, first we notice that it can be written as 
\be  
\Vn_{i+1/2}(\q_h^+ - \q_h^- ) = \int_0^1  \text{sign} \left( \mathbf{A}^{\!\! \mathbf{V}}_{\mathbf{n}}(\Q)  \left ( \Phi(s)  \right) \right )  \mathbf{A}^{\!\! \mathbf{V}}_{\mathbf{n}}(\Q) \left ( \Phi(s)  \right)  \partial_s \Phi(s) ds,  
\ee 
and then, we approximate the previous expression by a quadrature formula as follows:
\[
\Vn_{i+1/2}(\q_h^+-\q_h^-) =\sum_{j=1}^l \omega_j \text{sign} \left( \mathbf{A}^{\!\! \mathbf{V}}_{\mathbf{n}}(\Q) (\Phi(s_j)\right) \mathbf{A}^{\!\! \mathbf{V}}_{\mathbf{n}}(\Phi(s_j))\partial_s \Phi(s_j).   
\]
Now, we propose to approximate $\mathbf{A}^{\!\! \mathbf{V}}_{\mathbf{n}}(\Phi(s_j))\partial_s \Phi(s_j)$ by the following expression:
\[
\mathbf{A}^{\!\! \mathbf{V}}_{\mathbf{n}}(\Phi(s_j)) \partial_s \Phi(s_j) \approx \frac{\mathbf{A}^{\!\! \mathbf{V}}_{\mathbf{n},{\Phi_j}}}{2\epsilon_j}\left(\Phi(s_j+\epsilon_j)-\Phi(s_j-\epsilon_j)\right),
\]
where $\mathbf{A}^{\!\! \mathbf{V}}_{\mathbf{n},{\Phi_j}}=\mathbf{A}^{\!\! \mathbf{V}}_{\mathbf{n}}(\Phi(s_j-\epsilon_j),\Phi(s_j+\epsilon_j))$ is a Roe-matrix associated to the system (see  \cite{Pares2006} for details), that is a matrix satisfying
\begin{equation}
\mathbf{A}^{\!\! \mathbf{V}}_{\mathbf{n},{\Phi_j}}\left(\Phi(s_j+\epsilon_j)-\Phi(s_j-\epsilon_j)\right) =  \tilde{\F}(\Phi(s_j+\epsilon_j))-\tilde{\F}(\Phi(s_j-\epsilon_j)) 
  +\mathcal{P}\left(\Phi(s_j+\epsilon_j)-\Phi(s_j-\epsilon_j)\right),
\end{equation}
where $\mathcal{P}\left(\Phi(s_j+\epsilon_j)-\Phi(s_j-\epsilon_j)\right)$ is defined as in the previous section using the states $\Phi(s_j-\epsilon)$ and $\Phi(s_j+\epsilon)$.
Therefore, the viscosity term reads as follows:

\be
\label{viscosity_n_osher}
\Vn_{i+1/2}(\q_h^+-\q_h^-) = \sum_{j=1}^l \omega_j \text{sign}\left(\mathbf{A}^{\!\! \mathbf{V}}_{\mathbf{n}}(\Phi(s_j)\right) \frac{\mathcal{R}_j}{2\epsilon_j},
\ee
where 
\be
\mathcal{R}_j= \tilde{\F}(\Phi(s_j+\epsilon_j))-\tilde{\F}(\Phi(s_j-\epsilon_j))  +\mathcal{P} \left(\Phi(s_j+\epsilon_j)-\Phi(s_j-\epsilon_j)\right).
\ee

Note that if $\q_h^-$ and $\q_h^+$ lie on the same stationary solution $\Phi(s)=\Phi^E(s)$ and $\mathcal{R}_j=\mathbf{0}$, $j=1, \dots, l$, so $\Vn_{i+1/2}(\q_h^+-\q_h^-)$ vanishes. Therefore, the proposed numerical scheme is exactly well balacend for stationary solutions given by \eqref{eq.StatSolForm}.

Here, as quadrature rule, we propose the Romberg method with $l=3$ and 
\[
\begin{array}{l}
s_1=1/4, \ s_2=3/4, \ s_3=1/2, \\
\omega_1=\omega_2=2/3, \ \omega_3=-1/3 \\
\epsilon_1=\epsilon_2=1/4, \ \epsilon_3=1/2.
\end{array}
\]
With this choice, the  viscosity term $\Vn_{i+1/2}(\q_h^+ - \q_h^-)$ of the Osher-Romberg method results as follows:
\be
\begin{aligned}
	\Vn_{i+1/2}(\q_h^+ - \q_h^-) = \ 
	& \frac{4}{3} \text{sign}(\mathbf{A}^{\!\! \mathbf{V}}_{\mathbf{n}}(\Phi(1/4))) \left( \tilde{\F}(\Phi(1/2))-\tilde{\F}(\q_h^-)  +\mathcal{P}_{i+1/4}\left(\Phi(1/2)-\q_h^-\right)\right) \\
	+ & \frac{4}{3}\text{sign}(\mathbf{A}^{\!\! \mathbf{V}}_{\mathbf{n}}(\Phi(3/4))) \left(\tilde{\F}(\q_h^+)-\tilde{\F}(\Phi(1/2)) +\mathcal{P}_{i+3/4} \left( \q_h^+ - \Phi(1/2)\right) \right) \\
	-&\frac{1}{3}\text{sign}(\mathbf{A}^{\!\! \mathbf{V}}_{\mathbf{n}} (\Phi(1/2))) \left(\tilde{\F}(\q_h^+)-\tilde{\F}(\q_h^-) +\mathcal{P}_{i+1/2} \left(\q_h^+ - \q_h^-) \right)\right).
\end{aligned}
\ee
Note that the mayor drawback in the previous expression is that the complete eigenstructure of the Jacobian matrix $\mathbf{A}^{\!\! \mathbf{V}}_{\mathbf{n}}$ (\ref{eq.ExtendedJacobian}) should be computed as
\be
\text{sign}(\mathbf{A}^{\!\! \mathbf{V}}_{\mathbf{n}}) = \mbf{R}\, \text{sign}(\mbf{\Lambda})\, \mbf{R}^{-1},
\ee 
where $\mbf{\Lambda}$ is the diagonal matrix of the sign of the eigenvalues of $\mathbf{A}^{\!\! \mathbf{V}}_{\mathbf{n}}$, $\mbf{R}$ is the matrix of the right-eigenvectors and $\mbf{R}^{-1}$ its inverse. As counter part, the Osher-Romberg method is little dissipative and is stable under the standard CFL condition.

\bigskip

The extension to higher order of accuracy is also possible: for example, to obtain a second order method it suffices to slightly modify the MUSCL reconstruction procedure of Section \ref{sec.Reconstruction}. The basic idea is again to define the reconstruction operator as a combination of a smooth stationary solution together with a standard reconstruction operator to reconstruct the {\it fluctuations} with respect to the given stationary solution. 
The extension of the well balanced procedures to higher order of accuracy will be the object of another work. 
However for a complete presentation of the general framework one can refer to the work of Castro et al. in \cite{castro2008well}.

\subsection{Numerical Results}

In this section, first we want to show that the well balanced method really works in general situations and not only close to the equilibria of the system.
In this way, it will be clear that it can be applied in any context without corrupting the standard characteristics of the scheme, and it will perform better than classical schemes when near to a prescribed equilibrium.
Then, we will see that the coupling between our nonconforming techniques with the well balanced strategy allow us to study the vortex flow of Section \ref{ssec.Numresults_vortex} even for longer periods of time.

\subsubsection{Riemann problem}
\begin{figure}
	\centering
	\includegraphics[width=0.8\linewidth]{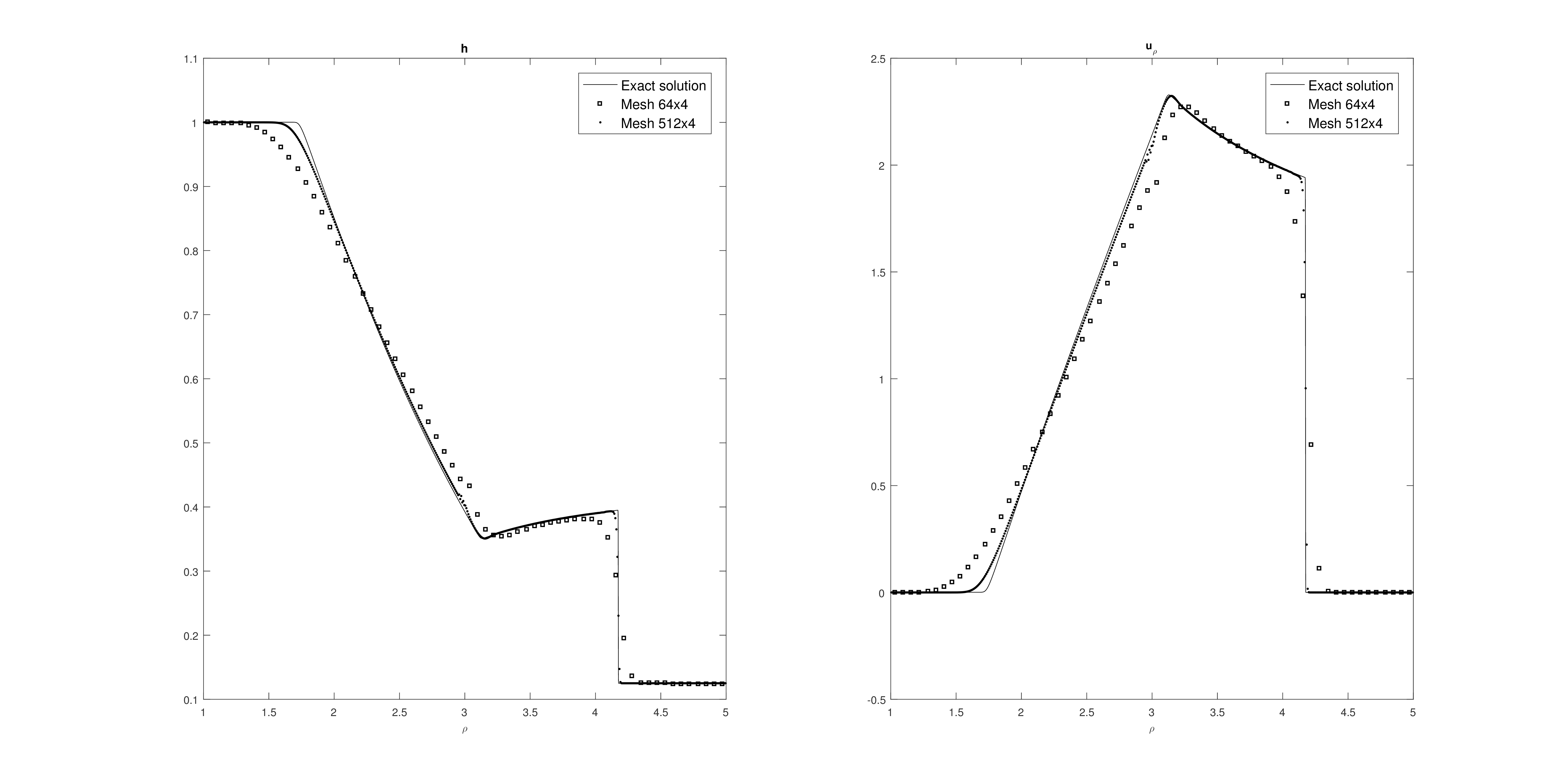}
	\caption{Comparison between the exact and the numerical solution for the Riemann problem. The numerical solution is obtained with the well balanced scheme of order one with two different meshes 
	(a coarser and a finer one). On the left we show the water level $h$ and on the right the radial velocity $u_r$ for $r \in [1,5]$ at a fixed angle $\phi = \frac{\pi}{4}$.} 
	\label{fig:wb_riemannpborder1}
\end{figure}

To show the correctness of our method we solve a classical Riemann problem with our well balanced Osher-Romberg ALE scheme. 
We consider the system of equation in \eqref{eq.ShallowWaterCylindrical}, and as computational domain $[r, \phi] =  [1,5] \times [0, 2 \pi] $. We impose the following initial conditions 
\be		
\begin{cases}
	h = 1, \  \text{if } r < r_m, \quad h = 0.125, \  \text{if } r \ge r_m, \\ 
	\urho = \uphi  = 0
\end{cases}
\ee
with $r_m = 3$. The results at the final computational time $t_f = 0.4$ are shown in Figure \ref{fig:wb_riemannpborder1}, where we report a cut along $\phi = \pi/4$.
The method, even if it is set up to preserve the smooth stationary profile described in Section \ref{ssec.Numresults_vortex}, converges properly to the reference solution of 
this problem, despite the presence of discontinuities.

\subsubsection{Steady vortex in equilibrium}

\begin{figure}
	\centering 
	\includegraphics[width=0.49\linewidth]{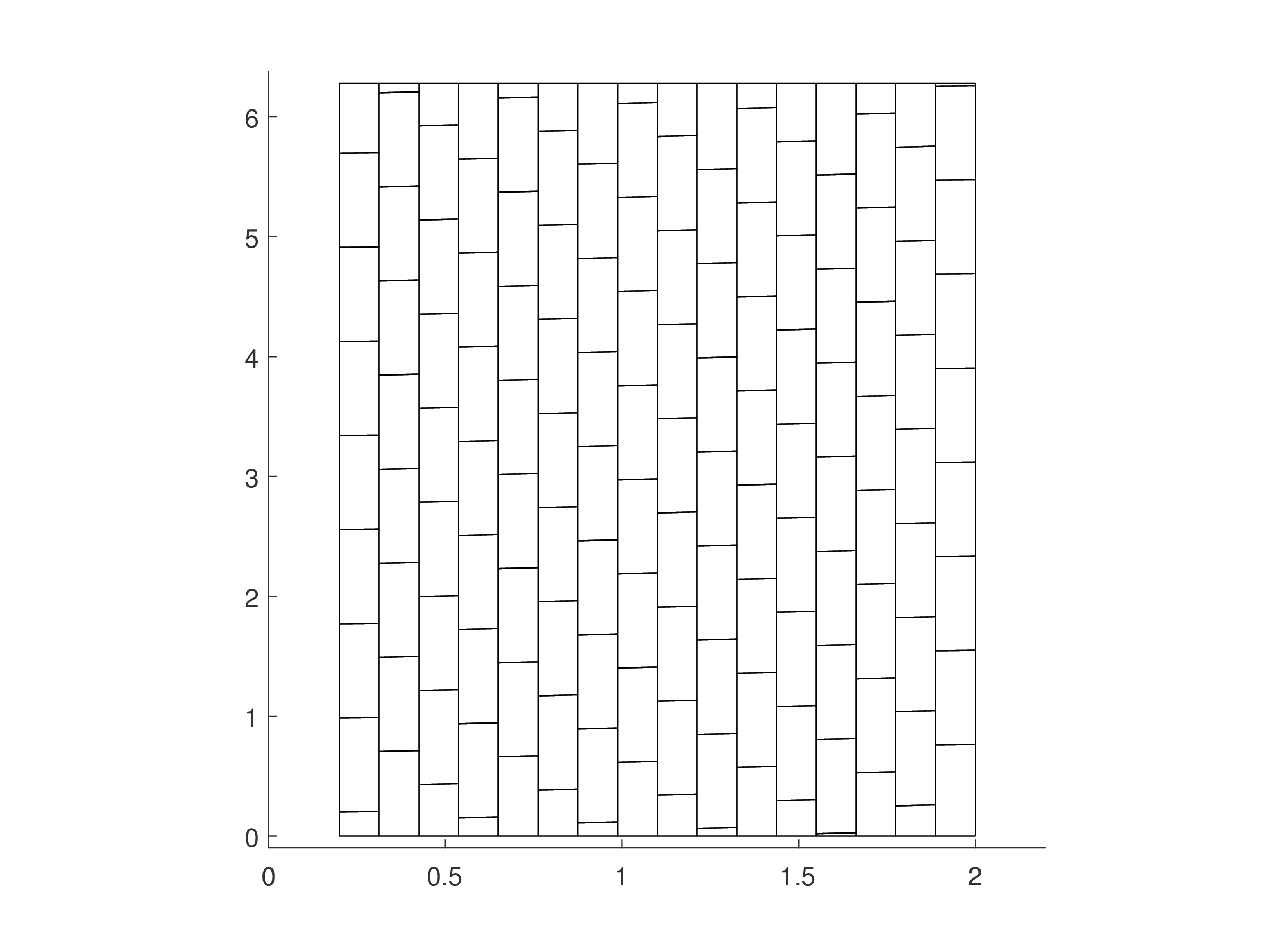}
	\includegraphics[width=0.49\linewidth]{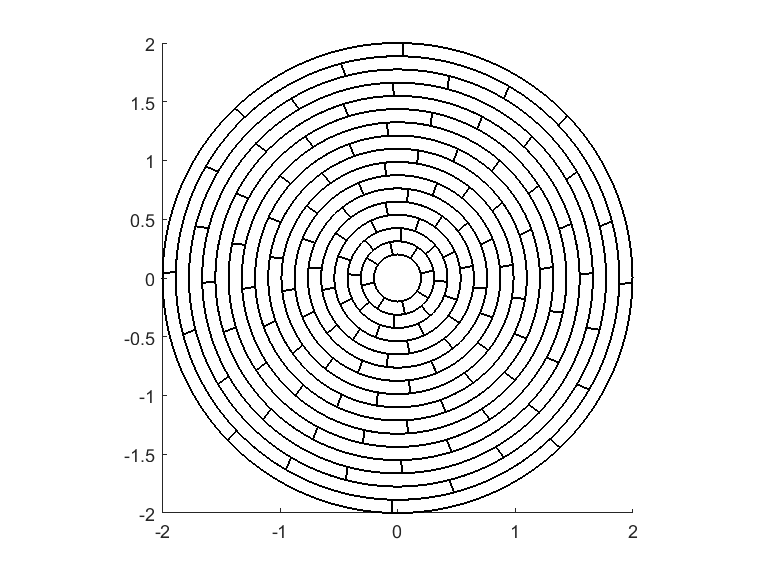}
	\caption{Stationary vortex in equilibrium obtained with well-balanced ALE schemes on moving nonconforming meshes. The mesh is shown at time $t=200$. On the left we report the grid in polar coordinates where the shear discontinuities lie  over straight lines. On the right the corresponding grid is shown in Cartesian coordinates.} 
	\label{fig.polarToCart_WB}
\end{figure}

\paragraph{Test A.} Let us consider again the test case of Section \ref{ssec.Numresults_vortex}, with the initial condition of \eqref{eq.InitalData_Equilibrium_cylindrical}. The coupling 
between our novel nonconforming ALE scheme together with the well balanced techniques gives us, even after a very long computational time, a good mesh quality 
(see Figure \ref{fig.polarToCart_WB}) and a numerical solution equal to the exact one up to machine precision (refer to Table \ref{tab.Wb_vortex}).  

Note that we have employed a mesh of squares with the constraints that interfaces lie over straight lines with constant radius. This automatically implies that each square of the mesh has two edges parallel to the $\phi-$axis: over this kind of edges the $\g$ component of the flux does not play any role, and so the method is well balanced simply because the $\f$ component of the flux is zero for stationary vortex-type solutions and \eqref{eq-P} has been proved to be discretized in the correct way. The other two edges are parallel between them, so at the equilibrium, fluxes through them cancel. 

\paragraph{Test B.} Moreover, to show that the method is able to preserve \textit{any} known stationary solution that satisfies the constraint in \eqref{eq.StatSolForm}, we have performed a similar 
test but starting from a different stationary condition 
\be 
\label{eq.InitalData_Equilibrium_2}
	 h(r, \phi,0) = \frac{r^2}{2g}, \qquad 
	 \urho(r, \phi,0) = 0, \qquad 
	 \uphi(r, \phi,0) = r, 
\ee
over the same computational domain $\Omega(r, \phi)= [0.2,2] \times [0, 2\pi]$. 
Even in this case the numerical solution remains close to the exact one up to machine precision for very long times, as also shown in Table \ref{tab.Wb_vortex}.

\begin{table} 
	\begin{center}
	\begin{tabular}{cc|cc|cc} 	
		\hline 
		\multicolumn{4}{c|}{Test A} & \multicolumn{2}{c}{Test B}   \\ 
		\hline 
		\multicolumn{2}{c|}{$ tend = 10 $} & \multicolumn{2}{c|}{ points  $16\times8$}   &  \multicolumn{2}{c}{ points  $16\times8$}    \\ 
		\hline
		points & \qquad error  \qquad \qquad & \quad time \qquad  &  error &  time & error \\ %& $\mathcal{O}2$   \\ 
		\hline
		12 $\!\times\!$ 6  & 1.42E-14 \quad & \quad 10    & 1.28E-14  &  10   & 2.11E-13  \\ 
		16 $\!\times\!$ 8  & 1.28E-14 \quad & \quad 50    & 3.74E-14  & 100   & 4.84E-13  \\ 
		24 $\!\times\!$ 12 & 3.04E-14 \quad & \quad 150   & 4.02E-14  & 150   & 3.25E-13  \\ 
		36 $\!\times\!$ 18 & 6.68E-14 \quad & \quad 200   & 4.88E-14  & 200   & 2.62E-13  \\ 
		\hline 
	\end{tabular}	
	\caption{Stationary vortex in equilibrium. Maximum error on the water level $h$ between the exact and the numerical solution obtained with the first order well balanced nonconforming ALE method. 
	         In the left column we show the error for Test A with finer and finer meshes with a fixed final time, in the central column we choose a coarse mesh and show the error for longer and longer 
					 times. In the right column, the results for Test B are shown. }
			\label{tab.Wb_vortex}
	\end{center}
\end{table}

%\begin{table} 
	%\begin{center}
		%\begin{tabular}{cc} 	
			%\hline
		    %time  & error \\
			%\hline
			%10    & 2.11E-13  \\
%%			50    & 2.10E-13  \\
			%100   & 4.84E-13  \\
			%150   & 3.25E-13  \\
			%200   & 2.62E-13  \\
			%\hline 
		%\end{tabular}	
		%\caption{Maximum error on the water level $h$ between the exact and the numerical solution obtained with the first order well balanced nonconforming ALE method over a coarse mesh of $16x8$ elements at different times.}
		%\label{tab.Wb_vortex2}
	%\end{center}
%\end{table}

\section{Conclusion} 
\label{sec.conclusion}

We have developed a robust second order direct ALE finite volume scheme on moving unstructured nonconforming meshes. The main focus was on straight slip-line interfaces, but the approach  can also easily be extended to general slide lines. In this paper, only some preliminary results for general slide lines have been shown, in order to provide a proof of concept. Further research in this direction is necessary. 
The presented results show that the method reaches its designed order of convergence and its overall accuracy. 
In particular, a high quality mesh is maintained and as a direct consequence the time step remains almost constant during the computation and the total amount of required computational effort is reduced,  despite the increased algorithmic complexity of the numerical scheme on nonconforming meshes. 
Furthermore, even with straight slip-lines, the proposed method is already able to deal with sufficiently complex situations. 
In particular, if coupled with the presented well-balanced techniques it can be also considered for practical applications, for example in the context of the compressible Euler equations of gasdynamics  with gravity, which are highly relevant in computational astrophysics, e.g. for the simulation of accretion discs and rotating gas clouds around compact objects like stars, neutron stars or black holes. 

%The aim of our future research will be also to remove the constraint of straight lines for the interface, adopting a distinction between master and slave nodes and improving the algorithms for 
%the update of the connectivity tables. 

In future research we plan to extend the presented method to better than second order of accuracy by extending the ADER-WENO and ADER-DG ALE schemes 
\cite{LagrangeMHD,Lagrange2D,Lagrange3D,LagrangeDG} to moving nonconforming unstructured meshes. Further research may also concern the incorporation of 
time accurate local time-steping (LTS) \cite{ALELTS1D,ALELTS2D,LTSTransport} into our algorithm.

%=============================================================================
%==========    A C K N O W L E D G M E N T S
\section*{Dedication} 

The new numerical method introduced in this paper is dedicated to Prof. \textbf{Eleuterio Francisco Toro} at the occasion of his 70th birthday and in honor of his groundbreaking 
scientific contributions to the field of shock capturing methods for computational fluid dynamics.  

\section*{Acknowledgments} 
The research presented in this paper has been partially financed by the European Research Council (ERC) under the 
European Union's Seventh Framework Programme (FP7/2007-2013) with the research project \textit{STiMulUs}, 
ERC Grant agreement no. 278267. This research has been also supported by the Spanish Government and FEDER
through the research project MTM2015-70490-C2-1-R and the Andalusian Government research projects
P11-FQM-8179 and P11-RNM-7069.
Moreover this project has received funding from the European Union's Horizon 2020 research and innovation Programme
under the Marie Sklodowska-Curie grant agreement no. 642768.

%=============================================================================
%==========  B I B L I O G R A P H Y
%\input{Bib}
\bibliographystyle{plain}
\bibliography{references}
%=============================================================================

%=============================================================================
%==========  A P P E N D I X
%\input{Appendix}
%=============================================================================

\end{document}